\documentclass[12pt]{article}

\usepackage{amsfonts}
\usepackage{epsfig}
\usepackage[percent]{overpic}
\usepackage{amsmath}

\usepackage{cite}

\usepackage[colorlinks]{hyperref}
\usepackage{amssymb}
\usepackage{graphics}
\usepackage{verbatim}
\usepackage{psfrag}
\usepackage{color}
\usepackage{amsthm}
\usepackage[all]{xy}
\usepackage{makeidx}
\usepackage{enumerate}
\usepackage[letterpaper,margin=0.7in]{geometry}
\usepackage{mathtools}

\newcommand{\eex}{\hbox{}\hfill\rule{0.8ex}{0.8ex}}

\usepackage[percent]{overpic}
\usepackage{caption}
\usepackage{float}
\usepackage{placeins}

\usepackage{multicol}
\usepackage{subcaption}
\usepackage{float}

\let\oldbibliography\thebibliography
\renewcommand{\thebibliography}[1]{%
  \oldbibliography{#1}%
  \setlength{\itemsep}{0.5mm}%
}

\theoremstyle{plain}
\newtheorem{thm}{Theorem}[section]

\newtheorem{lemma}[thm]{Lemma}
\newtheorem{proposition}[thm]{Proposition}
\newtheorem{theorem}[thm]{Theorem}
\newtheorem{definition}[thm]{Definition}
\newtheorem{ex}[thm]{Example}

\newtheorem*{remark}{\textit{Remark}}

\def\R{\mathbb{R}}
\def\C{\mathbb{C}}
\def\N{\mathbb{N}}

\def\P{\mathbb{P}}
\def\E{\mathbb{E}}

\def\cE{\mathcal{E}}

\def\cK{\mathcal{K}}
\def\cL{\mathcal{L}}

\def\cO{\mathcal{O}}

\def\cX{\mathcal{X}}

\def\txte{{\textnormal{e}}}
\def\txti{{\textnormal{i}}}

\def\txtD{{\textnormal{D}}}

\def\Re{{\textnormal{Re}}}
\def\Im{{\textnormal{Im}}}

\def\I{\infty}

\newcommand{\be}{\begin{equation}}
\newcommand{\ee}{\end{equation}}
\newcommand{\benn}{\begin{equation*}}
\newcommand{\eenn}{\end{equation*}}
\newcommand{\bea}{\begin{eqnarray}}
\newcommand{\eea}{\end{eqnarray}}
\newcommand{\beann}{\begin{eqnarray*}}
\newcommand{\eeann}{\end{eqnarray*}}

\newcommand{\myendex}{$\blacklozenge$\end{ex}}
\newcommand{\myendexerc}{$\lozenge$\end{exerc}}
\newcommand{\myendpexerc}{$\lozenge$\end{pexerc}}

\usepackage[inline]{enumitem}
\newcommand{\caputo}[2]{\mathfrak{\tilde{D}}^{#1}_{#2}}
\newcommand{\RL}[2]{\mathfrak{D}^{#1}_{#2}}

\newcommand{\der}[2]{D^{#1}_{#2}}

\newcommand{\dt}{{\Delta t}}
\newcommand{\bi}{{\rm i}}

\newcommand{\SchwartzTF}{\mathcal{S}} 

\newcommand{\dd}[1][x]{\,\operatorname{d}\!#1}
\newcommand{\diff}[2]{\partial_{#2} {#1}}
\newcommand{\difff}[3]{\partial^{#3}_{#2} {#1}}

\newcommand{\integral}[3]{\int_{#1} \, {#2} \, \dd[#3]\,}
\newcommand{\integrall}[4]{\int_{#1}^{#2} \, {#3} \, \dd[#4]\,}
\newcommand{\abs}[1]{|#1|}

\newcommand{\Levy}{\cL_{Levy}}  
\newcommand{\Riesz}[1][\gamma]{\operatorname{D}_0^{#1}}  
\newcommand{\RieszFeller}[1][\theta]{\operatorname{D}_{#1}^\gamma}  
\newcommand{\Fourier}{\mathcal{F}}

\newcommand{\profile}{\phi}
\newcommand{\speed}{s}
\newcommand{\um}{u_-}
\newcommand{\up}{u_+}
\newcommand{\upm}{u_\pm}
\newcommand{\Rd}{\R^d}

\newcommand{\xx}[1]{\,\text{ #1 }\,}
\newcommand{\Xx}[1]{\quad\text{ #1 }\,}

\usepackage{todonotes}

\newcommand{\blue}[1]{{\color{blue} {#1}}}

\usepackage[normalem]{ulem} 

\def \diag {\mathrm{diag}}

\def \sgn {\mathrm{sgn}}
\def \tr {\mathrm{tr}}
\def \dif {\mathrm{d}}

\def \u {{U}}

\def \w {{W}}
\def \A {{A}}

\def \F {{F}}
\def \I {{I}}
\def \P {{P}}
\def \M {{M}}

\def \fD {\mathfrak{D}}

\def \fL {\mathfrak{L}}
\def \BC {\mathfrak{B}}

\def \Lspace {\mathnormal{L}}
\def \Hspace {\mathnormal{H}}
\def \ACspace {\mathnormal{AC}}

\def \reg {{\rm reg}}
\def \ss {{\rm add}} 
\def \ca {{\rm lin}}

\def \AA{A$_{\rm add}$} 
\def \AL{A$_{\rm lin}$} 

\def\blue{\color{blue}}

\author{
Franz Achleitner
\thanks{TU Wien, Institute for Analysis and Scientific Computing, 
 1040 Vienna, Austria, \texttt{franz.achleitner@tuwien.ac.at}}~~~
Goro Akagi
\thanks{Tohoku University, Mathematical Institute and Graduate School of Science, 980-8578 Sendai, Japan, 
\texttt{goro.akagi@tohoku.ac.jp}}~~~
Christian Kuehn
\thanks{Technical University of Munich, Department of Mathematics, Research Unit ``Multiscale and Stochastic Dynamics'', 85748 Garching b.~M\"unchen, Germany, \texttt{ckuehn@ma.tum.de}}
\thanks{Technical University of Munich, Munich Data Science Institute (MDSI), 85748 Garching b.~M\"unchen, Germany}
\thanks{External Faculty, Complexity Science Hub Vienna, 1080 Vienna, Austria}~~~\\
Jens Markus Melenk
\thanks{TU Wien, Institute for Analysis and Scientific Computing, 
 1040 Vienna, Austria, \texttt{melenk@tuwien.ac.at}}~~~
Jens D.M. Rademacher
\thanks{Universit\"at Hamburg, MIN faculty, Department of Mathematics, Applied Dynamical Systems, 20146 Hamburg, Germany, \texttt{jens.rademacher@uni-hamburg.de}}~~
Cinzia Soresina
\thanks{Universit\"at Graz, Department of Mathematics and Scientific Computing, 8010 Graz, Austria, \texttt{cinzia.soresina@uni-graz.at}}~~
Jichen Yang
\thanks{Harbin Engineering University, College of Mathematical Sciences, 150001 Harbin, China, \texttt{jichen.yang@hrbeu.edu.cn}}
}

\begin{document}

\title{Fractional Dissipative PDEs}

\maketitle

\abstract{In this chapter we provide an introduction to fractional dissipative partial differential equations (PDEs) with a focus on trying to understand their dynamics. The class of PDEs we focus on are reaction-diffusion equations but we also provide an outlook on closely related classes of PDEs. To simplify the exposition, we only discuss the cases of fractional time derivatives and fractional space derivatives in the PDE separately. As our main tools, we describe analytical as well as numerical methods, which are generically necessary to study nonlinear dynamics. We start with the analytical study of steady states and local linear stability for fractional time derivatives. Then we extend this view to a global perspective and consider time-fractional PDEs and gradient flows. Next, we continue to steady states, linear stability analysis and bifurcations for space-fractional PDEs. As a final analytical consideration we discuss existence and stability of traveling waves for space-fractional PDEs. In the last parts, we provide numerical discretization schemes for fractional (dissipative) PDEs and we utilize these techniques within numerical continuation in applied examples of fractional reaction-diffusion PDEs. We conclude with a brief summary and outlook on open questions in the field.}

\section{Introduction}
\label{sec:intro}

Reaction-diffusion partial differential equations (PDEs) are a crucial modeling tool appearing in all the natural sciences, engineering and beyond~\cite{Britton2,Grindrod,Murray2}. In particular, reaction-diffusion PDEs are relevant at very small length scales up to very large macroscopic length scales. The existence, uniqueness and regularity theory for large classes of reaction-diffusion PDEs is quite classical by now and many aspects are well-understood~\cite{Amann4,Evans,Lunardi,Rothe,Smoller}. Yet, asking for a more precise description of the qualitative and quantitative dynamics of quite small classes of reaction-diffusion PDEs, a lot of questions appear that have led to very active research areas. There are few absolutely complete and comprehensive results about the dynamics of reaction-diffusion PDEs. This is expected as such a dynamical classification does not even exist for nonlinear ODEs formed by the reaction terms. Yet, for the case with enough dissipativity, there is a well-established strategy that can be applied. First, one uses dissipativity to establish boundedness of solutions, often obtaining the existence of a (finite-dimensional) attractor as a consequence of global dissipation~\cite{Robinson1,Temam}. Unfortunately, just having  the existence of an attractor yields relatively little information on the precise dynamics, so one then proceeds to study finer information via invariant sets and patterns~\cite{Henry,KuehnBook1,SchneiderUecker}. A first step in this direction is the analysis of steady states. Often one starts with homogeneous steady states and then proceeds to more complex steady states using complementary methods, e.g., bifurcation theory~\cite{ChowHale,Kielhoefer,Kuznetsov} and variational methods~\cite{Zeidler3}. Once steady states are found, one continues to investigate their local linear stability via classical methods such as spectral theory~\cite{Kato}. A next step is the investigation of non-stationary patterns such as traveling waves~\cite{KuehnBook1}, which are effectively invariant steady sets in a moving frame. Again, one then studies the stability of traveling waves~\cite{KapitulaPromislow,Sandstede1}. To cross-validate and globally extend dynamical results, it is often imperative to combine analytical approaches with numerical techniques, which are based upon spatial discretization combined with temporal simulation or numerical continuation methods~\cite{DankowiczSchilder,Doedel_AUTO2007,uecker_pde2path_2014}. 

To illustrate this classical paradigm, let us consider a typical example given by the PDE
\be
\label{eq:AC}
\partial_t u = \partial_x^2 u + u(p_1-u)(u-p_2)=:\partial_x^2u+f(u;p),\quad u_0(x)=u(0,x),
\ee
where $u=u(t,x)$, $x\in\Omega \subseteq \R$, $t\geq 0$, and $p=(p_1,p_2)\in\R^2$ are parameters. The reaction-diffusion PDE~\eqref{eq:AC} has many different names in different fields: Fisher equation with Allee effect~\cite{Murray2}, Allen-Cahn equation~\cite{AllenCahn}, Chafee-Infante equation~\cite{ChafeeInfante}, Nagumo equation~\cite{Nagumo}, real Ginzburg-Landau equation~\cite{CrossHohenberg}, Schl\"ogl model~\cite{Schloegl}, and (overdamped) $\Phi^4$ model~\cite{Hairer3}. This clearly demonstrates the importance of \eqref{eq:AC} as a key example, or even as a normal form~\cite{KirrmannSchneiderMielke}, for many dynamical phenomena. One easily checks that~\eqref{eq:AC} is dissipative due to the $-u^3$ term in the nonlinearity~\cite{Robinson1,Temam}. For the case of a bounded domain $\Omega=(a,b)$, $a<b$, one obtains a finite-dimensional attractor~\cite{FiedlerRocha,FiedlerScheel,Temam,Robinson1}. For concreteness, let us consider homogeneous Neumann boundary data
\be
\label{eq:ACNeumann}
\frac{\partial u}{\partial x}(a)=0=\frac{\partial u}{\partial x}(b).
\ee 
In this case, we immediately spot the three homogeneous equilibria $u_{*,0}(x)\equiv 0$, $u_{*,1}=p_1$, $u_{*,2}=p_2$. The standard configuration for the Allen-Cahn equation is to fix $p_2=-p_1$ and set $p_1^2=:p_{AC}$ so that the nonlinearity becomes $p_{AC}u-u^3$. Even in this simple case, the homogeneous state $u_{*,0}$ undergoes bifurcations as $p_{AC}$ is increased generating non-homogenous/patterned steady states in the process~\cite{Henry,Kielhoefer,KuehnBook1}. Note that linear stability of any steady state on a $\Omega=(a,b)$ reduces to the study of the eigenvalue problem 
\benn
(\txtD_uf)(u_{*,\cdot};p)U=\lambda U,\qquad U=U(x),
\eenn
where $\txtD_u$ denotes the Fr\'echet derivative of $f$ on a suitable function space $\cX$, e.g., one could consider $\cX=H^1(\Omega)=W^{1,2}(\Omega)$. Computation and visualization of the non-homogeneous steady states is possible locally near bifurcation points. Globally, the
situation is already very difficult and it helps to employ numerical tools~\cite{DankowiczSchilder,Doedel_AUTO2007,uecker_pde2path_2014}. A natural step is to spatially discretize the steady state problem, e.g., using finite elements, finite  differences or a spectral Galerkin method. Then one can apply path-following techniques to generate global bifurcation diagrams~\cite{DankowiczSchilder,Doedel_AUTO2007,uecker_pde2path_2014}.   
 
For the case of an unbounded domain $\Omega=\R$ and the standard Nagumo parameter choices $p_1=1$ and $p_2=p_{Nag}\in(0,1)$, it is well-known that the PDE~\eqref{eq:AC} has traveling wave solutions~\cite{Sandstede2,KapitulaPromislow,KuehnBook1} of the form 
\be
\label{eq:AnsatzTW}
u(t,x)=u(x-st)=u(\zeta) 
\ee
for a wave speed $s=s(p_{Nag})$. In fact, the ansatz~\eqref{eq:AnsatzTW} just results in a second-order ODE. One can prove that this ODE has a parametric family of heteroclinic orbits corresponding to traveling front solutions and a parametric family of homoclinic orbits corresponding to traveling pulse solutions. One can then prove using linearization and Sturm-Liouville theory that the traveling front is linearly (and in fact also nonlinearly) stable, while the traveling pulse is already linearly unstable~\cite{Sandstede2,KapitulaPromislow,KuehnBook1}. To visualize the traveling wave dynamics, one option is to just combine spatial discretization on a very 
large domain $(a,b)$ with $a\ll -1$ and $b\gg1$ with a time discretization. Alternatively, one could freeze the wave by using a suitable reference frame~\cite{BeynThuemmler}.

Although the theory of~\eqref{eq:AC} is very detailed, an important question for applications is its robustness. This not only means proving stability results in phase space but also refers to structural stability of the model itself, i.e., if the PDE class is changed and/or extended, do dynamical phenomena persist or do new effects appear? Before answering this question, we have to consider the input from applications. Which extensions to reaction-diffusion PDEs do occur in practice? One answer is provided by going back to the microscopic derivation of diffusion. In particular, the Laplacian $\Delta$ as well as the associated heat operator $\partial_t-\Delta$ can arise from taking a macroscopic limit of a random walk described by a stochastic process $z(t)$. If the stochastic process satisfies a power law for the mean-squared displacement
\benn
\E[z(t)^2]\sim 2\kappa_\alpha t^\alpha,
\eenn
then one easily calculates that a standard random walk has the exponent $\alpha=1$. The case $\alpha=1$ could also be called ``classical'' diffusion. Classical diffusion is a Markovian and a Gaussian/Wiener process. Note that $\kappa_1=:\sigma$ becomes the diffusion constant of the mean-field limit diffusion equation
\benn
\partial_t u = \sigma \Delta u
\eenn 
for the probability density $u(t,x)$ of the random walker to be at position $x$ at time $t$. Yet, if $\alpha>1$, then the behavior is super-diffusive (or ``ballistic''), while $\alpha\in(0,1)$ gives subdiffusive dynamics. The cases $\alpha\neq 1$ are also called anomalous diffusion~\cite{MetzlerKlafter} and lead to different macroscopic operators beyond the standard heat equation; we will define these fractional differential operators, as well as several additional ones, in Section~\ref{sec:def}. Once all these additional fractional operators arising from various applications have been collected, we can try to proceed to study dissipative fractional PDEs along the same lines as the classical example~\eqref{eq:AC} presented above. 

The remaining part of this review is structured as follows. In Section~\ref{sec:def} we fix the notation we use for this chapter. In Section~\ref{sec:spectral1} we study the extension of steady states and linear stability with a focus on fractional time derivatives. In Section~\ref{sec:gradient} we take a global perspective and consider time-fractional reaction-diffusion equations and their associated gradient flows. In Section~\ref{sec:spectral2} we briefly recap linear stability analysis and bifurcations for space-fractional derivatives. In Section~\ref{sec:tw}, we survey results for traveling waves for space-fractional PDEs. In Section~\ref{sec:discrete}, we give an overview for numerical discretization schemes for fractional (dissipative) PDEs, while in Section~\ref{sec:continuation}, we show how these techniques can be embedded into numerical continuation in practical examples of fractional reaction-diffusion PDEs arising in applications. In Section~\ref{sec:outlook}, we provide a brief summary and outlook towards interesting open problems. 

We stress that our review provides only a partial snapshot of the field of dynamics of fractional dissipative PDEs. It is probably already impossible to cover the field comprehensively as it has grown so quickly in recent years. Yet we hope that this chapter makes it much easier to find an entry point to the results and techniques in the field, even beyond the results we mention here. Furthermore, we believe that the snapshot of results discussed here clearly indicates that the general approach to study dissipative reaction-diffusion PDEs can be lifted from classical to fractional derivatives. Yet, the techniques and phenomena are often substantially different, which makes it an exciting area to explore. 

\section{Definitions and Notation}
\label{sec:def}

Let $\Omega$ now denote a domain in $\R^d$. Consider $x\in\Omega\subseteq \R^d$, $t\geq0$ and write  $u=u(t,x)\in\R^N$ for the (unknown) solution of the PDEs we study here; often we will just take $d=1$ and $N=1$ as the main basic case. $U$ will be used to highlight a linearized problem $U=U(t,x)\in\R^N$. $\Delta$ denotes the classical  Laplacian, $\partial_x^2$ the classical Laplacian in one dimension, and $\partial_t$ the classical partial time derivative.

\subsection{Time-Fractional Derivative(s)}
\label{ssec:deftimefrac}

First, we define the Riemann-Liouville fractional integral and derivative for the subdiffusion range $\alpha\in(0,1)$. Let $f(t) \in \Lspace^1(0,T)$ for any $T>0$. The fractional Riemann-Liouville integral of order $\alpha$ is defined as
	\begin{equation}
	\label{e:FracIntegral}
		(\fD_{0,t}^{-\alpha}f)(t) :=  (k_\alpha \ast f)(t) = \frac{1}{\Gamma(\alpha)}\int_0^t (t-s)^{\alpha-1}f(s)\, \dif s, \quad \alpha\in(0,1),
	\end{equation}
where the (Laplace) convolution kernel is $k_{\alpha}(t):=t^{\alpha-1}/\Gamma(\alpha)$, and $\Gamma(\alpha)$ is the Gamma function. We can also generalize the time dependence of the integrals yielding the natural definitions
\begin{align*}
\RL{-\alpha}{t_0,t} f(t):= \frac{1}{\Gamma(\alpha)} \int_{t_0}^t (t-s)^{\alpha-1} f(s)\,\dif s, \qquad 
\RL{-\alpha}{t,T} f(t):= \frac{1}{\Gamma(\alpha)} \int_{t}^T (s-t)^{\alpha-1} f(s)\,\dif s.
\end{align*}
Notably, for $\alpha=1$ we get $(\fD_{0,t}^{-\alpha}f)(t) = \int_0^tf(s)\dif s$. For $f:[0,T]\to \R$ the Riemann-Liouville fractional derivative of order $\alpha$ is defined as
\begin{equation}
\label{e:RLDerivative}
(\fD_{0,t}^{\alpha}f)(t) := \frac{\dif}{\dif t}(k_{1-\alpha}\ast f)(t) = \frac{1}{\Gamma(1-\alpha)}\frac{\dif}{\dif t} \int_0^t (t-s)^{-\alpha}f(s)\, \dif s, \quad \alpha\in(0,1),
\end{equation}
with the obvious generalizations to other endpoints $t_0$, $T$, which is then reflected in the notation by adding suitable indices. Notably, the Riemann-Liouville fractional derivative is nonzero on constants, 
\begin{equation*}
	(\fD_{0,t}^{\alpha}1)(t) = k_{1-\alpha}(t),
\end{equation*}
which tends to zero as $t^{-\alpha}$ for $t\to\infty$ and is unbounded for $t\to 0$. At $\alpha=0$, $\fD_{0,t}^{\alpha}$ is the identity operator, i.e., $(\fD_{0,t}^{\alpha}f)(t)=f(t)$, while for $\alpha=1$, $(\fD_{0,t}^{\alpha}f)(t)$ formally yields $(\dif/\dif t)\left(\delta\ast f\right)(t)=f'(t)$ with Dirac delta distribution $\delta(t)$ and $':=\dif/\dif t$. A simple sufficient condition for the existence of the Riemann-Liouville fractional derivative is as follows.

\begin{lemma}[{\cite[Lemma 2.2]{samko93}}]
\label{l:RLderivativeexist}
	Let $f(t) \in \ACspace([0,T])$ and $\alpha \in (0,1)$. Then $(\fD_{0,t}^{\alpha} f)(t)$ exists almost everywhere. Moreover $(\fD_{0,t}^{\alpha}f)(t) \in \Lspace^p(0,T)$, $1\leq p <1/\alpha$, and
	\begin{align*}
		(\fD_{0,t}^{\alpha} f)(t) = (k_{1-\alpha}\ast f')(t) + k_{1-\alpha}(t)f(0).
	\end{align*}
\end{lemma}

Here $\ACspace([0,T])$ is the space of functions $f$ which are absolutely continuous on $[0,T]$, i.e., $f(t) = c + \int_0^t g(s) \dif s$ for some $g \in \Lspace^1(0,T)$ and constant $c$. This lemma motivates the definition of the (Djrbashian-)Caputo fractional derivative. For $f\in \ACspace([0,T])$ the \emph{(Djrbashian-)Caputo fractional derivative} of order $\alpha$ is defined as
	\begin{equation}\label{e:Caputo}
		(\tilde\fD_{0,t}^{\alpha}f)(t) := (k_{1-\alpha} \ast f')(t) = \frac{1}{\Gamma(1-\alpha)} \int_0^t (t-s)^{-\alpha}f'(s)\, \dif s, \quad \alpha\in(0,1).
	\end{equation}
Hence, for sufficiently smooth $f:[0,T]\to\R$, Lemma \ref{l:RLderivativeexist} gives a well-known relation between Riemann-Liouville and Caputo derivatives
\begin{align*}
\tilde\fD_{0,t}^\alpha f(t) = \fD_{0,t}^\alpha(f(t) - f(0)).
\end{align*}
One can extend these definitions also beyond the case $\alpha\in(0,1)$ to any $\alpha>0$. Define $m := \lceil \alpha\rceil \in {\mathbb N}_0$ for $\alpha \not\in {\mathbb N}_0$ and $m:= \alpha+1$ for $\alpha \in {\mathbb N}_0$.
Then, the Riemann-Liouville derivatives are given by
\begin{align*}
\RL{\alpha}{0,t} f(t)&:= \frac{\dif^m}{\dif t^m} \RL{-(m-\alpha)}{0,t} f(t),  &
\RL{\alpha}{t,T} f(t)&:= (-1)^m \frac{\dif^m}{\dif t^m} \RL{-(m-\alpha)}{t,T} f(t).
\end{align*}
The left and right Caputo derivatives are given by
\begin{align*}
\caputo{\alpha}{0,t} f(t)&:= \RL{-(m-\alpha)}{0,t} f^{(m)}(t), &
\caputo{\alpha}{t,T} f(t)&:= (-1)^m \RL{-(m-\alpha)}{t,T} f^{(m)}(t),
\end{align*}
where one requires sufficient regularity, in terms of a generalization of Lemma \ref{l:RLderivativeexist} to ensure existence of the derivatives. The Riemann-Liouville and the Caputo fractional derivatives are related to each other: For sufficiently smooth $f$, one has
\begin{align*}
\RL{\alpha}{0,t} f(t) & = \caputo{\alpha}{0,t} u + \sum_{k=0}^{m-1} \frac{f^{(k)}(0) (t-0)^{k-\alpha}}{\Gamma(k+1-\alpha)}.
\end{align*}

\subsection{Space-Fractional Derivative(s)}
\label{ssec:defspacefrac}

We start with the case of the spatial domain $\Omega=\R$. In experiments asymmetric superdiffusive behavior has been observed in many applications~\cite{BeScMeWh01,DelCastilloNegrete:1998,dCNCaLy04,Hanert+etal:2011,MeJeChBa14} (see also the references in~\cite{KlaRadSok08,ScMeBa09} and~\cite[Part 3]{delia-etal20}). Using continuous time random walks (CTRW) with an asymmetric distribution of jump lengths allows one to derive evolution equations with asymmetric L\'evy operators~\cite[\S 3.3.3]{Mendez+etal:2010}. In this context, a very general class of space-fractional derivatives is derived microscopically from L\'evy operators. A subclass of L\'evy operators are the Riesz--Feller operators, which are the infinitesimal generators of stable L\'evy processes; see~\cite[\S 14]{Sato:1999} and \cite{DiROSeVa22}. In the  one-dimensional setting, Riesz--Feller operators can be characterized by two parameters: the index of stability~$\gamma$ and another parameter $\theta$ measuring the asymmetry. Riesz--Feller operators can be defined as Fourier multiplier operators 
 \[ \Fourier[\RieszFeller u](\xi) = \psi^\gamma_\theta (\xi)\ \Fourier [u] (\xi), \qquad \xi\in\R, \] 
 with $\Fourier$ denoting the Fourier transform and symbol
 \[ \psi^\gamma_\theta(\xi) = -|\xi|^\gamma \exp\left[\txti (\sgn(\xi)) \theta \tfrac{\pi}{2}\right]
     \Xx{for} 0<\gamma\leq 2 \xx{and} |\theta| \leq \min\{\gamma, 2-\gamma\}.
 \]
We will require enough regularity for $u=u(x)$ with $x\in\R$ so that we can define its Fourier transform. Special cases for Riesz-Feller operators are the second order derivative $\Riesz[2] = \difff{}{x}{2}$, the (negative) fractional Laplacians $\Riesz = -(-\partial^2_x)^{\gamma/2}$ for $0<\gamma\leq 2$ and $\theta=0$ (which are the only symmetric Riesz--Feller operators), and Weyl fractional derivatives $\RieszFeller[2-\gamma]$ for $0<\gamma<2$ and $\theta=2-\gamma$ (which are extremely asymmetric Riesz--Feller operators). 

For the higher-dimensional setting $\Omega=\R^d$, we will only need the (positive) fractional Laplacian $(-\Delta_x)^{\gamma/2}$, $\gamma\in(0,2)$. In full space fractional Laplacians $(-\Delta_x)^{\gamma/2}$ can be defined in various equivalent ways, see, e.g.,~\cite{kwasnicki17}: 
For example, again as Fourier multiplier operators
\[ 
 \Fourier[(-\Delta_x)^{\gamma/2} u](\xi) 
 = |\xi|^\gamma \ \Fourier [u] (\xi), \qquad \xi\in\R^d, 
\] 
with symbol~$|\xi|^\gamma$. Or, we can define them as singular integral operators in two different ways: 
\be 
\label{fractionalLaplacian_multiD:singular_integral}
 (-\Delta_x)^{\gamma/2} u (x)
 := c_{d,\gamma}  \integral{\R^d}{\frac{-u(x+y)-u(x-y)+2u(x)}{2 |y|^{d+\gamma}}}{y} = c_{d,\gamma} P.V. \int_{\R^d} \frac{u(x) - u(y)}{|x - y|^{d + \gamma}}\,dy, 
\ee
where P.V.\ indicates that a principal value integral is considered, and the normalizing constant
\be \label{c:d_alpha}
 c_{d,\gamma}
:= \frac{2^\gamma \Gamma\big(\tfrac{d+\gamma}2\big)}{\pi^{d/2} \big|\Gamma\big(-\tfrac{\gamma}2\big)\big|}
\ee
is such that the singular integral operators have the Fourier symbol~$|\xi|^{\gamma}$. On bounded domains, the situation becomes more complicated and one has to carefully distinguish different definitions for the fractional Laplacian, e.g., the spectral fractional Laplacian (SFL) and the integral fractional Laplacian (IFL). These operators coincide on ${\mathbb R}^d$, but they differ when restricted to bounded domains $\Omega$. To define the Dirichlet spectral fractional Laplacian let $\Omega\subset {\mathbb R}^d$ be a bounded Lipschitz domain and let $(\phi_i, \lambda_i)_{i \in {\mathbb N}} \subset H^1_0(\Omega) \setminus \{0\} \times {\mathbb R}$ be the  eigenpairs of
\begin{equation}
\label{eq:mel-EWP}
\Delta \phi  = \lambda \phi \quad \mbox{ in $\Omega$},  \quad \phi|_{\partial\Omega} = 0
\end{equation}
with the eigenvalues $\lambda_i \rightarrow -\infty$ sorted in descending order and normalized such that $\|\phi_i\|_{L^2(\Omega)} = 1$, $i \in {\mathbb N}$. Note that $(\phi_i)_{i\in {\mathbb N}}$ forms an orthonormal basis of $L^2(\Omega)$. The (positive) Dirichlet spectral fractional Laplacian $(-\Delta)^{\gamma/2}$ is defined as the operator $u \mapsto (-\Delta)^{\gamma/2}  u(x):= \sum_{i} (-\lambda_i)^{\gamma/2} (u,\phi_i)_{L^2(\Omega)}\phi_i(x)$. Note that $\gamma = 0$ reproduces the identity operator and $\gamma = 2$ the standard Laplacian. The operator is an isomorphism $\widetilde H^{\gamma/2}(\Omega) \rightarrow H^{-\gamma/2}(\Omega):= (\widetilde{H}^{\gamma/2}(\Omega))^\prime$, where for $\beta \in [0,1]$ the fractional Sobolev space
$\widetilde{H}^\beta(\Omega)$ is the closure of $C^\infty_0(\Omega)$ under the norm
\begin{equation}
\|u\|_{\widetilde{H}^\beta(\Omega)}:= 
\left(\sum_{i} (-\lambda_i)^\beta |(u,\phi_i)_{L^2(\Omega)}|^2\right)^{1/2}. 
\end{equation}
(See \cite[Chap.~3]{mclean00} for alternative, equivalent norms on $\widetilde{H}^\beta(\Omega)$.) The Dirichlet spectral fractional Laplacian is based on Dirichlet boundary conditions in (\ref{eq:mel-EWP}). Replacing the Dirichlet conditions in (\ref{eq:mel-EWP}) by Neumann boundary conditions leads to a similar operator, the Neumann spectral fractional Laplacian. Throughout the text we will abbreviate 
\begin{equation*}
    \Delta^{\gamma/2}:= - (-\Delta)^{\gamma/2}
\end{equation*}
for the integral and the spectral fractional (both Dirichlet and Neumann) Laplacian. Note that this is consistent with the 1D operator $\Riesz$ introduced earlier. 

\begin{remark}
In the literature there is sometimes confusion regarding the signs of ``the'' fractional Laplacian. This issue already arises for the classical Laplacian as some authors prefer to have it as ``negative'' operator while others prefer a ``positive'' operator. This is just a matter of sign conventions. In dynamical systems it is more common to have the classical Laplacian with nonpositive eigenvalues to indicate that there is decay for the heat semigroup in accordance with the principle that stability is associated with negative real parts for elements in the spectrum. So we will adopt this convention here, which explains our sign conventions.
\end{remark}


For further general references regarding fractional calculus, we refer to~\cite{diethelm04,samko93}, which cover the univariate case in significant  detail. We refer to \cite{delia-etal20} for a good survey of further theoretical and numerical aspects of fractional differential operators. 

\section{Steady States \& Linear Stability: Time Fractional PDEs}
\label{sec:spectral1}

In case of subdiffusion $\alpha\in(0,1)$, the waiting time probability density function from the continuous time random walk modelling is non-exponential and yields a non-Markovian process with memory. Arguing via Fourier and Laplace transforms, e.g., \cite{MetzlerKlafter,Mendez2014,Jin21}, one derives what we refer to here as the subdiffusion equation for $u = u(t,x)\in\R$, $x\in\R$, $t>0$,
\begin{equation}\label{e:subdiffeq}
	\partial_t u =  \fD_{0,t}^{1-\alpha} \sigma \partial_x^2 u,
\end{equation}
with diffusion coefficient $\sigma$ and time-fractional Riemann-Liouville derivative $\fD_{0,t}^{1-\alpha}$. This non-local convolution operator in time with power law kernel entails that the system has memory. Acting by the inverse of $\fD_{0,t}^{1-\alpha}$ on \eqref{e:subdiffeq} yields the distinct but mathematically equivalent form~\cite{SK2005}
\begin{align}\label{e:Caputodiff}
\tilde \fD_{0,t}^{\alpha}u = \sigma\partial_x^2 u,
\end{align} 
with the Caputo fractional derivative $\tilde \fD_{0,t}^{\alpha}$. However, care has to be taken in the functional analytic interpretation of \eqref{e:subdiffeq} and \eqref{e:Caputodiff}, in particular regarding continuity of solutions at $t=0$ \cite{Jin21}. Indeed, replacing the left hand side of \eqref{e:Caputodiff} with $\fD_{0,t}^\alpha(u-u_0)$ gives an alternative form, which requires less smoothness of $u$ with $u(0,x) = u_0(x)$~\cite{MetzlerKlafter,Zacher2019}. The Cauchy problems of  \eqref{e:subdiffeq} and \eqref{e:Caputodiff} are also equivalent to an integro-differential equation that was first investigated via Mellin transforms~\cite{SW1989}, which we do not discuss further here; more details on the comparisons can be found in \cite{MMPG2007}. In all cases, the integral memory means that equations with subdiffusion cannot be viewed as dynamical systems on the phase space of the natural initial condition $u(0,x)$. Starting at time $t>0$ requires suitably prescribing $u(s,x)$ on the elapsed time interval $0\leq s\leq t$, so that such equations may be viewed as equations with monotonically increasing delay. Moreover, as in general for non-autonomous equations, we expect that solutions can cross the initial state and stability properties are not necessarily determined by the spectrum of a time-independent linear operator. 

\subsection{The Linear Subdiffusion Equation}
\label{s:Subdiffusion}

In order to highlight essential differences between classical diffusion and subdiffusion we consider the basic subdiffusion equation in the Riemann-Liouville form \eqref{e:subdiffeq}. Compared to the heat equation, it is more subtle to see that solutions of \eqref{e:subdiffeq} remain positive, and to determine the growth and decay properties. We briefly outline some aspects and refer to, e.g., \cite{MetzlerKlafter} for further details and different perspectives involving the Fox-function. The Green's function of \eqref{e:subdiffeq} has the Fourier transform~\cite[eq.~49]{MetzlerKlafter}
\begin{align}\label{e:FourierSubSol}
	\Fourier\Phi(t,q) = E_\alpha(-\sigma q^2 t^\alpha),
\end{align}
where $E_\alpha(z) = \sum_{n=0}^\infty z^n/\Gamma(1+n\alpha)$ is the Mittag-Leffler function, and $q$ the Fourier-wavenumber. The function \eqref{e:FourierSubSol} has the asymptotic behavior~\cite[eq.~20]{Metzler2004}
\begin{align}\label{e:AsympML}
E_\alpha(-\sigma q^2 t^{\alpha}) \sim 
\begin{cases}
\exp(-\frac{\sigma q^2 t^\alpha}{\Gamma(1+\alpha)}), & t\ll (\sigma q^2)^{1/\alpha},\\
(\sigma q^2 t^\alpha \Gamma(1-\alpha))^{-1}, & t\gg (\sigma q^2)^{1/\alpha}.
\end{cases}
\end{align}
Here we see an effect of memory for nonzero wavenumber: exponential decay for relatively short times, and algebraic decay for long times, here with power $-\alpha$. The separation of these depends in particular on the wavenumber. The algebraic decay effectively replaces the exponential decay of Fourier modes for the classical heat equation, and illustrates that linear stability for processes with subdiffusion has a different character from that for classical diffusion. On the one hand, this behavior is also foreshadowed in Fourier-Laplace space, for which we view the right-hand side of \eqref{e:subdiffeq} realized as an operator with space-time domain $\ACspace([0,T];\Hspace^2(\R))$ into $\Lspace^1(0,T;\Lspace^2(\R))$. With the Laplace transformation $\fL$ we obtain
\begin{align}\label{e:GreenFL}
	(\fL\Fourier{\Phi})(s,q) = (s + \sigma q^2 s^{1-\alpha})^{-1},\quad \Re(s)>0,
\end{align}
where $s\in\C$ is the Laplace space variable. Indeed, one expects that the branch point of \eqref{e:GreenFL} at the origin creates algebraic decay in the inverse Laplace transform. This analysis will be the topic of Section~\ref{s:growthdecay} for reaction subdiffusion models and directly relates to decay in fractional linear ODE. This is always algebraic with rate $-\alpha$, in contrast to classical ODE, where decay is exponential for eigenvalues with negative real parts, cf.\ \cite{Matignon96,BGK21}. On the other hand, in $x$-space a scaling analysis shows that the parabolic scaling of the classical heat equation is replaced in \eqref{e:subdiffeq} by a fractional self-similar scaling: if $u(t,x)$ solves \eqref{e:subdiffeq} then so does $u(\varepsilon^{2/\alpha} t,\varepsilon x)$ for $\varepsilon\in \R$. Based on the similarity variable $x/t^{\alpha/2}$ a series expansion of the Green's function can be found (e.g.\ \cite[eq.~46]{MetzlerKlafter}, \cite[eq.~4.23]{Mainardi2001}), 
\begin{equation}
\label{e:GreensFunction}
\Phi(t,x) = \frac{1}{\sqrt{4\sigma t^\alpha}}\sum_{n=0}^{\infty}\frac{(-1)^n}{n!\Gamma\left(1-\alpha(n+1)/2\right)}\left(\frac{|x|}{\sqrt{\sigma t^\alpha}}\right)^n, \quad t>0. 
\end{equation}
Via this form one can find that $\Phi$ has an expression in terms of the Wright function~\cite{Kilbas2006} and is therefore positive~\cite{Mainardi2001} and decays algebraically locally uniformly in $x$ for $t\gg 1$ with power $-\alpha/2$~\cite{Kochubei1990,EK2004,YR21}. In addition, we remark that, using the method of separation of variables, it was shown in~\cite{MMAK2019} that solutions to the Cauchy problem \eqref{e:subdiffeq} with homogeneous Dirichlet boundary condition decay pointwise algebraically with power $-\alpha$ for $t\gg1$, i.e., faster than on the unbounded domain, which was also found to be optimal in the more general case \eqref{e:Zacher} below. A more detailed discussion of the Mittag-Leffler and Wright functions, their estimates and relation to fractional calculus can be found in \cite{Jin21}; we also refer to \cite{GLM2000,Mainardi2001,MMP2010} for further details of the Wright functions serving as the Green's functions to various (space-)time fractional diffusion equations.

We also briefly outline some aspects of the subdiffusion equation in the Caputo form \eqref{e:Caputodiff}. The formulation via $\tilde \fD_{0,t}^{\alpha}$ motivates to consider a more general form of time fractional evolution equations with a convolution $(k\ast v)(t) = \int_0^t k(t-s)v(s)ds$, $t\geq 0$ for a kernel $k$, generalizing \eqref{e:Caputodiff} to
\begin{equation}
\label{e:convsubdif}
	k \ast\partial_t u =  \sigma \partial_x^2 u.
\end{equation}
The theory for Cauchy-problems in \eqref{e:Caputodiff} and \eqref{e:convsubdif} is by now quite well developed and an excellent introduction is provided by \cite{Jin21}. For \eqref{e:convsubdif} we particularly refer to \cite{Vergara2015} and \cite{KSVZ2016}, which provide decay estimates on bounded and unbounded domains, respectively. For instance, in the subdiffusion case with homogeneous Dirichlet boundary condition the decay estimate of solutions has a different character, and exponential decay of Fourier modes is replaced by algebraic decay with rate $-\alpha$, consistent with the separable solutions to \eqref{e:subdiffeq}; moreover, \cite{Vergara2015} proves $\Lspace^2$-decay with a fixed rate $-\alpha$ for any space dimension. On unbounded domains, the character of algebraic decay persists, but the $\Lspace^2$-decay depends on the spatial dimension $d$: the rate is slower than that on the bounded domain for $d<4$ and is fixed at $-\alpha$ for $d\geq 4$, cf.~\cite{KSVZ2016}. Hence, the decay rate on the bounded domain provides a lower bound for the case of the unbounded domain in any space dimension. In contrast to the classical diffusion, the subdiffusion is so slow that increasing the spatial dimensions (particularly for $n>4$) and restricting to bounded domains does not increase the decay rate. We also refer to \cite{Zacher2019} for various analytical aspects of time-fractional diffusion equations. In addition, to the subdiffusion equation, various examples of kernels $k$ can be found in \cite{Vergara2015}. Lastly, we remark that the subdiffusion equations of the forms \eqref{e:Caputodiff} and \eqref{e:convsubdif} can be extended as ``parabolic-type'' equations by adding uniformly elliptic operators $L$ to the fractional derivatives: $\tilde\fD_{0,t}^\alpha + L$ or $k\ast\partial_t + L$; we refer to \cite{EK2004,Vergara2015,KSVZ2016} for details.

\subsection{Subdiffusion-Reaction Models}
\label{s:subdiffeqs}

Modelling additional reactions and also the extension to systems is a subtle task that can lead to rather different results. Most approaches come as subdiffusion-limited (S) or activation-limited (A) reactions. 

\paragraph{S: Subdiffusion limited.} In the more straightforward and broadly studied class (S), either the non-local operator acts on the diffusion terms as well as the reaction terms~\cite{Henry2006,SWT03,Seki2003,Yuste2004,LHW2007,Nec2008}, or one simply adds an extra source or sink to the right hand side of \eqref{e:Caputodiff}~\cite{VeZa17}, so that in the scalar case we obtain
\begin{align}\label{e:Zacher}
\partial_t u = \fD_{0,t}^{1-\alpha} [\sigma \partial_x^2 u + f(u)]
\quad \text{or}\quad\tilde \fD_{0,t}^{\alpha} u =  \sigma\partial_x^2 u + f(u),
\end{align}
with $f:\R\to\R$. Here the reaction is ``anomalous'' even without subdiffusive transport $\sigma=0$ and for spatially uniform states. Analogously, adding a nonlinear source term to \eqref{e:convsubdif} gives rise to a time fractional evolution equation 
\begin{equation}\label{e:convder}
	k \ast\partial_t u =  \sigma\partial_x^2 u + f(u).
\end{equation}
The theory for Cauchy-problems in \eqref{e:Zacher} is well developed (mostly on bounded domains) and an excellent introduction is provided by \cite{Jin21}, for which we also refer to the numerous references in this active field. The initial-boundary value problem of \eqref{e:convder} and even more general semilinear diffusion cases have been considered in \cite{VeZa17,KTT2020}, which provide well-posedness and stability theory. The existence problem of steady states in \eqref{e:Zacher} and \eqref{e:convder} is exactly the same as for classical reaction-diffusion systems, although the stability properties may differ. For instance, the linear stability principle has a different character, and exponential decay in the stable case is replaced by algebraic decay with rate $-\alpha$. Indeed, \eqref{e:Zacher} for linear $f(u)=bu$ and for a Fourier mode $\tilde q$ such that $b<\sigma \tilde q^2$ has the same form as the subdiffusion equation \eqref{e:subdiffeq} for the Fourier mode $q=\sqrt{\tilde q^2-b/\sigma}$, consistent with the power $-\alpha$ in \eqref{e:AsympML}. There are numerous studies of asymptotic stability for equations of the form \eqref{e:convder}; we refer to  \cite{VeZa17,KTT2020} as abstract studies and the references therein. In addition, some recent studies concern the asymptotic stability for modified forms of \eqref{e:convder}, e.g., with fractional Laplacian~\cite{DKT2023} or nonlinearity involving infinite delays~\cite{NTN23}. Another difference in stability properties is that spatially constant steady states tend to be spectrally more stable against non-steady perturbations in the presence of subdiffusion of class (S) than classical diffusion, e.g.\ \cite{Matignon96,BGK21,KP23}. We refer to \cite{BGK21} for a recent account of stability in fractional ODE.

Both aspects, i.e., non-expoential decay and (de)stabilizing effect, also arise for the less considered class (A) on which we focus in the following. Here the reaction kinetics is not affected by subdiffusive transport, i.e., any nonlocal time-dependent operator acts on the diffusion terms only  \cite{Henry2000,Henry2006,Fedotov2010,Nec2007,Nec2008}. Within this class of models we distinguish the two cases of added source/sink (\AA) and linear creation/annihilation (\AL) as detailed next. 

\paragraph{\AA: subdiffusion with added source/sink.} Akin to classical reaction-diffusion models, following \cite{Henry2000}, in this case one simply adds an extra source or sink to the right-hand side of \eqref{e:subdiffeq}. More generally, for a possibly nonlinear reaction term $f(u)$ one thus obtains 
\begin{equation}\label{e:Henry2000}
\partial_t u = \sigma \fD_{0,t}^{1-\alpha} \partial_x^2 u + f(u). 
\end{equation}
While for the classical diffusion at $\alpha=1$ the usual maximum principle holds,  the Green's function for \eqref{e:Henry2000} with negative linear reaction dynamics, e.g., $f(u)=-u$, partially takes negative values, cf.\ \cite{Henry2006}. Heuristically, the reason is that the sink can remove particles much faster than they can be transported from other positions. Hence, \eqref{e:Henry2000} is not a model for an absolute density $u$, see also \cite{Nepomnyashchy2016}. However, when $u$ represents a density perturbation from a saturated strictly positive state this phenomenon is acceptable, at least for short times. Within \eqref{e:Henry2000} this would concern $u\approx u_*$, where $u_*>0$ solves $f(u_*)=0$. A corresponding ad hoc model, cf.\ \cite{Henry2002}, for multiple species is
\begin{equation}\label{e:Henry2000-system}
\partial_t \u =\Sigma \fD_{0,t}^{1-\alpha} \partial_x^2 \u + \F(\u),
\end{equation}
where $\u=\u(t,x) \in\R ^N$ is the vector of the perturbation densities, $\Sigma=\diag(\sigma_1,\ldots,\sigma_N)\in \R^{N\times N}$ is the diffusion matrix with positive entries, and $\F :\R ^N\to\R ^N$ is a possibly nonlinear function that models the reaction kinetics. 

\paragraph{\AL: subdiffusion with linear creation/annihilation.} In this case, the addition or removal of particles through (linear) reactions occurs only during the waiting time \cite{Henry2006} and the resulting scalar equation with reaction rate $a$ reads
\begin{equation*}
\partial_t u = \sigma e^{at} \fD_{0,t}^{1-\alpha} ( e^{-at} \partial_x^2 u ) + a u. 
\end{equation*}
From the corresponding integral form one can infer that this model preserves positivity of solutions since the amount of removed particles is less than the amount of existing particles. For multiple species and reaction rate matrix $\A\in\R^{N\times N}$ the model 
\begin{equation}\label{e:Langlands2008}
\partial_t \u = \Sigma e^{\A t} \fD_{0,t}^{1-\alpha} ( e^{-\A t} \partial_x^2 \u ) + \A \u,
\end{equation}
has been proposed in~\cite{Langlands2008,Nepomnyashchy2016}, and rigorously derived in \cite{AL21}. Since the matrix exponentials are suitable for linear reactions only, the form of an appropriate model for nonlinear reactions is not clear. We also refer to \cite{Nepomnyashchy2016} for a discussion of the modelling of reaction subdiffusion equations.

\subsection{Spectra, Growth and Decay}
\label{s:growthdecay}

The Fourier-Laplace transform can be employed also for systems to derive dispersion relations of temporal and spatial modes, analogous to the classical reaction diffusion systems.  We illustrate decay and growth in subdiffusion from this viewpoint. In contrast to classical diffusion, the fractional powers in the dispersion relation yield branch points, and we already alluded to the related algebraic decay in the subdiffusion equation in Section~\ref{s:Subdiffusion}. Also generally, the branch points strongly impact the decay laws determined by an inverse Laplace transform, and moreover the decay properties of models (\AA) and (\AL) differ significantly, although the models are equal in absence of reactions. In addition, the solutions to the dispersion relation to the left of the most unstable branch point in the complex plane depend on the choice of branch cut, and following \cite{YR21} we therefore refer to these solutions as pseudo-spectrum; when there is no ambiguity in this section we refer to all solutions of dispersion relations (in the admissible domain) as pseudo-spectrum. Non-canonical branch cuts admit to relate pseudo-spectrum to regular spectrum of the classical diffusion equations through the limit $\alpha\to 1$. For illustration, we consider two-component systems and start with classical diffusion, i.e., a system
\begin{equation}
\label{e:RegAISystem}
	\partial_t \u = \Sigma \partial_x^2 \u + \F(\u), 
\end{equation}
where $\u=\u(t,x) \in\R ^2$ is the vector of density of the species  with $x\in\R$, $t>0$, $\F:\R ^2\to\R ^2$ represents the reaction kinetics, and $\Sigma =\diag(1,\sigma)$ is the diagonal matrix of positive constant diffusion coefficients. Here the first diffusion coefficient is normalized to $1$ by a suitable scaling of $x$ so that $\sigma$ is in fact the ratio of diffusion coefficients. Suppose that \eqref{e:RegAISystem} possesses a spatially homogeneous steady state $\u _*\in\R^2$, i.e., $\F (\u _*)=0$. The linearization of \eqref{e:RegAISystem} about $\u _*$ is given by
\begin{equation}
\label{e:RegLinearSystem}
	\partial_t \u = \Sigma \partial_x^2 \u +\A \u , \quad \A := \left(\txtD_\u \F(\u_*)\right)=
	\begin{pmatrix}
	a_{1}	&	a_{2}\\
	a_{3}	&	a_{4}
	\end{pmatrix}.
\end{equation}
The right-hand side of \eqref{e:RegLinearSystem} defines the linear operator $\cL :=\Sigma \partial_x^2 +\A $, with eigenvalue problem 
\begin{equation*}
	\cL \u  = \Sigma \partial_x^2 \u +\A \u  = s\u.
\end{equation*}
Upon Fourier transforming $x\in \R$ to Fourier-wavenumbers $q\in\R$, the eigenvalue problem becomes
\begin{equation*}
	\cL (\Fourier\u) = -q^2\Sigma (\Fourier\u) + \A (\Fourier\u) = s(\Fourier\u),
\end{equation*}
whose solvability condition as a matrix eigenvalue problem yields the dispersion relation
\begin{equation}\label{e:RegDR}
	D_{\reg}(s,q^2) := \det(s\I + q^2 \Sigma - \A ) = (s + q^2 - a_{1})(s + \sigma q^2 - a_{4}) -a_{2} a_{3} = 0.
\end{equation}
The solution set $\Lambda_\reg := \{ s \in \C : D_{\reg}(s,q^2) = 0 \text{ for a } q \in \R \}$ defines in particular the $\Lspace^2$-spectrum of $\cL$ with domain $(\Hspace^2(\R))^2$, e.g., \cite{Sandstede2002}, which we refer to as regular spectrum. The homogeneous steady state $\u _*$ is called strictly spectrally stable (unstable) if $\max(\Re(\Lambda_\reg))<0$ ($>0$). It is then also linearly and nonlinearly stable (unstable) for \eqref{e:RegLinearSystem} and \eqref{e:RegAISystem}, respectively \cite{Sandstede2002}.

\paragraph{Dispersion relation for \AA:} Proceeding in the same way for \eqref{e:Henry2000-system} gives 
\begin{equation*}
	(s\I + s^{1-\alpha}q^2 \Sigma - \A )\, \fL\Fourier\u = \Fourier\u_0, \quad \Re(s)>0, \ q \in \R,
\end{equation*}
with the constraint on $s$ due to the branch point of $s^{1-\alpha}$ at the origin. Hence, the analogue of \eqref{e:RegDR} is the dispersion relation on $\Omega_0^+ := \{s\in \C\setminus\{0\} : \arg(s)\in (-\pi/2,\pi/2) \}$ given by 
\begin{align*}
	D_\ss(s,q^2) & := \det(s\I +s^{1-\alpha}q^2\Sigma -\A )
	 = (s+s^{1-\alpha} q^2- a_{1})(s+s^{1-\alpha} \sigma q^2- a_{4})- a_{2}a_{3}=0.
\end{align*}
For any given $q\in\R$ we can extend the domain of $D_\ss$ on the principal branch such that the values on positive reals are preserved. Due to the discrete solution set, we can choose $\theta_1(q)\in(0,\pi/2)$ such that there is no solution of the dispersion relation on the branch cut 
\[
\BC_{0}^{\theta_1(q)} := \{ s\in\C : \Im(s)/\Re(s) = \tan(\theta_1(q)),\, \Re(s)< 0\},
\]
i.e., $D_\ss(s,q^2) \neq 0$ for  $s\in\BC_{0}^{\theta_1(q)}$. For simplicity we suppress the dependence of $\theta_1$ on $q$, and obtain the extended domain of $D_\ss$ as
	\[
	\Omega_0:=\{s \in \C\setminus\{0\} : \arg (s) \in (-\pi+\theta_1,\pi+\theta_1),\,\theta_1 \in (0,\pi/2)\}.
	\]
The pseudo-spectrum $\Lambda_\ss$ is the corresponding set of roots of $D_\ss$ on $\Omega_0$. 

\paragraph{Dispersion relation for \AL:} For \eqref{e:Langlands2008} we follow \cite{Langlands2008,YR21,AL21} and make the simplifying assumption that $\A$ is diagonalizable, so that with the eigenvalues $\mu_1, \mu_2$ of $\A$, the coordinate change
\begin{equation*}
	\P =
	\begin{pmatrix}
	a_{2}			&	a_{2}\\
	\mu_1-a_{1}	&	\mu_2-a_{1}
	\end{pmatrix},\quad
	\P ^{-1}=
	\begin{pmatrix}
	\frac{a_{1}-\mu_2}{a_{2}(\mu_1-\mu_2)}	&	\frac{1}{\mu_1-\mu_2}\\
	-\frac{a_{1}-\mu_1}{a_{2}(\mu_1-\mu_2)}	&	-\frac{1}{\mu_1-\mu_2}
	\end{pmatrix}
\end{equation*}
gives $\bar{\A}:=\P ^{-1}\A \P =\mathrm{diag}(\mu_1,\mu_2)$. We order the eigenvalues so that $\Re(\mu_1)\geq \Re(\mu_2)$. Setting $\u =\P \w $ turns \eqref{e:Langlands2008} into
\[
	\partial_t\w   = \bar{\Sigma }e^{\bar{\A}t}\fD_{0,t}^{1-\alpha}(e^{-\bar{\A}t}\partial_x^2 \w )+\bar{\A}\w , 
	\quad \bar{\Sigma }= \begin{pmatrix}
d_{1} & d_{2}\\
d_{3} & d_{4}
\end{pmatrix}:=\P ^{-1}\Sigma \P.
\]
Fourier-Laplace transforming then gives the dispersion relation
\begin{equation*}
\begin{split}
 D_\ca(s,q^2):=\ & \det(s\I -\bar{\A}+(s\I -\bar{\A})^{1-\alpha}\bar{\Sigma }q^2) 
  = \det((s\I -\bar{\A})^{1-\alpha})\det((s\I -\bar{\A})^{\alpha}+\bar{\Sigma }q^2)\\
  =\ & (s-\mu_1)^{1-\alpha}(s-\mu_2)^{1-\alpha} \big(((s-\mu_1)^\alpha+d_{1}q^2)((s-\mu_2)^\alpha+d_{4}q^2)-d_{2}d_{3}q^4\big) = 0,
\end{split}
\end{equation*}
where $\Re(s-\mu_1),\,\Re(s-\mu_2)>0$. As above, we may extend the domain to $s\in\Omega_\ca\subset \C$ by choosing branch cuts that connect $\infty$ to the finite branch points $\mu_1$ and $\mu_2$, respectively. It turns out, cf.\ \cite{YR21}, that $s=\infty$ is a branch point of $(s-\mu_1)^{\alpha}(s-\mu_2)^{\alpha}$ for $\alpha\in(0,1/2)\cup(1/2,1)$; for $\alpha=1/2$ we can think of $\infty$ as an intermediate point on the branch cut connecting $\mu_1$ and $\mu_2$. For simplicity, we assume real $\mu_1>\mu_2$; the general case can be found in \cite{YR21}. In order to keep non-integer powers of positive reals on the positive real axis, we then choose branch cuts $\BC_{\mu_j}^{\theta_j}$, $j=1,2$, such that 
\[
\Omega_\ca := \{s \in \C\setminus\{\mu_1,\mu_2\} : \arg(s-\mu_j) \in (-\pi+\theta_j,\pi+\theta_j),\  \theta_j\in (0,\pi/2),\ j=1,2 \}.
\] 
As for (\AA), given $q\in\R$ we select $\theta_1,\theta_2$ such that the solutions of $D_\ca(\cdot,q^2)=0$ do not lie on the branch cuts. The pseudo-spectrum $\Lambda_\ca$ is the corresponding set of roots of $D_\ca$  on $\Omega_\ca$.

\begin{remark}
We illustrate the special nature of rational $1-\alpha=n/m\in\mathbb{Q}$ for the non-trivial factor of $D_\ca$, 
\begin{align*}
D_{\ca2}(s,q^2) &:= \big((s - \mu_1)^\alpha + d_{1} q^2\big) \big((s - \mu_2)^\alpha + d_{4} q^2\big) - d_{2} d_{3} q^4 = 0, \quad s \in \Omega_\ca.
\end{align*}
In the rational case it can be recast as a polynomial with respect to the two variables $z_j=(s-\mu_j)^{1/m}$, $j=1,2$, of degree $m-n$ in each variable, given by
\[
\tilde D_{\ca2}(z_1,z_2,q^2) = (z_1^{m-n} + d_1q^2)(z_2^{m-n} + d_4q^2) - d_2d_3q^4.
\] 
By resultant theory, the polynomial system $\tilde D_{\ca2}(z_1,z_2,q^2) = 0$, $z_1^m + \mu_1 = z_2^m + \mu_2$ has a finite number of roots $(z_1,z_2)$, which gives a finite number of roots $s$ for $D_{\ca2}(s,q^2)$, including multiplicities.
\end{remark}

\paragraph{Growth and decay:} Via an inverse Laplace transform, the location of pseudo-spectrum translates into growth and decay properties in the time-domain for each Fourier-wavenumber. For case (\AA), we quote a simplified form of \cite[Theorem 4.2]{YR21}, where the following set of solutions with maximal real part is relevant:
\[
S^+:=\{(s,q)\in \overline{\Omega_0^+}\times \R : D_\ss(s,q^2)=0\text{ and } \Re(s) \mbox{ is maximal for } q \}.
\]

\begin{theorem}[Case \AA]\label{t:ILT-ss}
	Let $\alpha = \tfrac n m\in(0,1), \,n,m\in\N$, be a reduced fraction,  $\lambda:=\sup(\Re(\Lambda_\ss))$, and $\Fourier\u(t,q)$ a solution to the Fourier transformed linearization of \eqref{e:Henry2000-system} in $U_*$ with  $\Fourier\u(0,q)=\Fourier\u_0(q)$. 
	\begin{itemize}
		\item[(1)] If $\lambda \geq 0$ then for any $(s_0,q_0)\in S^+$ we have $\Fourier\u(t,q_0) = C (q_0,\Fourier\u_0(q_0)) t^{k-1}e^{s_0 t} + o(t^{k-1}e^{\Re(s_0)t})$ with $C(q_0,\Fourier\u_0(q_0))\in \C^2$ nonzero for almost all $\Fourier\u_0$ and $k$ the multiplicity of $s_0$ as a root of $D_\ss(\cdot,q_0^2)=0$.
		\item[(2)] Assume $\det(\A)\neq 0$. If (i) $\lambda=0$ and $S^+=\emptyset$ or (ii) $\lambda < 0$ or (iii) $\Lambda_\ss=\emptyset$, then for any $q\neq 0$ there is a $q$-independent invertible $\M\in\R^{2\times 2}$ such that $\Fourier\u(t,q) =  q^2\M \Fourier\u_0(q) t^{\alpha-2} + \cO(t^{\alpha-2-1/m})$. In addition, the entries of $\M$ approach $0$ as $\alpha\to 1$ and  $\infty$ as $\det(\A)\to 0$. 
	\end{itemize}
\end{theorem}

Concerning item (1) we note that $S^+\neq \emptyset$ if $\lambda>0$. In item (2) the condition $\det(\A)\neq 0$ has been added compared to \cite[Theorem 4.2]{YR21}, where this was assumed a priori, cf.\ Section~\ref{s:Turinginstability}. In case $q=0$ we obtain $\partial_t (\Fourier\u) =  \A (\Fourier\u)$, whose solutions are either constant, grow algebraically as in item (1) for $s_0=0$, or decay exponentially. Theorem~\ref{t:ILT-ss} means that the growth of Fourier modes is the same as for classical diffusion, $\alpha=1$, but the decay properties differ: in the stable case each Fourier mode decays algebraically as $t^{\alpha-2}$. Notably, this is in contrast to the decay law $t^{-\alpha}$ for the subdiffusion-limited class (S) discussed in Section~\ref{s:subdiffeqs}, see also~\cite{Matignon96,LHW2007,BGK21}. For the case (\AL) with $\mu_1>\mu_2$ we quote a simplified form of \cite[Theorem 5.3]{YR21}. Here we denote $\Omega_\ca^- := \{s\in\Omega_\ca:\Re(s)<\Re(\mu_1)\}$, $\Omega_\ca^{0+} := \Omega_\ca \setminus \Omega_\ca^-$ and 
\begin{align*}
S^+&:=\{(s,q)\in\Omega_\ca^{0+}\times (\R\setminus\{0\}): D_\ca(s,q^2)=0,\text{ and }\Re(s)\text{ maximal for } q\},\\
S^-&:=\{(s,q)\in\Omega_\ca^-\times (\R\setminus\{0\}): D_\ca(s,q^2)=0, \text{ and }\Re(s)\text{ maximal for } q\}.
\end{align*}

\begin{theorem}[Case \AL]\label{t:ILT-ca}
Assume $\mu_1>\mu_2$. Let $\alpha=\tfrac n m\in(0,1),\,n,m\in\N$, be a reduced fraction, $\lambda:=\sup(\Re(\Lambda_\ca))$, and $\Fourier\u(t,q)$ a solution to the Fourier transform of \eqref{e:Langlands2008} with  $\Fourier\u(0,q)=\Fourier\u_0(q)$. 
	\begin{itemize}
	\item[(1)] If $\lambda\geq\mu_1$, then for any $(s_0,q_0)\in S^+$ we have $\Fourier\u(t,q_0) = C(q_0,\Fourier\u_0(q_0))t^{k-1} e^{s_0t} + o(t^{k-1}e^{\Re(s_0)t})$ with $C(q_0,\Fourier\u_0(q_0))\in\C^2$ nonzero for almost all initial data and $k$ the multiplicity of $s_0$ as a root of $D_{\ca2}(\cdot,q_0^2)=0$.
	\item[(2)] If (i) $\lambda=\mu_1$ and $S^+=\emptyset$ or (ii) $\lambda\in(\mu_2,\mu_1)$, then $S^-\neq \emptyset$ and for any $(s_0,q_0)\in S^-$ there exist $q$-independent nonzero $c_j\in\C,\,j=1,2$, and a singular $\M_1(q)\in\C^{2\times 2}$ with one-dimensional kernel such that for any initial data satisfying $\P^{-1}\Fourier\u_0(q_0)\notin\ker(\M_1(q_0))$ we have
\begin{align*}
\Fourier\u(t,q_0) =&\ q_0^{-2} \P
	\begin{pmatrix}
c_1 & 0\\ 0 & c_2 t^{-1} 
\end{pmatrix}\M_1(q_0)\P^{-1}\Fourier\u_0(q_0)\, e^{\mu_1 t} t^{-\alpha}+ \cO(t^{-\alpha-1/m} e^{\mu_1 t}).
\end{align*}
	\item[(3)] If (i) $\lambda\leq\mu_2$ or (ii) $\Lambda_\ca=\emptyset$, then for any $q\neq 0$ and any initial data satisfying $\P^{-1}\Fourier\u_0(q)\notin\ker(\M_1(q))$, with $c_j,\M_1(q)$ as in item (2), it holds that 
\begin{align*}
\Fourier\u(t,q) =&\ q^{-2} \P
	\begin{pmatrix}
c_1 & 0\\ 0 & c_2 t^{-1} 
\end{pmatrix}\M_1(q)\P^{-1}\Fourier\u_0(q)\, e^{\mu_1 t} t^{-\alpha} +  \cO(t^{-\alpha-1/m} e^{\mu_1 t}).
\end{align*}
	\end{itemize}
In addition, $c_j\to 0$ as $\alpha\to1$ for $j=1,2$, and the entries of $\M_1(q)$ are bounded for any $q\neq 0$.
\end{theorem}

Concerning item (1), if $\lambda>\mu_1$, then $S^+$ is guaranteed to be non-empty. Items (2,3) can be understood as thresholds for decay rates. In particular, exponential decay for any Fourier mode cannot be faster than that of the kinetics. For initial data in the kernel of $\M_1(q)$, the leading order terms of the corresponding solution is rather involved; we refer to \cite[Theorem 5.3]{YR21} for further details. The statement also holds for $\mu_1=0$ and scalar equations, and thus includes the subdiffusion equation \eqref{e:subdiffeq}. The decay law $t^{-\alpha}$ (for rational $\alpha$) indeed coincides with that derived in Section~\ref{s:Subdiffusion} via \eqref{e:AsympML}, and with that for \eqref{e:Zacher}, e.g., in \cite{Matignon96,Vergara2015}. As for Theorem~\ref{t:ILT-ss} the method of proof only allows for rational $\alpha$. In contrast to Theorem \ref{t:ILT-ss}, the coefficients decay as $q^{-2}$ in Theorem \ref{t:ILT-ca} so that some smoothing of initial data is conceivable.

In summary, a pseudo-spectrum with positive real parts readily implies exponential instability with the same spectral rate that would arise in classical diffusion with the same real parts. Negative real parts, however, do not imply the exponential decay of the classical case: In case (\AA) there is only algebraic decay (in systems at least for rational $\alpha$) of power $\alpha-2$, and indeed the origin is always a branch point of the dispersion relation. In case (\AL) the branch points lie at the eigenvalues of $\A$, which admits exponential decay if the pseudo-spectra have negative real parts. However, the decay cannot be faster than the rate of the most unstable eigenvalue of $\A$, and there is an algebraic correction with power $-\alpha$. From a technical point of view, the essence of two distinct decays is that in case (\AA) the branch point does not solve the dispersion relation, while it does in case (\AL). As mentioned, also in model class (S) strictly stable spectrum is expected to give only algebraic decay, which is proven by energy methods for Dirichlet boundary conditions (together with further results) in \cite{Vergara2015,VeZa17}; we also refer to Section \ref{sec:gradient} for more details on global/gradient-flow methods based upon energy dissipation functionals.

In the results for subdiffusion models of class (A), $\alpha$ is taken rational since then the inverse Laplace transform can be computed explicitly via residues. Although the results suggest uniformity, there seems to be no proof. Moreover, the estimates for a selected Fourier wavenumber, and although naturally locally uniform and with asymptotic rates, there appears to be no additional results in $x$-space. Both issues may be approachable by energy methods. For the case (\AA), i.e., \eqref{e:Henry2000-system}, constants grow as $q^2$, thus requiring sufficient smoothness in $x$-space. In contrast, in the case (\AL), i.e., \eqref{e:Langlands2008}, the analogous coefficients decay as $q^{-2}$, which indicates some smoothing. However, results on well-posedness for these models seem not to be available in the literature. 

\subsection{Turing-type instability}
\label{s:Turinginstability}

Concerning the onset of instability, models \eqref{e:Henry2000-system} and  \eqref{e:Langlands2008} differ significantly from each other and from the case of normal diffusion $\alpha=1$. We quote some results for which we refer to \cite{YR21} and the references therein, in particular \cite{Nec2007,Nec2008,Henry2002,Henry2005}. These concern the case of two species systems, $N=2$, and when parameters are such that $\alpha=1$ admits a so-called Turing instability for a critical diffusion ratio $\sigma>\sigma_\reg$. For this instability the onset at $\sigma=\sigma_\reg$ occurs with a unique finite wavenumber $q=q_\reg$, while $\A$  has eigenvalues with strictly negative real parts, i.e.,  $\tr(\A ) = a_{1} + a_{4} <0$ and $\det(\A ) = a_{1}a_{4} - a_{2}a_{3} >0$. Without loss of generality we may assume $a_1>0$, $a_4<0$, $a_2a_3<0$, which gives $\sigma_\reg>1$ via the well known instability condition $\sigma a_1+a_4>0$, e.g., \cite{Murray2}. Turing instability for class (S) is well studied: roughly speaking, the subdiffusion does not have an effect on the steady instability threshold; we refer to \cite{LHW2007,Nec2008,KP23} for further details. However, this is not the case for models of class (A).

\paragraph{\AA: destabilizing subdiffusion.} For \eqref{e:Henry2000-system}, if the spectrum for $\alpha=1$ is Turing unstable, then the spectrum for all $\alpha\in(0,1)$ is unstable. Somewhat counterintuitively, this means that, in terms of the diffusion ratio $\sigma$, the reaction-subdiffusion in case (\AA) is always less stable than classical diffusion. In particular, the diffusion instability threshold $\sigma_\text{add}$ for (\AA) is smaller than that for classical diffusion. A more significant deviation from parabolic equations is that the spectrum becomes unstable through infinite wavenumbers, 
in the following sense: 
for any $\sigma$ the grow rate at $q=\pm\infty$ is zero, and for
$\sigma>\sigma_\text{add}$ the unstable modes have wavenumber $q$ with $|q|\in [q(\sigma),\infty)$, where $q(\sigma)\to\infty$ as $\sigma\searrow \sigma_\text{add}$. In particular, the signature finite wavenumber selection at the onset of a Turing instability is absent. At the same time, in the limit $\alpha\to 1$ the pseudo-spectrum converges locally uniformly as a set in $\C$ to the regular spectrum. Hence, there is also an instability threshold in terms of $\alpha$. We refer to \cite{YR21} for some numerical explorations.

\paragraph{\AL: stabilizing subdiffusion.}
For \eqref{e:Langlands2008} the analytical comparison of regular and pseudo-spectrum is constrained to $s\approx 0$. In this region, the pseudo-spectrum is real and more stable than normal diffusion, and pseudo-spectrum becomes unstable via finite wavenumber at a diffusion ratio $\sigma_\text{lin}>\sigma_\reg$, if $\alpha\geq\alpha_\A$ with a certain threshold $\alpha_\A$. For $\alpha<\alpha_\A$ the instability is completely suppressed. In the limit $\alpha\to 1$ the pseudo-spectrum again converges locally uniformly as a set in $\C$ to the regular spectrum. This means that for $\alpha\approx1$ the critical pseudo-spectrum is real and therefore it is globally more stable than regular spectrum. Numerical computations suggest that this relation holds broadly in parameter space{\blue, cf. \cite{YR21}}.

\section{Gradient Flows: Time Fractional PDEs}
\label{sec:gradient}

{
\renewcommand{\d}{\text{\rm d}}
\newcommand{\vep}{\varepsilon}
\newcommand{\lam}{\lambda}
\newcommand{\disp}{\displaystyle} 
\newcommand{\epsi}{\varepsilon} 
\newcommand{\e}{{\rm e}}

In the last section, we have seen how to analyze dissipative subdiffusion equations locally by focusing on spectra and bifurcations. In this section, we describe the natural complementary view on this class of differential equations via global gradient flow techniques. Gradient flows have been vigorously studied in the field of nonlinear PDEs. An abstract form of gradient flows in an appropriate space ${\cX}$ (e.g., Hilbert and Banach spaces) reads,
\begin{equation}\label{GF}
\partial_t u(t) = - \nabla E(u(t)) \ \mbox{ for } \ t > 0, \quad u(0) = u_0,
\end{equation}
where $\partial_t = \d/\d t$, $u : [0,\infty) \to {\cX}$ is a gradient flow emanating from an initial datum $u_0$ and $\nabla E$ denotes a functional derivative in a proper sense (e.g., Fr\'echet derivative and subdifferential) of an energy (or entropy) functional $E : {\cX} \to \bar\R:=[-\infty,\infty]$. There are many parabolic equations (systems), e.g., Cahn-Hilliard system, Allen-Cahn equation, nonlinear diffusion equations such as porous medium and fast diffusion equations, many of phase transition models, which admit variational structures of gradient flow type. One of crucial features of gradient flows is the decrease (in time) of the energy $t \mapsto E(u(t))$. Moreover, any elements of the $\omega$-limit set of each gradient flow are equilibria, i.e., critical points of the functional $E$. On the other hand, the non-emptiness of the $\omega$-limit set is nontrivial, and moreover, even if it is nonempty, it is rather delicate whether or not the $\omega$-limit set is singleton. 

This section is concerned with a fractional gradient flow, which is a variant of gradient flows with time-fractional derivatives instead of the classical one, say
\begin{equation}\label{fracGF}
\partial_t^\alpha (u - u_0) (t) = - \nabla E(u(t)) \ \mbox{ for } \ t > 0,
\end{equation}
where $\alpha \in (0,1)$ and $\partial_t^\alpha:=\fD_{0,t}^{\alpha}$ is a shorthand notation for the Riemann-Liouville derivative throughout this section.

\subsection{Brief History}

The study of fractional gradient flows may date back to the 1970s, when the following nonlinear Volterra equation was studied extensively:
\begin{equation}\label{Vol}
u(t) + (a*\xi)(t) = g(t), \quad \xi(t) \in A(u(t)), \quad t > 0,
\end{equation}
where $A : {\cX} \to {\cX}$ is a (possibly multivalued) nonlinear operator in a Hilbert or Banach space ${\cX}$, $g : (0,T) \to {\cX}$ is given and $a : [0,+\infty) \to [0,+\infty]$ is a kernel function of an appropriate class. Actually, the convolution of both sides of \eqref{fracGF} with the conjugate kernel $k_{\alpha}$ yields \eqref{Vol} with the choice $a = k_{\alpha}$, $A = -\nabla E$ and $g(t) \equiv u_0$, since $k_{\alpha} * k_{1-\alpha}$ is identically equal to $1$. In those days, equation \eqref{Vol} was often studied in a spirit of the nonlinear semigroup theory (see, e.g.,~\cite{Komura,CL,HB1}).

Barbu~\cite{B75,B77,B-dn} studied an abstract nonlinear Volterra equation \eqref{Vol}, which arises in the study of mechanical systems with memory effects, under the assumptions that $a$ is of class $W^{1,1}_{\rm loc}([0,+\infty))$ (hence it is differentiable almost everywhere and finite at $t = 0$) and positive, $g \in W^{1,2}(0,T;{\cX})$ and $A = \partial \varphi$ in a Hilbert space ${\cX}$. Differentiating both sides of \eqref{Vol} in time, we see that
$$
u'(t) + a(0) \xi(t) + \big( a' * \xi \big)(t) = g'(t), \quad \xi(t) \in A(u(t)), \quad 0 < t < T,
$$
which can also be regarded as a nonlocal (in time) perturbation problem of a (local) nonlinear evolution equation. We also refer the reader to~\cite{CrNo78,Gri78,LoSt78,KiSt79,Cl80,Aiz93}. However, these results cannot be applied to the Riemann-Liouville kernel $k_\beta$, for $k_\beta'(0+)$ is always infinite. As for singular kernels such as the Riemann-Liouville kernel $k_\beta$, Cl\'ement and Nohel~\cite[Theorem 3.1]{ClNo81} studied \eqref{Vol} for \emph{completely positive kernels} $a \in L^1_{\rm loc}([0,+\infty))$ and proved existence of a generalized solution, which is a limit of a certain class of approximate solutions. Moreover, the abstract result is also applied to a couple of nonlocal (in time) PDEs (see also~\cite{ClNo79}). We further refer the reader to ~\cite{Gri79} and~\cite[Theorems 2 and 3]{Ki80}, which are closely related to the main focus of this section. On the other hand, in these works, no application to the time-fractional problem is  explicitly mentioned. Evolution equations including nonlocal derivatives (e.g., Riemann-Liouville derivative) are studied in~\cite{Cl84,Gri85,CP90} by finding out that nonlocal differential operators are $m$-accretive in Bochner spaces (see Section~\ref{Ss:result} below). On the other hand, most  existence results are established for generalized solutions and, to the best of our knowledge, there had been no result which can guarantee the existence of strong (in time) solutions for time-fractional abstract evolution equations of gradient flow type prior to  Zacher's pioneer work~\cite{Za09}, which is still concerned with only some linear evolution equations and will be more detailed below. One of the reasons for such a stagnation might reside in the lack of useful chain-rule formulae, which play an indispensable role in the analysis of (classical) gradient flows \eqref{GF}.

Let $V$ and $V^*$ be a Hilbert space and its dual space, respectively, and let $H$ be a pivot Hilbert space such that the Gelfand triplet holds,
$$
V \hookrightarrow H \simeq H^* \hookrightarrow V^*.
$$
Based on the Galerkin approximation, Zacher~\cite{Za09} proved existence and uniqueness of strong (in time) solutions to the linear evolution equation,
\begin{equation}\label{Zacher}
\dfrac{\d}{\d t} \big( \left[ k*(u-u_0) \right](t), v \big)_H + b(t,u(t),v) = \langle f(t), v \rangle_V \ \mbox{ for all } v \in V,
\end{equation}
where $(\cdot,\cdot)_H$ and $\langle \cdot, \cdot \rangle_V$ stand for an inner product of $H$ and a duality pairing between $V$ and $V^*$, respectively, $k$ is a completely positive kernel, $f : (0,T) \to V^*$ and $u_0 \in H$ are given and $b(t,\cdot,\cdot)$ is a (time-dependent) bounded coercive bilinear form defined on $V$. Key devices for proving the existence of strong (in time) solutions to \eqref{Zacher} consist  not only of the $m$-accretivity of the nonlocal time-differential operator but also of a nonlocal energy identity (which is a sort of the chain-rule formula),
\begin{align}
\MoveEqLeft{
\left( \dfrac \d {\d t} \left(h * u \right)(t), u(t) \right)_H
= \dfrac 1 2 \dfrac \d {\d t} \left( h * \|u\|_H^2 \right)(t)
+ \dfrac 1 2 h(t) \|u(t)\|_H^2 \nonumber
}\\
 &- \dfrac 1 2 \int^t_0 h'(s) \|u(t)-u(t-s)\|_H^2 \, \d s \label{FLI}
\quad \mbox{ for a.e. } t \in (0,T)
\end{align}
for any $h \in W^{1,1}(0,T)$ and $u \in L^2(0,T;H)$ (cf.~Proposition \ref{P:chain} below). Moreover, in~\cite{VeZa08}, convergence of nonlocal (in time) gradient flows in Euclidean spaces $\R^d$ to equilibria is proved under certain assumptions, which seem to exclude Riemann-Liouville and Caputo derivatives but can treat their regularizations. Let us also briefly review other contributions of R.~Zacher and his collaborators (e.g., V.~Vergara) without any claim of completeness. In~\cite{Za08}, boundedness (both in space and time) of solutions to the time-fractional PDE,
\begin{equation}\label{linDif}
 \partial_t^\alpha (u-u_0) - \mathrm{div}\left(A(x,t)\nabla u\right) = g(x,t),
\end{equation}
where $A(x,t)$ is a uniformly elliptic but possibly discontinuous matrix field and $u_0$ and $g$ are prescribed, as well as its quasilinear variants (e.g., $A$ also depends on $u$) is proved by means of De Giorgi's iteration technique (see also~\cite{Za13}). This result was extended to degenerate and singular diffusion equations of $p$-Laplacian type in~\cite{VeZa10}, and moreover, a weak Harnack inequality is also proved for \eqref{linDif} in~\cite{Za10}. Furthermore, an optimal decay estimate is established for some time-fractional diffusion equations in~\cite{Vergara2015}, where a chain-rule formula (cf.~Proposition \ref{P:chain}) as well as an asymptotic analysis of solutions to time-fractional ODEs play a crucial role (see also~\cite{AV19}). On the other hand, existence of solutions was assumed and is still open as a question in these papers. We shall discuss the existence of solutions in Section~\ref{Ss:Appl}. In~\cite{VeZa17}, the method developed in~\cite{Vergara2015} is also extended to a finite-time blow-up analysis of solutions to time-fractional semilinear diffusion equation. Furthermore, a series of results on the time-fractional (linear) diffusion equation is conveniently summarized in~\cite{Za19}, where an optimal decay estimate (see also~\cite{KSVZ2016}) and an $L^p$-maximal regularity result (see also~\cite{Za05,Za06}) are particularly interesting. 

We close this section by pointing the reader to two monographs on time-fractional evolution equations. A theory based on an integral equation, which is  standard for the case $\alpha=1$, rather than energy methods is well established by Gal and Warma in~\cite{GW20}, where various time-fractional semilinear parabolic problems are comprehensively studied. Furthermore, we also quote~\cite{KRY20} as a nice summary of abstract differential equations with fractional time-derivatives. In particular, a novel Sobolev framework for fractional time-derivatives is developed, and moreover, decay properties of solutions to time-fractional PDEs are discussed.

\subsection{Abstract Theory for Fractional Gradient Flows}
\label{Ss:result}

In this subsection, we consider a fractional gradient flow of subdifferential type in a Hilbert space. To this end, we set up a couple of notions: Let $H$ be a real Hilbert space and let $\phi : H \to (-\infty,+\infty]$ be a proper (i.e., $\phi \not\equiv +\infty$) lower-semicontinuous convex functional with \emph{effective domain} $D(\phi) := \{w \in H \colon \phi(w)<+\infty\}$. Then the subdifferential operator $\partial \phi : H \to 2^H$ (here $2^H$ denotes the power set of $H$) is defined by
\begin{align*}
\partial \phi(w) = \{\xi \in H \colon \phi(v) - \phi(w) \geq (\xi, v-w)_H \ \mbox{ for } v \in D(\phi)\} \ \mbox{ for } w \in D(\phi)
\end{align*}
with domain $D(\partial\phi) := \{w \in D(\phi) \colon \partial \phi(w) \neq \emptyset \}$. The subdifferential is a generalized notion of the Fr\'echet derivative for convex functionals, and it is necessary to establish an abstract theory which is applicable to PDEs having gradient flow structures. The Cauchy problem for a fractional gradient flow reads,
\begin{equation}\label{cp}
-\partial_t^\alpha (u - u_0)(t) \in  \partial \phi(u(t)) \ \mbox{ for } \ t > 0,
\end{equation}
where $u_0 \in D(\phi)$ is a prescribed initial datum. We recall the Riemann-Liouville kernel
$$
k_\beta(t) := \frac{t^{\beta-1}}{\Gamma(\beta)}
$$
of order $\beta \in (0,1)$ and $\partial_t^\alpha (u-u_0) = \fD_{0,t}^{\alpha}(u-u_0) = \partial_t [k_{1-\alpha} * (u-u_0)]$ for $\alpha \in (0,1)$ (see also \eqref{e:FracIntegral}). We are then concerned with strong solutions to \eqref{cp} defined as

\begin{definition}[Strong solution]
\label{D:sol}
For $0 < T < \infty$, a function $u \in L^2(0,T;H)$ is called a strong solution to \eqref{cp} on $[0,T]$ if the following conditions are all satisfied:
\begin{enumerate}
 \item[(i)] $k_{1-\alpha} * (u- u_0) \in W^{1,2}(0,T;H)$ and $[k_{1-\alpha} * (u-u_0)](0) = 0$,
 \item[(ii)] $u(t) \in D(\partial \phi)$ for a.e.~$t \in (0,T)$ and 
$$
- \partial_t^\alpha (u-u_0)(t) \in \partial \phi(u(t))
$$
for a.e.~$t \in (0,T)$.
\end{enumerate}
Moreover, a function $u \in L^2_{\rm loc}([0,\infty);H)$ is called a global strong solution to \eqref{cp} if for every $0 < T < \infty$ the restriction $u|_{[0,T]}$ of $u$ to $[0,T]$ is a strong solution to \eqref{cp} on $[0,T]$.
\end{definition}

Define the space $H^{\alpha,2}_0(0,T;H)$ by
$$
H^{\alpha,2}_0(0,T;H) = \{ w|_{[0,T]} \colon w \in H^{\alpha,2}(\R;H) \ \mbox{ and } \ \mathrm{supp}\,w \subset [0,\infty)\},
$$
where $H^{\alpha,2}(\R;H) \subset L^2(\R;H)$ denotes the Bessel potential function space of order $\alpha$. Then the main result of this section is stated as follows:

\begin{theorem}[{\cite[Theorem 2.3, Corollary 2.7]{A19}}]
\label{T:abst}
For every $u_0 \in D(\phi)$, the Cauchy problem admits a unique global strong solution $u \in L^2_{\rm loc}([0,\infty);H)$ such that the following holds true\/{\rm :}
\begin{enumerate}
 \item[\rm (i)] Set $\mathcal{F}(t) := \left[ k_{1-\alpha} * \left( \phi(u(\cdot)) - \phi(u_0) \right) \right](t)$ for $t \geq 0$. Then the function $t \mapsto \mathcal F(t)$ is continuous on $[0,\infty)$, and it holds that
\begin{align*}
\int^t_s \left\| \partial_t^\alpha (u-u_0)(\tau) \right\|_H^2 \, \d \tau
+ \mathcal{F}(t) \leq \mathcal{F}(s)
\end{align*}
for all $0 \leq s \leq t < \infty$.
 \item[\rm (ii)] It holds that $u \in H^{\alpha,2}_0(0,T;H)$ for $0 < T < \infty$, and moreover, 
       \begin{equation}\label{regu-sol}
	u \in \begin{cases}
	       L^\infty(0,T;H) &\mbox{ if } \ \alpha \in (0,1/2],\\
	       C ([0,T];H) &\mbox{ if } \ \alpha > 1/2
	      \end{cases}
       \end{equation}
       for $0 < T < \infty$. In particular, if $\alpha > 1/2$, then the initial condition $u(0) = u_0$ holds true in a classical sense.
 \item[\rm (iii)] It holds that $\phi(u(\cdot)) \in L^\infty(0,\infty)$ and $k_{\alpha} * \|\partial_t^\alpha (u-u_0)\|_H^2 \in L^\infty(0,\infty)$, and moreover,
\begin{align*}
 \frac 12 \left( k_{\alpha} * \|\partial_t^\alpha (u-u_0)\|_H^2 \right)(t) + \phi(u(t)) \leq \phi(u_0)
\end{align*}
for a.e.~$t > 0$. 
\item[\rm (iv)] Let $u_j$ be a global strong solution to \eqref{cp} with an initial datum $u_{0,j}$ for $j = 1,2$. Then it holds that
\begin{align}\label{contr}
\|u_1(t) - u_2(t)\|_H \leq \|u_{0,1}-u_{0,2}\|_H
\end{align}
for a.e.~$t > 0$.
\end{enumerate}
\end{theorem}

\begin{remark}\label{R:abst}
\begin{enumerate}
 \item[(i)] Even for $u_0 \in \overline{D(\phi)}^H$, one can still construct a function satisfying all the requirements except $k_{1-\alpha} * (u-u_0) \in W^{1,2}(0,T;H)$ for strong solutions to \eqref{cp} in Definition \ref{D:sol}. We can still check $k_{1-\alpha} * (u-u_0) \in W^{1,2}_{\rm loc}((0,T];H)$, which may not, however, imply the uniqueness of such ``solutions'' immediately.
 \item[(ii)] From the theorem above, we know that $\phi(u(t)) \leq \phi(u_0)$ for a.e.~$t > 0$; however, it may still be open to question whether $t \mapsto \phi(u(t))$ is nonincreasing or not. On the other hand, the nonlocal quantity $\mathcal{F}(t)$ is surely continuous and nonincreasing in $t$.
 \item[(iii)] The assertions except (iv) of Theorem~\ref{T:abst} can also be extended to \emph{semiconvex} functionals $\phi : H \to (-\infty,+\infty]$, that is, the case where
$$
\phi(w) + \frac{\lambda}2 \|w\|_H^2
$$
is convex for some $\lambda \in \R$. In particular, $\phi$ may not be convex when $\lambda > 0$. Then the contraction property like \eqref{contr} no longer holds, but a Lipschitz property still holds true, i.e., for each $T > 0$ there exists a constant $C_T > 0$ depending only on $T$, $\alpha$ such that
$$
\|u_1(t) - u_2(t)\|_H \leq C_T\|u_{0,1}-u_{0,2}\|_H
$$
for a.e.~$t \in (0,T)$. See~\cite[Section 5]{A19} for more details.
\end{enumerate}
\end{remark}

We briefly review two key ingredients of the proof for Theorem \ref{T:abst} in~\cite{A19}, that is, an operator theoretic frame for fractional derivatives and a fractional chain-rule inequality. Let $\mathcal{H} = L^2(0,T;H)$ and define an operator $\mathcal{B} : D(\mathcal{B}) \subset \mathcal{H} \to \mathcal{H}$ by
\begin{align*}
 D(\mathcal{B}) &= \{w \in \mathcal{H} \colon k_{1-\alpha} * w \in W^{1,2}(0,T;H), \ (k_{1-\alpha}*w)(0)=0 \},\\
 \mathcal{B}(u) &= \partial_t \left( k_{1-\alpha} * w\right) \ \mbox{ for } \ w \in D(\mathcal{B}).
\end{align*}
Then $\mathcal{B}$ is maximal monotone in $\mathcal{H} \times \mathcal{H}$, and moreover, the Yosida approximation $\mathcal{B}_\lam$ of $\mathcal{B}$ is represented as 
$$
\mathcal{B}_\lambda(w) = 
\frac{\d}{\d t}\left( k_{1-\alpha,\lambda} * w \right) \quad \mbox{ for } \ w \in \mathcal{H},
$$
where $k_{1-\alpha,\lambda}$ is a nonnegative nonincreasing function of class $W^{1,1}_{\rm loc}([0,\infty))$ such that
$$
k_{1-\alpha,\lam} \to k_{1-\alpha} \quad \mbox{ strongly in } L^1(0,T)
\ \mbox{ as } \ \lambda \to 0_+ \ \mbox{ for any } \ T > 0.
$$
We now move on to a fractional chain-rule inequality. Let us start with

\begin{proposition}[{\cite[Proposition 3.4]{A19}}]
\label{P:chain}
Let $T > 0$ and let $h$ be a nonnegative function of class $W^{1,1}(0,T)$. Let $\psi : H \to (-\infty,+\infty]$ be a proper lower-semicontinuous convex functional. Let $u \in L^1(0,T;H)$ be such that $\psi(u(\cdot))$ is integrable over $(0,T)$. Then for any $t \in (0,T)$ satisfying $u(t) \in D(\partial \psi)$ it holds that
\begin{align*}
\MoveEqLeft{
\left( \dfrac{\d}{\d t} [h*(u-u_0)](t), g \right)_H
}\\
&= \dfrac{\d}{\d t} [h*\psi(u(\cdot))](t) + h(t)\left[(u(t)-u_0,g)_H - \psi(u(t))\right]\\
&\quad + \int^t_0 h'(\tau) \left[ (u(t-\tau)-u(t),g)_H + \psi(u(t)) - \psi(u(t-\tau)) \right] \, \d \tau
\end{align*}
for any $u_0 \in H$ and $g \in \partial \psi(u(t))$. In addition, if $u_0 \in D(\psi)$, $h \geq 0$ and $h' \leq 0$, then it holds that
$$
\left( \dfrac{\d}{\d t} [h*(u-u_0)](t), g \right)_H
\geq \dfrac{\d}{\d t} \left[ h* \big( \psi(u(\cdot)) - \psi(u_0) \big) \right](t).
$$
\end{proposition}

The particular choice $\psi(w) = (1/2)\|w\|_H^2$ yields \eqref{FLI} immediately. We usually apply the proposition mentioned above for the regularized kernel $h = k_{1-\alpha,\lambda} \in W^{1,1}(0,T)$ for $T > 0$. Noting that $\mathcal{B}_\lam(w) \to \mathcal{B}(w)$ in $\mathcal{H}$ as $\lambda \to 0_+$ for $w \in D(\mathcal{B})$, which is a well-known property of the Yosida approximation, we find that
\begin{align*}
\MoveEqLeft{
\left( \dfrac{\d}{\d t} [k_{1-\alpha}*(u-u_0)](t), g \right)_H
}\\
 &= \left( \mathcal{B}_\lam(u-u_0)(t), g \right)_H + \left( \mathcal{B}(u-u_0)(t) - \mathcal{B}_\lam(u-u_0)(t), g \right)_H\\
 &\geq \frac{\d}{\d t} \left[ k_{1-\alpha,\lam} * \big( \psi(u(\cdot))-\psi(u_0) \big) \right](t) \\
 &\quad + \left( \mathcal{B}(u-u_0)(t) - \mathcal{B}_\lam(u-u_0)(t), g \right)_H,
\end{align*}
provided that $u - u_0 \in D(\mathcal{B})$. Let $u \in D(\mathcal{B})$ and $\xi \in \mathcal{H}$ be such that $u(t) \in D(\partial \psi)$ and $\xi(t) \in \partial \psi(u(t))$ for a.e.~$t \in (0,T)$. Substitute $g = \xi(t)$ in the above. The integration of both sides over $(s,t)$ yields
\begin{align*}
\MoveEqLeft{
\int^t_s \left( \dfrac{\d}{\d t} [k_{1-\alpha}*(u-u_0)](\tau), \xi(\tau) \right)_H \, \d \tau
}\\
&\geq \left[ k_{1-\alpha,\lam} * \big( \psi(u(\cdot))-\psi(u_0) \big) \right](t)
- \left[ k_{1-\alpha,\lam} * \big( \psi(u(\cdot))-\psi(u_0) \big) \right](s)\\
&\quad - \left\|\mathcal{B}(u-u_0) - \mathcal{B}_\lam(u-u_0)\right\|_{\mathcal H} \|g\|_{\mathcal H}
\end{align*}
for $0 \leq s \leq t \leq T$. Passing to the limit as $\lam \to 0_+$, we obtain the following fractional chain-rule inequality,
\begin{align}
\MoveEqLeft{
\int^t_s \left( \partial_t^\alpha (u-u_0)(\tau), \xi(\tau) \right)_H \, \d \tau
}\nonumber\\
&\geq \left[ k_{1-\alpha} * \big( \psi(u(\cdot))-\psi(u_0) \big) \right](t)
- \left[ k_{1-\alpha} * \big( \psi(u(\cdot))-\psi(u_0) \big) \right](s)\label{cr}
\end{align}
for a.e.~$0 < s \leq t < T$, whenever $u - u_0 \in D(\mathcal{B})$, $\xi \in L^2(0,T;H)$ and $[u(t), \xi(t)] \in G(\partial \psi)$ for a.e.~$t \in (0,T)$. In particular, if $\psi(u) \in L^\infty(0,T)$, then \eqref{cr} holds for all $0 \leq s \leq t \leq T$.

\subsection{Applications to Nonlinear Time-Fractional PDEs}
\label{Ss:Appl}

In this subsection, we discuss applications of Theorem \ref{T:abst} to a couple of nonlinear PDEs. In what follows, we let $\Omega$ be an open subset of $\R^d$ with (smooth) boundary $\partial \Omega$.

\smallskip
\noindent
{\bf (1) Parabolic $p$-Laplace equation.} We consider the following Cauchy-Dirichlet problem:
\begin{alignat}{4}
\partial_t^\alpha (u-u_0) &= \Delta_p u \quad && \mbox{ in } \Omega \times (0,\infty),\label{pde1}\\
u&= 0 &&\mbox{ on } \partial \Omega \times (0,\infty),\label{bc1}
\end{alignat}
where $1 < p < \infty$, $0 < \alpha < 1$ and $\Delta_p$ denotes the so-called $p$-Laplacian given by
$$
\Delta_p u = \mathrm{div}\,\left(|\nabla u|^{p-2}\nabla u\right).
$$
Here we remark that the initial condition is already included in the PDE \eqref{pde1} due to the formulation of strong solutions for time-fractional gradient flows (see also Definition \ref{D:sol}) 
(this will also be applied to the following two other problems). We set $H = L^2(\Omega)$ and $\phi : H \to [0,\infty]$ as follows:
$$
	   \phi(w) = \begin{cases}
		      \frac 1p \int_\Omega |\nabla w|^p\, \d x &\mbox{ if } \ w \in W^{1,p}_0(\Omega),\\
		      \infty &\mbox{ otherwise.}
		     \end{cases}
$$
Then $\phi$ is proper, lower-semicontinuous, and convex in $H$. Hence, thanks to Theorem \ref{T:abst}, the Cauchy-Dirichlet problem \eqref{pde1}--\eqref{bc1} admits a unique global $L^2$-solution $u \in L^\infty(0,\infty;W^{1,p}_0(\Omega)) \cap H^{\alpha,2}_{0,{\rm loc}}([0,\infty);L^2(\Omega))$. We also refer the reader to~\cite{AB23} for the $1$-Laplacian (i.e., $p = 1$) case.

\smallskip
\noindent
{\bf (2) Nonlinear diffusion equation.} Let us consider the following Cauchy-Dirichlet problem:
\begin{alignat}{4}
\partial_t^\alpha (u-u_0) &= \Delta \left( |u|^{m-1}u \right) \quad && \mbox{ in } \Omega \times (0,\infty),\label{pde2}\\
u&= 0 &&\mbox{ on } \partial \Omega \times (0,\infty),\label{bc2}
\end{alignat}
where $0 < m < \infty$ and $0<\alpha<1$. In order to reduce the Cauchy-Dirichlet problem to the abstract gradient flow \eqref{cp}, we set $H = H^{-1}(\Omega)$ and $\phi : H \to [0,\infty]$ as follows:
$$
\phi(w) = \begin{cases}
	   \frac 1{m+1} \int_\Omega |w|^{m+1}\, \d x &\mbox{ if } \ w \in L^{m+1}(\Omega),\\
	   \infty &\mbox{ otherwise,}
	  \end{cases}
$$
which turns out to be a proper lower-semicontinuous convex functional in $H$. Then using Theorem~\ref{T:abst}, the Cauchy-Dirichlet problem \eqref{pde2}--\eqref{bc2} admits a unique global $H^{-1}$-solution $u \in L^\infty(0,\infty;L^{m+1}(\Omega)) \cap H^{\alpha,2}_{0,{\rm loc}}([0,\infty);H^{-1}(\Omega))$ (cf.~see also~\cite{WWZ21}).

\smallskip
\noindent
{\bf (3) Allen-Cahn equation.} The Allen-Cahn equation reads,
\begin{alignat}{4}
 \partial_t^\alpha (u-u_0) &= \Delta u - \beta(u) + \lambda u
 \quad &&\mbox{ in } \Omega \times (0,\infty),\label{pde3}\\
 u&=0 \quad &&\mbox{ on } \partial\Omega \times (0,\infty),\label{bc3}
\end{alignat}
where $0 < \alpha < 1$, $\beta : \R \to \R$ is a maximal monotone function and $\lambda \in \R$. We set $H=L^2(\Omega)$ again and 
$$
\phi(w) = \begin{cases}
	   \frac12\int_\Omega|\nabla u|^2 \, \d x + \int_\Omega \hat \beta(u) \, \d x - \frac \lambda 2 \int_\Omega|u|^2 \, \d x &\mbox{ if } \ u \in D(\phi),\\
	   \infty &\mbox{ otherwise,}
	  \end{cases}
$$
where $\hat \beta$ is a primitive function of $\beta$, that is, $\beta = \partial \hat \beta$ and $D(\phi) = \{u \in H^1_0(\Omega) \colon \hat \beta(u) \in L^1(\Omega)\}$. Then $\phi$ is (not convex but still) semiconvex and lower-semicontinuous in $H$. Hence thanks to Remark~\ref{R:abst},  the Cauchy-Dirichlet problem \eqref{pde3}--\eqref{bc3} admits a unique global $L^2$-solution $u \in L^\infty(0,\infty;H^1_0(\Omega)) \cap H^{\alpha,2}_{0,{\rm loc}}([0,\infty);L^2(\Omega))$.

\subsection{Towards Dynamics of Fractional Gradient Flows}

Recall that due to the nonlocal nature of fractional derivatives in time, one cannot expect a classical semigroup property of the solution operator $S(t) : u_0 \mapsto u(t)$ associated with the fractional gradient flow \eqref{cp}. On the other hand, we may think of a trajectory approach as follows: We denote
$$
I := \left\{ [0,t] \colon t \geq 0 \right\}
$$
and
\begin{alignat*}{4}
\mathfrak{Sol}
= D(\phi) \cup \{ u : [0,t] \to D(\phi) \colon  u \mbox{ is a strong solution of \eqref{cp} on } [0,t]\\
\mbox{for some initial datum } u_0 \in D(\phi) \mbox{ and } t > 0\}.
\end{alignat*}
Moreover, we define an operator $U : I \times D(\phi) \to \mathfrak{Sol}$ by
$$
U([0,t],u_0) := \begin{cases}
		 u_0 &\mbox{ if } \ t = 0,\\
		 u &\mbox{ if } \ t > 0,
		\end{cases}
$$
where $u$ denotes the strong solution of \eqref{cp} on $[0,t]$ with the initial datum $u_0$, and furthermore, we generate a map $\overline{U([0,t],u_0)}$ from $(-\infty,0]$ into $D(\phi)$ by
$$
\overline{U([0,t],u_0)}(\cdot)
:= \begin{cases}
    U([0,t],u_0)(\cdot+t) &\mbox{ if } \ \cdot \in [-t,0],\\
    u_0 &\mbox{ if } \ \cdot \not\in [-t,0]
   \end{cases}
$$
for $t > 0$ and $\overline{U(\{0\},u_0)}(\cdot)\equiv u_0$ for $t = 0$. Define a phase space $\mathfrak{X}$ by
\begin{align*}
\mathfrak{X} = \left\{
\overline{U([0,t],u_0)} \colon u_0 \in D(\phi), \ t \geq 0
\right\}.
\end{align*}
For each $\tau \geq 0$, we define an evolution operator $T_\tau : \mathfrak{X} \to \mathfrak{X}$ by
$$
T_\tau(\overline{U([0,t],u_0)}) = \overline{U([0,t+\tau],u_0)}
$$
for each $\overline{U([0,t],u_0)} \in \mathfrak{X}$. Then we immediately observe that
$$
T_0(\overline{U([0,t],u_0)}) = \overline{U([0,t],u_0)},
$$
and moreover, for $\tau_1,\tau_2 \geq 0$, we find that
\begin{align*}
\MoveEqLeft{
 T_{\tau_2} \circ T_{\tau_1}(\overline{U([0,t],u_0)})
}\\
 &= T_{\tau_2}(\overline{U([0,t+\tau_1],u_0)})\\
 &= \overline{U([0,t+\tau_1+\tau_2],u_0)}
 = T_{\tau_1+\tau_2}(\overline{U([0,t],u_0)}).
\end{align*}
Hence the semigroup property $T_{\tau_1} \circ T_{\tau_2} = T_{\tau_1 + \tau_2}$ holds true. Furthermore, we can also define a functional $\mathcal E : \mathfrak{X} \to \R$ by
\begin{align*}
\mathcal E(\overline{U([0,t],u_0)})
&:= \left[
k_{1-\alpha} * \big( \phi(U([0,t],u_0)(\cdot)) - \phi(u_0)\big)
\right](t)\\
&= \left[
k_{1-\alpha} * \big( \phi(u(\cdot)) - \phi(u_0) \big)
\right](t),
\end{align*}
which is decreasing in $t$ due to Theorem \ref{T:abst}. Therefore we obtain
$$
\mathcal E(\overline{U([0,t],u_0)}) \leq \mathcal E(\overline{U([0,s],u_0)})
$$
if $s \leq t$. That is, $\mathcal E$ is a Lyapunov functional for the dynamical system $(\mathfrak{X},(T_\tau)_{\tau \geq 0})$.

\section{Steady States \& Linear Stability: Space Fractional PDEs}
\label{sec:spectral2}

In Section~\ref{sec:spectral1}, we have presented spectra and bifurcations for dissipative subdiffusion equations. In this section, we focus on one-dimensional space-fractional reaction--diffusion problems. The space-fractional case is technically a bit easier to deal with as one still studies Markovian dynamics and the spectral theory for the space-fractional Laplacian is a bit more classical, or it is even encoded immediately in its definition for the spectral fractional Laplacian (SFL) as introduced in Section~\ref{ssec:defspacefrac}. We proceed to study the effect of the fractional derivative order on the spectra, steady states and bifurcations, also considering three well-known PDEs in which the standard Laplacian is replaced by a spectral fractional Laplacian. A generic fractional reaction--diffusion system can be formulated as
\begin{equation}
    \partial_t u= -D(-\Delta)^{\gamma/2} {u} + F({u};p)=:D\Delta^{\gamma/2} {u} + F({u};p),
    \label{eq:reaction-diffusion-ss-frac}
\end{equation}
where we recall the shorthand notation
\begin{equation*}
\Delta^{\gamma/2}:=-(-\Delta)^{\gamma/2}
\end{equation*}
for the fractional diffusion operator with fractional order~$\gamma\in(0,2)$, considered on a bounded domain~$\Omega\subset \mathbb{R}^d$ and with suitable boundary conditions specified below for each example. As for classical diffusion, a key role in this context is played by the eigenvalues of the fractional diffusion operator. A natural goal, similar to Section~\ref{sec:spectral1}, is to study patterns and we are going to give a brief summary of the spectral properties for different definitions and analytical stability results for bifurcations and Turing instability.

\subsection{Spectra for Space-fractional PDEs}

As discussed in Section~\ref{sec:def}, the several definitions of the fractional Laplacian in~$\mathbb{R}^d$ are all equivalent in that setting.  This is not the case on a bounded domain,  so different definitions also have different spectral properties.  

Regarding its spectral properties, the spectral fractional Laplacian of fractional order $\gamma$ turns out to be the simplest case, as the eigenvalues and the eigenfunctions are the $(\gamma/2)$-power of the eigenvalues of the Laplacian and the very same eigenfunctions of the Laplacian, respectively.  In particular, we recall the spectral definition of the fractional Laplacian~$\Delta^{\gamma/2}$~\cite{capella2010regularity} (cf.~Section~\ref{ssec:defspacefrac})
\begin{equation}
\label{eq:spectral_definition}
\Delta^{\gamma/2} u(x) = -\sum_{j=1}^{\infty} (-\lambda_j)^{\gamma/2} u_j \phi_j(x), \qquad u_j:=\langle u,\phi_j \rangle_{L^2_\Omega}=\int_{\Omega}u \phi_j~\textnormal{d}x
\end{equation}
where~$\lambda_j$ and~$\phi_j$ the eigenvalues and the normalized eigenfunctions, respectively, of the standard Laplacian~$\Delta$ on~$\Omega$, i.e.,~the solutions of the eigenvalue problem~$\Delta \phi_j=\lambda_j\phi_j$ with specific boundary conditions~\cite{abatangelo2019getting}. For the one-dimensional case $\Omega = (0,L)$} the Dirichlet spectral fractional Laplacian has eigenpairs
\begin{equation}
(\phi_j, \tilde{\lambda}_j) = (\sin\left( {j\pi x}/{L} \right), -\left({j \pi}/{L} \right)^{\gamma}),\quad j\geq 1,
 \label{eq:eig_hD}
\end{equation}
while the Neumann spectral fractional Laplacian has eigenpairs 
\begin{equation}
(\phi_j, \tilde{\lambda}_j) = (\cos\left( {j\pi x}/{L} \right), -\left({j \pi}/{L} \right)^{\gamma}), \quad j\geq 0.
 \label{eq:eig_hN}
\end{equation}
For other definitions, as for instance the integral fractional Laplacian,  obtaining results on the spectral properties is not as straightforward as for the spectral definition.  In~\cite{servadei2014spectrum}, the differences between these two operators, namely the spectral and the integral definitions, have been emphasized by studying their spectral properties. In particular, with respect to Dirichlet boundary conditions the first eigenvalue of the integral fractional Laplacian is strictly less than that of the spectral one. We refer also to \cite{bucur2016nonlocal,duo2019comparative,musina2014fractional} for more detailed studies.

\subsection{Analytical Results for Turing Bifurcations}

Pattern formation in reaction--diffusion systems with fractional diffusion has been investigated in different contexts: predator--prey interactions~\cite{owolabi2018modelling},  a semiarid vegetation model~\cite{tian2015turing},  a model for coral reefs~\cite{Somathilake2018}, and chemical systems~\cite{golovin2008turing, zhang2014turing},  just to name a few. As for classical diffusion problems,  also in the fractional setting we expect to obtain conditions, depending on the fractional exponent $\gamma$, that can lead to the appearance of non-homogeneous solutions and how the fractional order influence the possibility of pattern formation (and eventually the type) is an interesting topic. However, the interplay between the fractional order and the other parameters can be nontrivial. In particular it is not obvious if the presence of the fractional operator leads to a common effect in this class of problems or it depends on the particular structure of the equation (nonlinear part). Considering the spectral definition, the modification of the classical Turing instability analysis is quite straightforward, since the fractional order enters as exponent of the eigenvalues in the diffusion matrix (hence in the dispersion relation). With other definitions, the first step is to obtain the corresponding linearized problem. To this end, the Fourier transformed $k-$space is the right framework, since it clarifies the relation between properties of the spectrum and the structure of the solutions; cf.~Section~\ref{sec:spectral1}.\medskip

We now consider three important benchmark problems: the Allen--Cahn equation~\cite{AllenCahn}, the Swift--Hohenberg equation~\cite{CrossHohenberg}, and the Schnakenberg system~\cite{schnakenberg1979simple} in a non-local setting. In particular, we consider the spectral definition.

\textbf{The fractional Allen--Cahn equation.} 
The model with the standard Laplacian, proposed in \cite{AllenCahn}, is a well-studied equation for various polynomial nonlinearities~\cite{alfaro2008singular, alama1997stationary, rabinowitz2003mixed, ward1996metastable}. In particular, the classical Allen--Cahn PDE ($\gamma=2$) with a cubic nonlinearity is also known as the (real) Ginzburg--Landau PDE \cite{aranson2002world}, an amplitude equation or normal form for bifurcations from homogeneous states \cite{mielke2002ginzburg,schneider1996validity}. Moreover, it is also common to consider an additional quintic term added to the Ginzburg--Landau PDE \cite{kapitula1998instability,kuehn2015numerical}. 

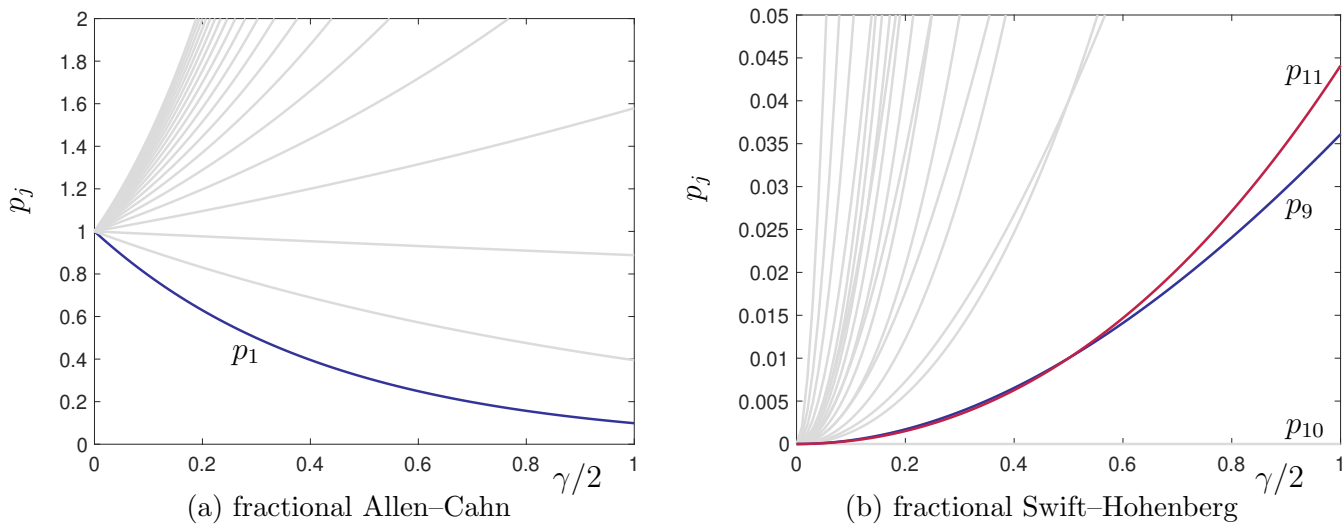
\begin{figure}
 \centering
 \begin{subfigure}[b]{0.43\textwidth}
         \centering
         \begin{overpic}[width=\textwidth]{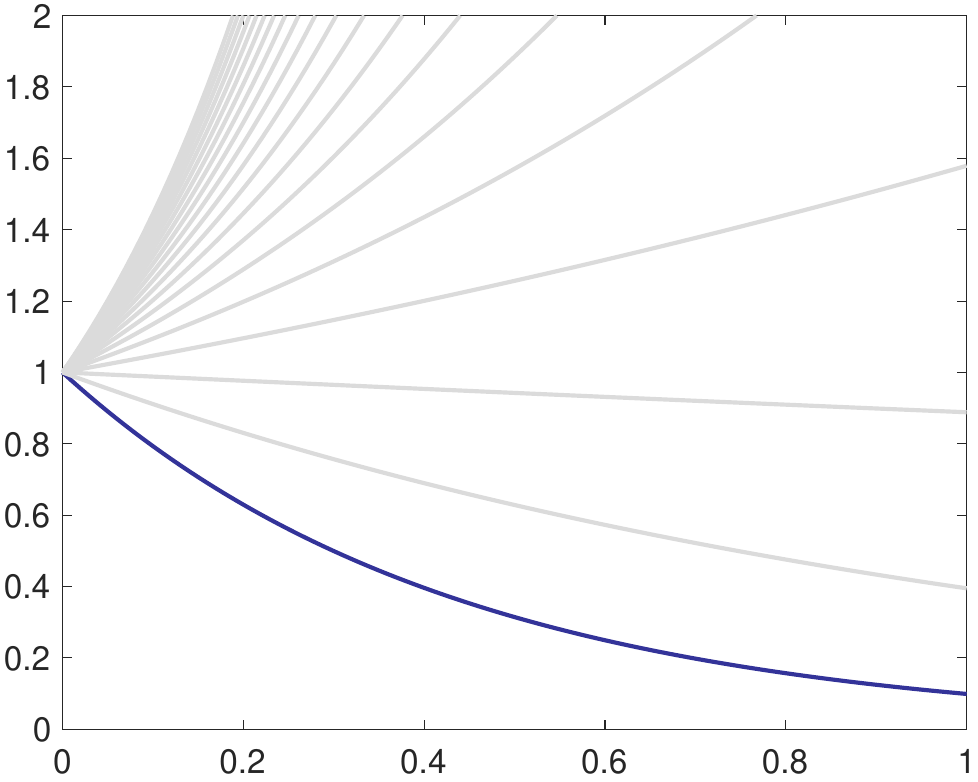}
         \put(-8,45){\rotatebox{90}{$p_j$}}
         \put(85,-2){$\gamma/2$}
         \put(30,20){$p_1$}
         \end{overpic}
         \caption{fractional Allen--Cahn}\label{fig:ac-nsc}
 \end{subfigure}
 \hspace{1cm}
 \begin{subfigure}[b]{0.45\textwidth}
         \centering
         \begin{overpic}[width=\textwidth]{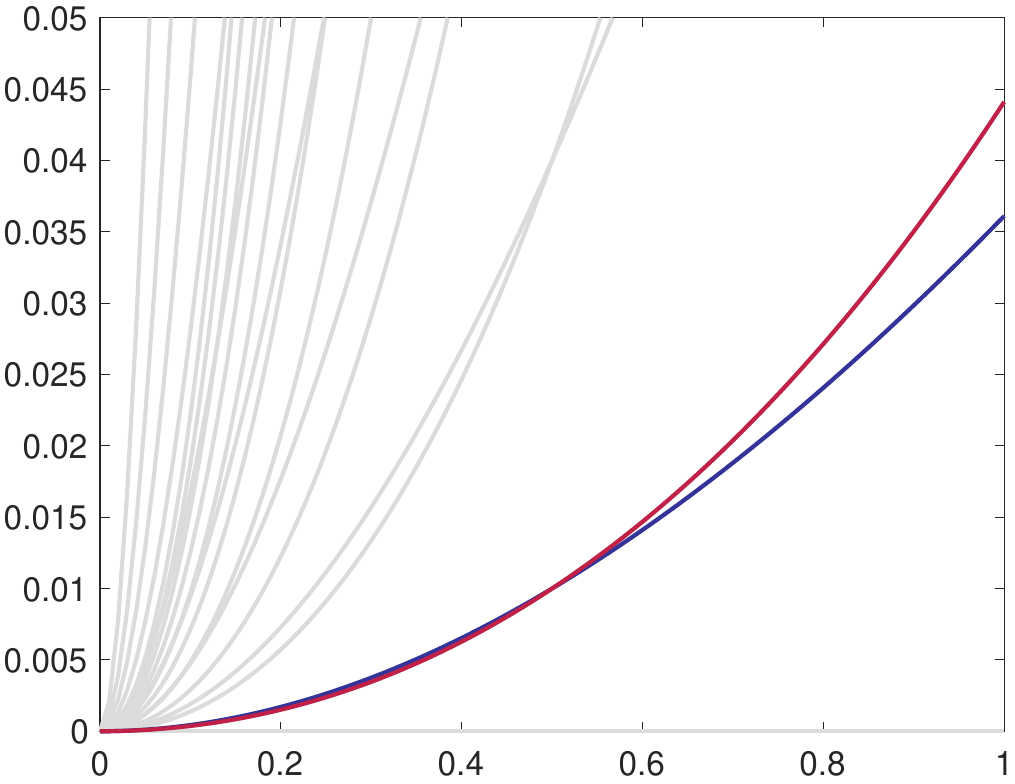}
         \put(-8,45){\rotatebox{90}{$p_j$}}
         \put(85,-2){$\gamma/2$}
         \put(90,43){$p_9$}
         \put(90,65){$p_{11}$}
         \put(90,7){$p_{10}$}
         \end{overpic}
         \caption{fractional Swift--Hohenberg}\label{fig:sh-nsc}
 \end{subfigure} \\[0.5cm]
\caption{Effect of the fractional order $\gamma$ on the bifurcation points $p_j$ on the homogeneous branch for the fractional Allen--Cahn equation (left) and the Swift--Hohenberg equation (right). Each line corresponds to a different wavenumber, using the analytical formula~\eqref{eq:AC-bifurcation-points-s} for Allen--Cahn with $L=10$ and formula~\eqref{eq:subsequent-bif} for Swift--Hohenberg with~$L=10\pi$.}
\label{fig:nsc}
\end{figure}

Different versions of the fractional Allen--Cahn have been recently investigated, see for instance~\cite{AKMR21, akagi2016fractional, bueno2014fourier}. Following~\cite{ehstand2021numerical}, we take with the fractional Allen--Cahn equation with cubic--quintic nonlinearities
\begin{equation}\label{eq:AC-frac}
\partial_t u = \Delta^{\gamma/2} u   + p u + u^3 - q u^5, \quad \quad u\in \mathbb{R},\ \ p,q \in \mathbb{R},
\end{equation}
on a bounded interval $\Omega \subset \mathbb{R}$ of length $L$ with homogeneous Dirichlet boundary conditions and we consider $p$ as the main bifurcation parameter. As for the standard Allen--Cahn equation ($\gamma=2$), the homogeneous solution $u_*=0$ is a stationary steady state of equation~\eqref{eq:AC-frac}. Furthermore, using the Fourier property of the fractional Laplacian, the bifurcation values from the homogeneous states can be analytically computed~\cite{ehstand2021numerical} and are located at
\begin{equation}\label{eq:AC-bifurcation-points-s}
p_j = \left( \frac{j\pi}{L} \right )^{\gamma}, \quad j \geq 1.
\end{equation}
These values depend on the fractional order $\gamma$, and in particular all bifurcation points move towards $p=~1$ as $\gamma$ decreases. The homogeneous steady state $u_*=0$ is stable for $p<p_1$ and unstable otherwise. The destabilization of the homogeneous steady state then happens at $p_1$. Thus, decreasing the fractional order $\gamma$ pushes the stability region of the homogeneous state to greater values of the bifurcation parameter $p$. As a last observation, from equation~\eqref{eq:AC-frac} we can see that the fractional order has no influence on the order of the bifurcation points, as shown in Figure~\ref{fig:ac-nsc}.\\

\textbf{The Fractional Swift--Hohenberg Equation.} 
The classical Swift--Hohenberg equation ($\gamma=2$) \cite{swift1977hydrodynamic} is a widely studied model in pattern dynamics~\cite{CrossHohenberg,Schneider5,Thieleetal}, and also a standard test case for deriving reduced amplitude/modulation equations~\cite{KuehnBook1,SchneiderUecker,KirrmannSchneiderMielke,ColletEckmann1,vanHarten}. 

Recently, the Swift--Hohenberg equation with nonlocal reaction terms has been investigated~\cite{KuehnThrom,MorganDawes,ehstand2021numerical}. Here, as in~\cite{ehstand2021numerical}, we consider the fractional Swift--Hohenberg equation with competing cubic--quintic nonlinearities
\begin{equation}\label{eq:SH-frac}
\partial_t u =  -(1+\Delta^{\gamma/2})^2 u +p u+ q u^3 - u^5, \quad \quad u\in \mathbb{R},\ \ p,q \in \mathbb{R},\ \ q>0,
\end{equation}
on a bounded interval $\Omega \subset \mathbb{R}$ of length $L$ with homogeneous Dirichlet boundary conditions. The parameter $p$ is taken as the main bifurcation parameter. As in the previous model, the homogeneous state $u_*$ is a stationary solution to the fractional equation~\eqref{eq:SH-frac}.  According to~\eqref{eq:eig_hD}, the bifurcations from the homogeneous states occur at\\[-0.2cm]
\begin{equation}\label{eq:subsequent-bif}
p_j = \left(1-\left(\frac{j\pi}{L}\right)^{\gamma}\right)^2, \qquad j\geq 1.
\end{equation}
For simplicity, we fix $L = m\pi$, $m \in \mathbb{N}_{>0}$.  Thus, the homogeneous state $u_*$ is stable for $p < 0 $ and first becomes unstable at $p_c=0$, independently of the fractional order $\gamma$ and with wavenumber $m$. In addition, note that as $\gamma$ tends to zero, the bifurcation points on the homogeneous branch $p_j$ accumulate at $p_c=0$. Interestingly and in contrast to the previous case, as the fractional order decreases, it may happen that bifurcation points shift their position, as shown in Figure~\ref{fig:sh-nsc} for different wavenumbers.\\

\begin{figure}
 \centering
 \begin{subfigure}[b]{0.45\textwidth}
         \centering
         \begin{overpic}[width=\textwidth]{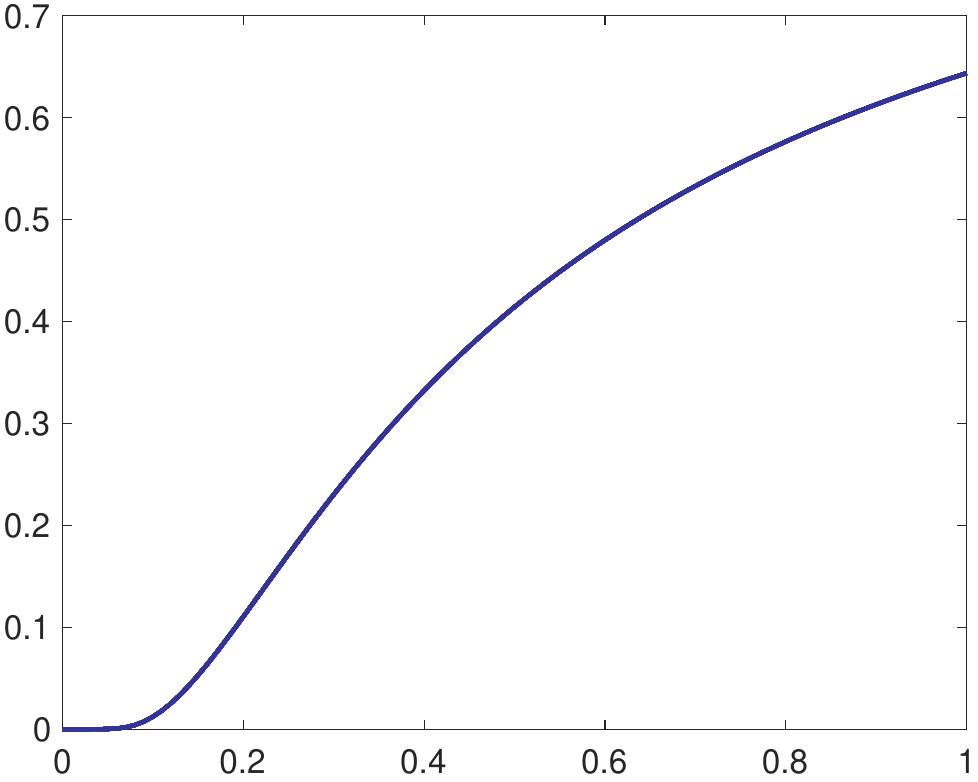}
         \put(-8,45){$k_c$}
         \put(85,-2){$\gamma/2$}
         \end{overpic}
         \caption{Critical value~$k_c(\gamma/2)$ }
         \label{fig:s-k}
 \end{subfigure}
 \hspace{1cm}
 \begin{subfigure}[b]{0.45\textwidth}
         \centering
         \begin{overpic}[width=\textwidth]{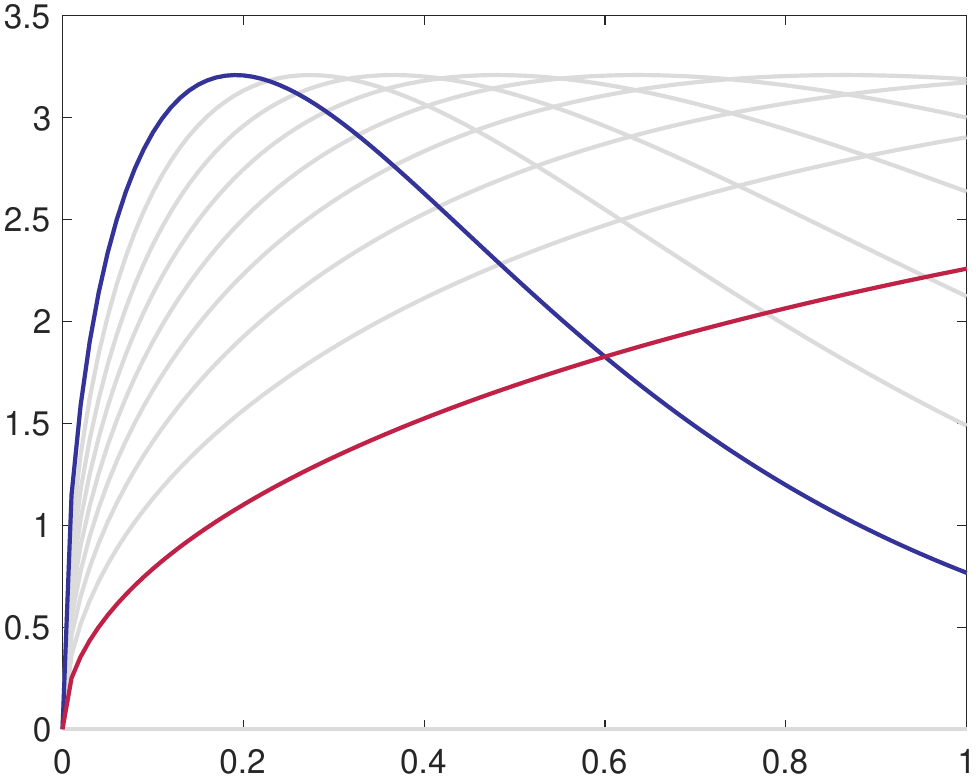}
         \put(-8,45){\rotatebox{90}{$p_j$}}
         \put(85,-2){$\gamma/2$}
         \put(92,18){$p_1$}
         \put(92,55){$p_{9}$}
         \put(92,7){$p_{0}$}
         \end{overpic}
         \caption{Bifurcation values~$p_j(\gamma/2)$ }
         \label{fig:s-nsc}
 \end{subfigure} \\[0.5cm]
\caption{Effect of the fractional order $s$ on critical value~$k_s$ from equation~\eqref{eq:kc-schnak-frac} (left) and on the bifurcation points $\mu_j$ on the homogeneous branch from equation~\eqref{eq:bifurcations-schnak-frac} (right) for the Schnakenberg system. In the left panel, each line corresponds to a different wavenumber, using the analytical formula~\eqref{eq:bifurcations-schnak-frac} with~$L=10\pi$ and $d=60$.}
\label{fig:Schna_nsc}
\end{figure}
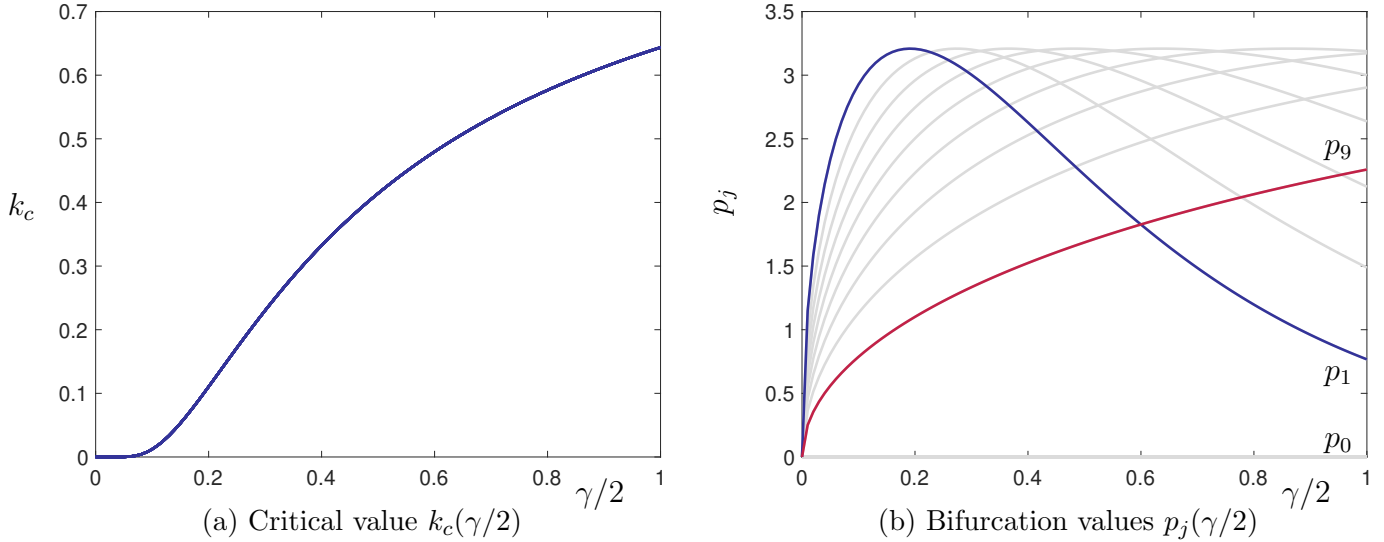

\textbf{The Fractional Schnakenberg System.} 
In its classical formulation, the Schnakenberg system~\cite{schnakenberg1979simple} is one of the prototype reaction--diffusion systems exhibiting Turing patterns~\cite{Murray2}. The fractional Schnakenberg system has been investigated in~\cite{ehstand2021numerical} with a focus on the bifurcation structure and on the snaking branch between periodic branches appearing with additional nonlinearities in the reaction~\cite{Uecker2014,uecker_pde2path_2014}. As in~\cite{ehstand2021numerical}, we focus on the fractional Schnakenberg system with a simple reaction part
\begin{equation}\label{eq:Schnak-frac}
\partial_t\begin{pmatrix}  u_1 \\ u_2 \end{pmatrix} = \begin{pmatrix} 1 & 0 \\ 0 & d \end{pmatrix} \begin{pmatrix} \Delta^{\gamma/2} u_1\\ \Delta^{\gamma/2} u_2 \end{pmatrix} + \begin{pmatrix} -u_1 + u_1^2u_2 \\ p - u_1^2u_2  \end{pmatrix}, \quad \quad  u_1,u_2\in \mathbb{R}, \ \ p, d\in \mathbb{R}, \ \ d>1.
\end{equation}
on a bounded interval $\Omega \in \mathbb{R}$ of length $L$ together with homogeneous Neumann boundary conditions. The usual bifurcation parameter is $p$ but we also want to understand the effect of the fractional order~$\gamma$. As for the Schnakenberg system with standard Laplacian ($\gamma=2$), the spatially homogeneous state $u_*=(p, {1}/p)$ is a stationary solution to~\eqref{eq:Schnak-frac}. Moreover, according to~\eqref{eq:eig_hN}, we can compute analytically the bifurcation values on the homogeneous branch, which are located at
\begin{equation}\label{eq:bifurcations-schnak-frac}
p_j = k_j^{\gamma/2}\sqrt{\frac{d(1-k_j^{\gamma})}{1+k_j^{\gamma}}}, \quad  k_j = \frac{j\pi}{L}, \quad j\geq 0,
\end{equation}
provided that $1-k_j^{\gamma}>0$ so that the associated wavenumber lies in the Turing instability region. Moreover, the critical parameter $p_c$ at which the first instability occurs is
\begin{equation}\label{eq:mu_schnak}
p_c = \sqrt{d(3-\sqrt{8})},
\end{equation}
thus not affected by the fractional diffusion. However, in contrast to the previous  cases, the critical value~$k_c$ at the destabilisation of the homogeneous state depends on the fractional order~$\gamma$
\begin{equation}\label{eq:kc-schnak-frac}
k_c(\gamma/2) = \left(\sqrt{\frac{d-\mu_c^{2}}{2d}}\right)^{\frac{2}{\gamma}}=\left(\sqrt{\sqrt{2}-1}\,\right)^{\frac{2}{\gamma}},
\end{equation}
which is plotted in Figure~\ref{fig:ac-nsc}. Note that we scaled here by a factor of two, as above, to have the horizontal axis of the bifurcation parameter normalized to $[0,1]$. We observe that the wavelength of non-homogeneous stationary solutions bifurcating from the first bifurcation point increases as $\gamma$ decreases. To keep the same wavenumber at the primary bifurcation point, the domain can be ``tuned'' (enlarged) to the primary bifurcation for each fractional order $\gamma$. In Figure~\ref{fig:s-nsc}, we show the effect of the fractional order $\gamma$ on the bifurcation points $p_j$ on the homogeneous branch. As for the Swift--Hohenberg equation, two bifurcation points corresponding to $j=1$ and $j=9$ shift their position as $\gamma$ decreases. Note that additional nonlinearties as in~\cite{Uecker2014,uecker_pde2path_2014} do not affect~$p_c$ and~$k_c$. 

\section{Traveling Waves: Space Fractional PDEs}
\label{sec:tw}

For earlier reviews on front propagation in reaction-diffusion equations with anomalous diffusion (including subdiffusive systems), see e.g.~\cite{VolNecNep13} and the references therein.

\begin{remark}
Subdiffusion limited reaction-diffusion equations like~\eqref{e:Zacher}, as any typical non-autonomous equation, do not support (non-trivial) traveling wave solutions of the form~\eqref{eq:AnsatzTW} with nonzero speed.
For example, consider the subdiffusion limited reaction-diffusion equation in~\eqref{e:Zacher} with \emph{(Djrbashian-)Caputo fractional derivative}~$\tilde\fD_{0,t}^{\alpha}$ of order~$\alpha\in(0,1)$.
Computing the (Djrbashian-)Caputo fractional derivative with respect to $t\geq 0$ of the traveling wave ansatz~$u(t,x)=u(x-st)=u(\zeta)$ with $\zeta=x-st$ and $s<0$ yields
	\begin{align*}
		\tilde\fD_{0,t}^{\alpha} \big[u(x-st)\big]
		&= \frac{1}{\Gamma(1-\alpha)} \int_0^t (t-\sigma)^{-\alpha} \big(\partial_\sigma u(x-s\sigma)\big)\, \dif \sigma \\
		&= \frac{1}{\Gamma(1-\alpha)} \int_0^t (t-\sigma)^{-\alpha} \big(-s u'(x-s\sigma)\big)\, \dif \sigma
	\intertext{and using $\xi:=x-s\sigma$, $\zeta=x-st$, $\delta:=\zeta-\xi$,}
		&= \frac{(-s)^\alpha}{\Gamma(1-\alpha)} \int_x^{x-st} \frac{u'(\xi)}{(x-st-\xi)^{\alpha}} \, \dif \xi
		= \frac{(-s)^\alpha}{\Gamma(1-\alpha)} \int_x^\zeta \frac{u'(\xi)}{(\zeta-\xi)^{\alpha}} \, \dif \xi \\
        &= \frac{(-s)^\alpha}{\Gamma(1-\alpha)} \int_0^{-st} \frac{u'(\zeta-\delta)}{\delta^{\alpha}} \, \dif \delta .
	\end{align*} 
The integral is rewritten using the traveling wave variables $\zeta$ and $\xi$, but still depends either on the spatial variable~$x$ or on the time variable~$t$.
Therefore, a subdiffusion limited reaction-diffusion equation~\eqref{e:Zacher} with Caputo fractional derivative can only support traveling wave solutions of the form~$u(t,x)=u(x-st)$ that are either stationary ($s=0$) or steady states (if $s\ne0$). 
However, in the natural long-term memory limit (for fixed $\zeta$), 
\begin{align*}
 \lim_{t\to\infty} \tilde\fD_{0,t}^{\alpha} \big[u(x-st)\big]
&= \lim_{t\to\infty} \frac{(-s)^\alpha}{\Gamma(1-\alpha)} \int_0^{-st} \frac{u'(\zeta-\delta)}{\delta^{\alpha}} \, \dif \delta
= \frac{(-s)^\alpha}{\Gamma(1-\alpha)} \int_0^{\infty} \frac{u'(\zeta-\delta)}{\delta^{\alpha}} \, \dif \delta \\
&= \frac{(-s)^\alpha}{\Gamma(1-\alpha)} \int_{-\infty}^{\zeta} \frac{u'(\delta)}{(\zeta-\delta)^{\alpha}} \, \dif \delta ,
\end{align*}
the traveling wave equation of~\eqref{e:Zacher} turns into an integro-differential equation of the form,
\begin{equation} \label{e:TW_Caputo}
  \frac{(-s)^\alpha}{\Gamma(1-\alpha)} \int_{-\infty}^{\zeta} \frac{u'(\delta)}{(\zeta-\delta)^{\alpha}} \, \dif \delta 
= \sigma\partial_\zeta^2 u + f(u), \qquad \zeta\in\R,  
\end{equation}
cf.\ \cite{NecVolNep10,VolNecNep13}. For the subdiffusive Allen-Cahn equation, a formal analysis of bistable fronts in \cite{NecVolNep10} predicts algebraic decay in one of the tails due to memory effects. For the time-fractional Fisher--KPP equation~\eqref{e:Zacher} with Caputo derivative~$\tilde\fD_{0,t}^{\alpha}$, the existence of invading (i.e., monostable) fronts satisfying~\eqref{e:TW_Caputo} has been proven in~\cite{Ishii23_arxiv} (for complementary results see~\cite{AchCueHit14,CueAch17}). Some numerical simulations (see e.g., \cite[\S2]{Ishii23_arxiv}) suggest that solutions of the initial value problem for~\eqref{e:Zacher} with initial datum given by a Heaviside step function converge to traveling wave solutions of~\eqref{e:TW_Caputo}.  
However, we are not aware of an analytical proof of the convergence of solutions of~\eqref{e:Zacher} to traveling wave fronts of~\eqref{e:TW_Caputo} for $\alpha\in(0,1)$. 
\end{remark}

\medskip
In the remaining part of this section, we are going to focus on the existence of traveling wave solutions for reaction-diffusion equations with anomalous diffusion of the form
\be \label{RD_Levy}
 \diff{u}{t} 
=\Levy u + f(u) \,,\quad x\in\R^d \,,\quad t\in(0,\infty) \,,
\ee
where $f\in C^\infty(\R)$ is a nonlinear function and $\Levy$ is a L\'evy operator (i.e.,~the generator of a L\'evy process). 

\subsection{Traveling Wave Solutions in One-Dimensional Settings}

A particularly important subclass arises in the Riesz-Feller case leading to evolution equations
 \begin{equation} \label{RD_RieszFeller}
  \diff{u}{t} = \RieszFeller u + f(u) \,,\quad x\in\R \,,\quad t\in(0,\infty) \,,
 \end{equation}
 where $f\in C^\infty(\R)$ is a nonlinear function of bistable type, {i.e.}, $f$ has precisely three roots $u_1 < u_2 < u_3$ in the interval $[u_1,u_3]$ such that 
\begin{subequations} \label{f_bistable}
 \begin{align} 
  & f(u_1)=f(u_2)=f(u_3)=0\,, \quad f'(u_1)<0\,, \quad f'(u_3)<0 , 
\\
  & f<0 \text{ on } (u_1,u_2) \text{ and } f>0 \text{ on } (u_2,u_3) ,
 \end{align}
\end{subequations} 
 and we recall that $\RieszFeller$ denotes a Riesz--Feller operator for some fixed parameters $1<\gamma\leq 2$ and $\abs{\theta} \leq \min\{\gamma,2-\gamma\}$; see Section~\ref{ssec:defspacefrac}. To study the traveling wave problem, it is necessary to extend the nonlocal operator to $C^2_b(\R)$. The following integral representations may be used to accomplish this task:

\begin{theorem} 
\label{thm:RieszFeller:extension}
If $0<\gamma<1$ or $1<\gamma<2$ and $|\theta|\leq \min\{\gamma,2-\gamma\}$,
 then for all $f\in\SchwartzTF(\R)$ and $x\in\R$ 
\begin{align} 
\label{eq:RieszFeller1}
\nonumber \RieszFeller f(x) = \tfrac{c_1 - c_2}{1-\gamma} f'(x) &+
  c_1 \integrall{0}{\infty}{ \tfrac{f(x+y)-f(x)-f'(x)\,y {\bf 1}_{(-1,1)}(y)}{y^{1+\gamma}} }{y} \\
  &+ c_2 \integrall{0}{\infty}{ \tfrac{f(x-y)-f(x)+f'(x)\,y {\bf 1}_{(-1,1)}(y)}{y^{1+\gamma}} }{y}
\end{align}
where ${\bf 1}_{(-1,1)}(\cdot)$ is an indicator function and some constants $c_1, c_2 \geq 0$ with $c_1+c_2 >0$.
\end{theorem}

The result follows from~\cite[Theorem 31.7]{Sato:1999}; see also~\cite[Theorem~2.4]{AK15AC}. A traveling front solution of~\eqref{RD_RieszFeller} is a traveling wave solution~\eqref{eq:AnsatzTW} connecting distinct endstates~$\upm\in\{u_1,u_3\}$ in the sense that $\lim_{\zeta\to\pm\infty} \profile(\zeta) = \upm$. The profile~$u=\profile(\zeta)$ of a traveling wave solution satisfies (the traveling wave equation)
\be 
\label{TWE:RieszFeller}
 -\speed \profile'(\zeta) = \RieszFeller{}[\profile](\zeta) + f(\profile), \qquad
 \zeta\in\R,
\ee
where $\RieszFeller$ has to be understood as its extension~\eqref{eq:RieszFeller1} to $C^2_b$-functions. The existence of traveling front solutions of~\eqref{RD_RieszFeller} has been proved in case of $\Riesz[2]=\partial_x^2$~\cite{Aronson+Weinberger:1974,Fife+McLeod:1977}, fractional Laplacians~$\Riesz = -(-\partial^2_x)^{\gamma/2}$, $\gamma\in(0,2)$ modeling symmetric super-diffusion~\cite{Cabre+Sire:2015,Chmaj:2013,Gui+Zhao:2014,Palatucci+etal:2013}, and Riesz--Feller operators $\RieszFeller$ with $1<\gamma<2$ and $\abs{\theta} \leq \min\{\gamma,2-\gamma\}$ in~\cite{AK15AC}. If a traveling front solution of~\eqref{eq:AC} with profile $\profile$ exists then its wave speed $\speed$ satisfies
\be
\label{eq:WS} 
 \speed
= -\tfrac{\integrall{\um}{\up}{f(w)}{w}}{\integrall{-\infty}{\infty}{ (\profile'(x))^2 }{x} }\,.
\ee 
Thus the potential $F(u):=F(\um)+\integrall{\um}{u}{f(v)}{v}$ indicates which stable state, either $\um$ or $\up$, will replace the other one. In case of a balanced potential, i.e., $F(\um)=F(\up)$ or $\integrall{\um}{\up}{f(v)}{v}=0$, a stationary traveling front will exist, {i.e.}, both stable states will co-exist. In contrast, in case of an unbalanced potential, $\integrall{\um}{\up}{f(v)}{v}\ne 0$, the stable state with smaller potential value will replace the one with larger potential value, also called the metastable state. It is important to note that many of the following results are restricted to balanced bistable functions $f$. For bistable reaction-diffusion equations~\eqref{RD_RieszFeller} with fractional Laplacian~$\Riesz$ where $0<\gamma<2$, 
\be
\label{RD_fracLaplacian}
  \diff{u}{t} = \Riesz u + f(u) \,,\quad 
  x\in\R \,,\quad t\in(0,\infty) \,,
\ee
the first preliminary investigations of traveling wave solutions highlighted some similarities and differences to the classical case~\eqref{eq:AC}, see~\cite{Nec+etal:2008,Volpert+etal:2010,Zanette:1997}: Using the simple reaction function~$f$ of McKean~\cite{McKean:1970}, traveling wave solutions of~\eqref{RD_fracLaplacian} can be constructed explicitly~\cite{Zanette:1997}. This model may serve as a benchmark example. For~\eqref{RD_fracLaplacian} with $\gamma\in (1,2)$ and bistable reaction function $f(u)=u(1-u^2)$, a variational formulation can be derived such that~\eqref{RD_fracLaplacian} is the associated Euler-Lagrange equation~\cite{Nec+etal:2008}. For general reaction functions, traveling front solutions~\eqref{eq:AnsatzTW} with profile~$\profile$ and wave speed~$\speed$ satisfy
 \begin{equation} 
 \label{eq:WS_fracLaplacian} 
-\speed\integrall{-\infty}{\infty}{ (\profile'(\zeta))^2 }{\zeta} 
=\integrall{-\infty}{\infty}{\Riesz{} [\phi](\zeta)\ u'(\zeta)}{\zeta} 
 +\integrall{\um}{\up}{f(w)}{w}  .
 \end{equation} 
For monotone profiles~$u$, the first term on the right-hand side vanishes such that the wave speed~$\speed$ satisfies again~\eqref{eq:WS}, see~\cite{Volpert+etal:2010}. However, the tail behavior of profiles is different: traveling wave profiles~$u$ approach the endstates at an algebraic rate $1/\abs{\zeta}^\gamma$ for $1\leq\gamma<2$, in contrast, to an exponential rate for $\gamma=2$, see~\cite{Volpert+etal:2010}. The existence and uniqueness of traveling front solutions of fractional Allen--Cahn equation~\eqref{RD_fracLaplacian} was proven via different methods~\cite{Cabre+Sire:2015,Chmaj:2013,Palatucci+etal:2013,Gui+Zhao:2014} (which we discuss below). The stability of traveling front solutions has been proven via the squeezing method using sub- and supersolutions~\cite{AK15AC,MaNiWa19}, and via spectral analysis~\cite{HoYu22}. Next, we discuss the different approaches to prove the existence of traveling front solutions.

\subsubsection{Approximation of the Diffusion Operator}

Fractional Laplacians~$\Riesz$ can be approximated by a family of convolution operators $J_\epsilon \ast \profile -\profile$ with $J_\epsilon\in L^1(\R)$ such that $J_\epsilon \ast \profile -\profile \longrightarrow \Riesz \profile$ for $\epsilon\to 0$ in a suitable sense. The associated (traveling wave) equations $-\speed\profile' = J_\epsilon \ast \profile -\profile + f(\profile)$ exhibit for all sufficiently small $\epsilon>0$ a unique monotone solution $(\profile_\epsilon,\speed_\epsilon)$ with $\profile_\epsilon'>0$ see~\cite{Bates+etal:1997, Chen:1997}. The limit $\profile = \lim_{\epsilon\to 0} \profile_\epsilon$ exists and yields the profile~$\profile$ of a traveling front solution of~\eqref{RD_fracLaplacian}.

\subsubsection{Using an Extension Problem}
\label{sec:FA-extension-problem}

The operator square root of the Laplacian $\Riesz[1/2]$ can be obtained from the harmonic extension problem to the upper half space as the operator that maps the Dirichlet boundary condition to the Neumann condition. For fractional Laplacians and other integro-differential operators, a similar characterization is available~\cite{caffarelli-silvestre07,Kw22} (see also Section~\ref{sec:mel-IFL} and the related discussion of the SFL in (\ref{eq:mel-CST})). From those characterizations, properties of these integro-differential equations can be derived from purely local arguments in the
extension problems~\cite{caffarelli-silvestre07}. For example, the stationary traveling wave problem 
\be 
\label{eq:RD:Cabre} 
 0 = \Riesz w + f(w) \Xx{in} \R 
\ee
for $\gamma\in(0,2)$ and a bistable function~$f$ can be related to a nonlinear boundary value problem for a partial differential equation. This allows one to prove that stationary traveling front solutions~$w(x)=\profile(x)$ of~\eqref{eq:RD:Cabre} exist if and only if $f$ is bistable with a balanced potential; and are unique (up to translations) under the additional assumption~$f'(\upm)<0$, see~\cite[Theorem 2.4]{Cabre+Sire:2015}. Moreover, the asymptotic behavior of the profile of traveling front solutions can be analyzed, see~\cite[Theorem 2.7]{Cabre+Sire:2015}.

\subsubsection{Minimizers of a Nonlocal Energy Functional}

Consider the nonlocal energy functional 
\be \label{functional:Palatucci}
 \cE(w,\Omega) 
:= \cK(w,\Omega) + \integral{\Omega}{F(w(x))}{x}
\ee
where $\cK(w,\Omega)$ is (related to) the squared Gagliardo--Sobolev--Slobodeckij $H^s$ semi-norm of $w$, and $F$ is a double-well potential with $F(\upm)=0$. In the one-dimensional setting, stationary traveling front solutions~$w(x)=\profile(x)$ of~\eqref{eq:RD:Cabre} are local minimizers of the functional $\cE(w,\R)$. For a bistable function $f\in C^1(\R)$ with balanced potential, there exists a unique (up to translations) nontrivial global minimizer~$w$ of~\eqref{functional:Palatucci}. This minimizer~$w$ is a traveling front solution of~\eqref{eq:RD:Cabre}, and is unique (up to translations) in the class of monotone solutions of~\eqref{eq:RD:Cabre}, see~\cite[Theorem 2]{Palatucci+etal:2013}.

\subsubsection{Homotopy Arguments}

In~\cite{Cabre+Sire:2015,Palatucci+etal:2013}, the existence of \emph{stationary} traveling front solutions of~\eqref{RD_fracLaplacian} (i.e.) with \emph{balanced} potentials have been shown. For general bistable reaction functions~$f$ (with possibly unbalanced potential), the existence of (unique) traveling front solutions of~\eqref{RD_fracLaplacian} can be proven via a homotopy from the balanced case~\cite{Gui+Zhao:2014}.

\subsubsection{Semigroup Approach and Squeezing Method using Sub- and Super-solutions}

For reaction-diffusion equations supporting a comparison principle, the squeezing method using sub- and supersolutions allows one to prove the existence, uniqueness and local asymptotic stability with exponential decay of traveling wave solutions. The squeezing method was formulated---first---for a class of nonlinear (nonlocal) evolution equations~\cite{Chen:1997} including~\eqref{eq:AC} and
\be \label{eq:RD:Bates}
 \diff{u}{t} = J \ast u - u + f(u) = \integral{\R}{J(x-y)\ u(y,t)}{y} - u(t,x) + f(u(t,x)) \,, \quad x\in\R\,, \quad t>0 \,,
\ee
for non-negative functions $J\in L^1(\R)$ and general bistable functions $f\in C^2(\R)$; and---later---extended to other classes of reaction-diffusion equations, such as~\eqref{RD_RieszFeller} see~\cite{AK15AC}.

\begin{remark}[Metastability]
The Allen--Cahn equation~\eqref{eq:AC} exhibits metastable behavior, i.e., some solutions appear to be stationary on extremely long---yet still transient---time scales, see e.g., \cite{We21} and references therein. The metastable behavior of the fractional Allen--Cahn equation~\eqref{RD_fracLaplacian} has been studied via asymptotic analysis and numerical simulations, see~\cite{AKMR21}.
\end{remark}

For an earlier review on numerical methods to construct traveling front solutions for the fractional Allen--Cahn equation~\eqref{RD_fracLaplacian}, see e.g.~\cite{AK15AC2}.

\subsection{Traveling Wave Solutions in Multi-Dimensional Settings}

Next, we discuss the existence of traveling wave solutions for reaction-diffusion equations~\eqref{RD_Levy} with bistable function~$f$ in the multi-dimensional setting. Consider the classical bistable reaction-diffusion equation
\be \label{RD_Laplacian}
 \diff{u}{t} 
=\Delta u + f(u) \,,\quad x\in\R^d \,,\quad t\in(0,\infty) \,,
\ee
and bistable reaction function~$f$ satisfying~\eqref{f_bistable}. Here, we briefly discuss the theory of traveling front solutions to~\eqref{RD_Laplacian} and refer to~\cite{Xin00,Ha16,Ta21} for further information. Planar traveling waves in $\R^d$ are of the form
\be 
\label{TWS:planar}
 u(t,x) :=\profile(x\cdot e-\speed t)
\ee
for some vector~$e\in\R^d$ with $\|e\|=1$ and $\speed\in\R$. If a planar traveling wave~\eqref{TWS:planar} is a solution of~\eqref{RD_Laplacian} then its profile~$\profile$ satisfies a second-order ODE 
\be 
\label{TWE:planar}
 \profile'' +\speed\profile'(\zeta) +f(\profile(\zeta)) =0
\qquad \text{for } \zeta\in\R ,
\ee
just like in the one-dimensional setting. A planar traveling front is a planar traveling wave~\eqref{TWS:planar} which satisfies~$\lim_{\zeta\to\pm\infty} \profile(\zeta) =\upm$ for some $\upm\in\R$. The existence of planar traveling fronts of~\eqref{RD_Laplacian} and the properties of its profile~$\profile$ and speed~$\speed$ follow again from the analysis of ODE~\eqref{TWE:planar}, see~\cite{Aronson+Weinberger:1974,Fife+McLeod:1977}. For bistable reaction function~$f$ satisfying~\eqref{f_bistable}, $\um=u_3$ and $\up=u_1$, there exists a monotone profile $\profile\in\C^1(\R)$ and a unique speed~$\speed$ satisfying~\eqref{eq:WS}. The level sets of a planar traveling front at level~$c\in\R$ and time~$t\geq0$ are
\begin{equation} 
\label{spTF:levelSet}
 \{x\in\Rd:\ u(t,x)=c\} 
=\{x\in\Rd:\ x\cdot e =\phi^{-1}(c) +\speed t \} ,
\end{equation}
i.e., the level sets are parallel hyperplanes which are orthogonal to the direction of propagation~$e$. Hence, they are called planar. Moreover, planar traveling fronts are only unique up to shifts. In particular, these traveling fronts are invariant in the moving coordinate frame with speed~$\speed$ in the direction~$e$. In the one-dimensional setting, bistable traveling fronts are asymptotically stable. In the multi-dimensional setting, the stability w.r.t. perturbations of the initial data is more subtle, but see, e.g., \cite{Xin92,LevXin92,Jon83,Kap97,Ta21,Volpert+etal:1994}. Robustness/stability of traveling front solutions w.r.t. perturbations of the nonlinear bistable reaction function has been studied also~\cite{WaTa23}.

For classical bistable reaction-diffusion equations~\eqref{RD_Laplacian} in the multi-dimensional setting, other kinds of traveling front solutions exist, see the reviews in e.g.~\cite{Ha16,Ta21}: For example, traveling fronts with non-planar level sets exist such as (i) (axisymmetric) conical-shape, (ii) (non-axisymmetric) pyramidal shape, and (iii) exponential/parabolic shape~\cite{CheGuoHamNinRoq07}, to name a few examples. Again, these traveling wave solutions (resp. their profiles) are invariant in a moving coordinate frame with constant speed.

The notion of transition fronts has been introduced to encompass all traveling front solutions which converge to $\upm$ far away from their level sets uniformly in time~\cite{BeHa12,Ha16}. There exist transition fronts of the classical bistable reaction-diffusion equation~\eqref{RD_Laplacian} which are not invariant in any moving coordinate frame with constant speed~\cite[Theorem 2.9]{Ha16}. Thus, traveling fronts are a (strict) subset of transition fronts. Nonetheless, for (bistable) transition fronts a~\emph{global mean speed} can be defined which is independent of the shape of level sets of the transition fronts, see~\cite{Ha16}. 

\begin{remark}[Rigidity results for the Allen-Cahn equation]
In case of bistable reaction functions~$f$ with balanced potential~$F$, Gibbons' conjecture may be stated as follows \cite{FarVal11,DipVal19}: Let $u$ be a bounded solution of the stationary Allen-Cahn equation
\be 
\label{RD_Laplacian_stationary}
 0 =\Delta u +u -u^3 ,\quad x\in\R^d ,
\ee
such that, say w.r.t.~$x_d$,
\be \label{Gibbons:uniformLimits}
 \lim_{x_d\to\pm\infty} u(x',x_d) =\upm \qquad
 \text{uniformly w.r.t. } x':=(x_1,\ldots,x_{d-1})\in\R^{d-1}
\ee
then $u(x)=\phi(x_d)$ for some $\phi:\R\to\R$. Gibbons' conjecture for~\eqref{RD_Laplacian_stationary} is true for any $d\in\N$, see the references given in~\cite{FarVal11,DipVal19}. De Giorgi's conjecture about the rigidity of the Allen-Cahn equation~\eqref{RD_Laplacian_stationary} assumes the monotonicity of the bounded solution~$u$ in one direction, say~$\partial u/\partial x_d >0$; instead of the uniform limits in~\eqref{Gibbons:uniformLimits}. Variants of De Giorgi's conjecture have been proven for $d\leq 8$, whereas for $d\geq 9$ counterexamples have been constructed, see the references given  in~\cite{CheGuoHamNinRoq07,FarVal09,Sav10,ChaWei18}. The analogous conjecture for the parabolic equation~\eqref{RD_Laplacian} does not hold, due to the existence of traveling wave solutions with exponential/parabolic shaped level sets~\cite{CheGuoHamNinRoq07}.
\end{remark}

Consider the fractional bistable reaction-diffusion equation
\be 
\label{RD_fracLaplacian_multiD}
 \diff{u}{t} 
=-(-\Delta)^{\gamma/2} u + f(u) \,,\quad x\in\R^d \,,\quad t\in(0,\infty) \,,
\ee
with fractional Laplacians $(-\Delta_x)^{\gamma/2}$, $\gamma\in(0,2)$ in $\R^d$. If a planar traveling wave~\eqref{TWS:planar} with $\phi\in C^2_b(\R)$ is a solution of~\eqref{RD_fracLaplacian_multiD} then its profile~$\phi$ satisfies a one-dimensional traveling wave equation~\eqref{TWE:RieszFeller} with $\theta=0$, i.e.
\be 
\label{TWE:fracLaplacian}
 -\speed \profile'(\zeta)
=-(-\partial_\zeta^2)^{\gamma/2}[\profile](\zeta) + f(\profile), \qquad
 \zeta\in\R.
\ee
For~\eqref{TWS:planar} with $\phi\in C^2_b(\R)$, the traveling wave equation~\eqref{TWE:fracLaplacian} is derived by using the identity
\begin{align*}
& -(-\Delta)^{\gamma/2} \big[\phi(x\cdot e -st)\big]
\\
&=c_{d,\gamma}  \integral{\R^d}{\frac{\phi((x+y)\cdot e -st)+\phi((x-y)\cdot e-st)-2\phi(x\cdot e-st)}{2 |y|^{d+\gamma}}}{y}
\\
&=c_{1,\gamma} \Bigg( \integral{\R}{\frac{\phi(\zeta +y_e)+\phi(\zeta-y_e)-2\phi(\zeta)}{2 |y_e|^{1+\gamma}}}{y_e} \Bigg)\Bigg|_{\zeta=x\cdot e-st}
\\
&=-(-\partial_\zeta^2)^{\gamma/2}[\profile](\zeta) ,
\end{align*}
where $y_e := y\cdot e$, see also~\cite{GraRyz15}. Thus, if there exist traveling front solutions of the fractional Allen-Cahn equation~\eqref{RD_fracLaplacian} in the one-dimensional setting, then these solutions of~\eqref{RD_fracLaplacian}---with profile~$\profile$ and speed~$\speed$---can be used to construct families of planar traveling front solutions of~\eqref{RD_fracLaplacian_multiD} in the multi-dimensional setting: For $e\in\R^d$ with $\|e\|=1$, the function~$u(x):=\phi(x\cdot e-st)$ is a planar traveling wave front solution of~\eqref{RD_fracLaplacian_multiD}.

For the fractional Allen--Cahn equation~\eqref{RD_fracLaplacian_multiD} with $\gamma\in(1,2)$, the existence (and stability) of~\emph{non-planar} traveling wave solutions has been proven in~\cite{ChaWei17,MaWan21,MaNiuWan22}. Stability of steady states for the fractional bistable reaction-diffusion equation has been studied in~\cite{CheYua19}.

\begin{remark}[Rigidity results for the fractional Allen-Cahn equation]
Following Gibbons and De Giorgi, analogous conjectures for
\be 
\label{RD_fracLaplacian_stationary}
 0 =-(-\Delta)^{\gamma/2} u +u -u^3 ,\quad x\in\R^d ,
\ee
have been studied: The fractional Gibbons' conjecture holds for any $d\in\N$, see~\cite{FarVal11} and references therein. The fractional De Giorgi's conjecture, i.e., De Giorgi's conjecture for~\eqref{RD_fracLaplacian_stationary}, has been partially proven, see a review of results in~\cite{ChaWei18,DipVal19,DipSerVal20,DipVal23}. Notably, the parameter regimes for $\gamma\in(0,1)$ and $\gamma\in[1,2)$ are very different, see the discussion in~\cite{DipVal19,DipSerVal20}.
\end{remark}


\section{Discretization and Simulation: Fractional PDEs}
\label{sec:discrete}

Eventually, the analytical/pen-and-paper methods discussed in the previous sections will reach their limits for fractional dissipative PDEs similarly to the case of standard differential operators. Hence, we have to develop numerical schemes to understand in more detail the dynamics we are interested in. Two aspects make the numerical treatment of fractional (partial) differential equations more challenging than that of standard (partial) differential equations. First, fractional differential operators are inherently nonlocal operators which brings about increased computational complexity due to fully populated matrices or the necessity to track the history in time-dependent problems. Second, solutions to fractional differential equations typically exhibit singularities even for smooth input data, which often lowers convergence rates, unless specialized techniques or locally refined meshes are employed.

\subsection{Numerical Methods for Fractional Operators in 1D}
\label{sec:mel-numerical-methods-1D}

Several classes of numerical methods for fractional differential operators have been proposed in the univariate setting. They can be broadly divided into 2 groups: the first group is related to discretizing the operator in the style of a finite difference method, which leads to methods with a time-marching flavor whereas the second class are Galerkin/collocation methods, which typically involve solving a non-trivial linear system of equations. Detailed discussions can be found in \cite{handbook_vol_3,li-chen18} and in particular the excellent monograph \cite{jin-zhou23} to which we will frequently refer. The simplest model of a time-fractional differential equation (FODE) takes the form
\begin{equation}
\label{eq:mel-FODE}
\der{\alpha}{0,t} u + \lambda u= f(t), \quad \mbox{ $t \in (0,T)$}, 
\end{equation}
where $\lambda \in {\mathbb R}$ and $\der{\alpha}{0,t}$ is either the Riemann-Liouville or the Caputo  fractional differential operators introduced in Section~\ref{ssec:deftimefrac} and $\alpha \not\in {\mathbb N}$. More complex time-fractional PDEs (FPDE) of the form
\begin{equation}
\label{eq:mel-FPDE} 
\der{\alpha}{0,t} u(t,x)  + A u(x) = f(t,x)
\quad \mbox{ $t \in (0,T), \quad x \in \Omega \subset {\mathbb R}^d$}
\end{equation}
for spatial operators $A$ can successfully be treated numerically by combining separate discretizations of $\der{\alpha}{0,t} u$ and $Au$ in a way analogous to the procedure in classical time-dependent problems, viz., the  famous ``method of lines'' in which first a spatial discretization is applied and then the resulting (fractional) ODE system is discretized or the well-known  ``Rothe method''/``method of transverse lines'' in which first the temporal problem is discretized and afterwards the resulting spatial problems are solved numerically. In both approaches for (\ref{eq:mel-FPDE}), it is useful to have a good understanding of  discretizations for the scalar case (\ref{eq:mel-FODE}). Indeed, a diagonalization of the operator $A$ (or its discretization) leads to problems  of the form (\ref{eq:mel-FODE}) so that (\ref{eq:mel-FODE}) may be thought of as a simple model for (\ref{eq:mel-FPDE}). In this section, we will not discuss the discretization  of the spatial operator $A$, and we will mostly restrict ourselves to the case $\alpha \in (0,1)$.

\begin{remark}
\label{remk:mel-singularities}
A hallmark of solutions of (\ref{eq:mel-FODE}) (and thus also (\ref{eq:mel-FPDE})) is that they have a (weak) singularity at the initial time $t = 0$~\cite{diethelm04}, \cite{stynes-handbook}; see also the remarks in Section~\ref{sec:gradient}. Consider the case of Caputo derivatives and $\alpha \in (0,1)$. For smooth $f$ and $\lambda = 0$, solutions  of (\ref{eq:mel-FODE}) can be expected to have the form $u(t) = t^\alpha \tilde u(t) + c$ for $c \in {\mathbb R}$ and smooth $\tilde u$ since a particular solution of $\caputo{\alpha}{0,t} u = t^\beta$ is given by $u(t) = \frac{\Gamma(\beta+1)}{\Gamma(\beta+\alpha+1)} t^{\alpha+\beta}$, \cite[Ex.~1]{stynes-handbook}. For $\lambda \ne 0$, the solution structure is more involved in that the behavior $O(t^\alpha)$ is only the leading order behavior as $t\rightarrow 0$. For example, the function $u(t) = E_{\alpha,1}(-\lambda t^\alpha)$ with the standard (entire) Mittag-Leffler function $E_{\alpha}$ solves the homogeneous equation $\caputo{\alpha}{0,t} u + \lambda u = 0$, \cite[Thm.~{4.3}]{diethelm04}. Hence, for smooth $f$ one should expect solutions of (\ref{eq:mel-FODE}) to have singularity components $t^{j \alpha}$ for $j=1,2,\ldots$.
\end{remark}

\subsubsection{The L1 Scheme and Related Methods}
\label{sec:mel-L1}

The definitions of $\RL{\alpha}{0,t}$ and $\caputo{\alpha}{0,t}$ may motivate various approximations. For example, for $\alpha \in (0,1)$ a discretization of $\caputo{\alpha}{0,t} u$ in the points $0 = t_0 < t_1 < \cdots< t_n$ may be taken to be
\begin{align}
\label{mel:eq:L1}
\caputo{\alpha}{0,t} u(t_n)  & = \frac{1}{\Gamma(1-\alpha)} \int_{\tau=0}^{t_n} (t_n - \tau)^{-\alpha} u^\prime(s)\,d\tau
\approx \frac{1}{\Gamma(1-\alpha)} \sum_{k=0}^{n-1} \int_{\tau=t_{k}}^{t_{k+1}} (t_n - \tau)^{-\alpha} \frac{u(t_{k+1}) - u(t_k)}{t_{k+1}-  t_k}\,d\tau. 
\end{align}
This approximation is known as the L1 method. For the case of equidistant points $t_k =  k \dt $ the sum takes the form of a (discrete) convolution:
\begin{align}
\label{mel:eq:L1-uniform}
\caputo{\alpha}{0,t} u(t_n)  &
\approx \frac{1}{\Gamma(1-\alpha)} \sum_{k=0}^{n-1} b_{n-k-1} \left[ u(t_{k+1})-  u(t_k)\right],
\qquad b_k:= \frac{\dt^{-\alpha}}{\Gamma(2-\alpha)} \left[ (k+1)^{1-\alpha} - k^{1-\alpha} \right].
\end{align}
A numerical method involving $\caputo{\alpha}{0,t}$ is obtained by replacing $u(t_k)$ with approximations $u_k$. When solving numerically a FODE of the form (\ref{eq:mel-FODE}) with initial condition $u(0) = u_0$, the discretization (\ref{mel:eq:L1-uniform}) of $\caputo{\alpha}{0,t}$ leads to a numerical method that has the form of a time-stepping scheme, i.e., the approximations $u_1$, $u_2,\ldots$ are determined in turn.  In contrast to classical ODEs, the cost of computing $u_n$ is not $O(n)$ but $O(n^2)$ due to the ``memory'' effect, i.e., the nonlocal nature of the fractional derivative. It is possible to reduce the cost to from $O(n^2)$ to $O(n \log n)$ (see Section~\ref{sec:mel-fast-methods} below).

\begin{remark}
\label{remk:mel-L2}
The L1 scheme is obtained in (\ref{mel:eq:L1}) by replacing $u$ by its piecewise linear interpolant $\Pi_{1} u$. For $\alpha \in (0,1)$, the L2-scheme is obtained by approximating $u$ by the linear interpolant on the first interval $(0,t_1)$ and on the interval $(t_{j-1}, t_{j})$ by the quadratic polynomial $\Pi_2 u$ interpolating in three consecutive nodes $t_{j-1}$, $t_j$, $t_{j+1}$, and by the quadratic interpolant in the nodes $t_{n-2}$, $t_{n-1}$, $t_n$ for the last interval $(t_{n-1}, t_n)$, \cite[Chap.~{4}]{jin-zhou23}. Likewise, an approximation of the Caputo derivative for $\alpha > 1$ can be based on locally higher order polynomial interpolation, \cite[Sec.~{3}]{li-chen18}. A variant, which includes the fractional form of the Crank-Nicolson scheme (if $\theta =1/2$), is known as Alikhanov's method (also known as $L2-L1_\sigma$) and approximates the fractional derivative at intermediate points $t_{j+\theta} = \theta t_j + (1-\theta) t_{j+1}$ using again polynomial approximations on the intervals $(t_j, t_{j+1})$ and the terminal interval $(t_n, t_{n+\theta})$, see \cite[Sec.~{4.1}]{jin-zhou23}, \cite[Sec.~{3}]{li-chen18}.
\end{remark}

The L1 scheme has a truncation error $O(\dt^{2-\alpha})$ for \emph{smooth} $u$, \cite[Lem.~{4.2}]{jin-zhou23}; however, since by Remark~\ref{remk:mel-singularities} the solution of FODEs typically have at least a (weak) singularity near the initial point $t=0$, lower convergence rates are observed, \cite{stynes-oriordan-gracia17}. For example, for a semidiscrete analysis of (\ref{eq:mel-FPDE}) with an elliptic operator $A$ on uniform temporal meshes \cite[Sec.~{5.1}]{jin-zhou23} shows convergence $O(\dt)$ in the discrete $\ell^\infty$ norm (see, e.g., \cite{mclean-mustapha15,stynes-oriordan-gracia17} for the fully discrete error analysis). For smooth right-hand sides $f$ in (\ref{eq:mel-FODE}) or (\ref{eq:mel-FPDE}), better convergence rates can be obtained using one of the following two devices. The first one is to use the ``corrected L1 scheme'', which is based on a uniform mesh and modifies the initial 2 steps of the time-marching formulation; for example,  errors $O((\dt/t_n)^{2-\alpha})$ are then recovered for the case of (\ref{eq:mel-FPDE}) and an elliptic operator $A$, \cite[Thm.~{5.3}]{jin-zhou23}. The second device is to resort to graded meshes better to resolve the singularity at $t = 0$. Currently successful techniques to obtain error estimates for the L1 scheme and related ones on uniform and non-uniform meshes include two strategies: The first one uses new variants of Gronwall's Lemma \cite[Thm.~{6.3}]{jin-zhou23} that are adapted to FODEs and originate from earlier work on Volterra integral equations, \cite{brunner04}. The second technique exploits that the L1 scheme is associated with an M-matrix and obtains estimates via suitably constructed barrier functions, \cite[Thm.~{6.6}]{jin-zhou23}.

\subsubsection{Convolution Quadrature}
\label{sec:mel-CQ}

The temporal fractional differential operators have a form that makes them amenable to convolution quadrature (CQ), going back to~\cite{lubich88-I,lubich88-II}; see also \cite{banjai-sayas22} and \cite[Chap.~{3}]{jin-zhou23} for an overview. CQ realizes numerically operators of the form (we follow the notation introduced by Lubich)
\begin{equation}
\label{eq:mel-K(partialt)}
 (K(\partial_t) u)(t):= \int_{s=0}^t k(t-s) u(\tau)\,d\tau
\end{equation}
in that good approximations at uniformly spaced points $t_n = n h$ are obtained in the form
$$
K(\partial^\dt_t) u(t_n) = \sum_{j=0}^n w_{n-j} u(t_j)
$$
for suitable weights $w_j$ (more below). Important properties of CQ are:
\begin{enumerate*}[label=\alph*)]
\item determining the weights $w_j$ requires only the Laplace transform of $k$;
\item CQ include the standard (numerical) differentiation and integration operators;
\item CQ has the structure of an algebra under composition, which greatly simplifies the development of methods for operators that are the composition of, e.g., an integral operator of the form~(\ref{eq:mel-K(partialt)}) and a differentiation operator as it appears in the definition of the operators $\caputo{\alpha}{0,t}$ and $\RL{\alpha}{0,t}$.
\end{enumerate*}

\begin{ex}{(Derivation of CQ)}
\label{example:mel-CQ}
Several ways to motivate the form of CQ are given in \cite[Sec.~{2}]{banjai-sayas22}. One way, taken from \cite{lubich88-I}, derives CQ from ODE solvers (here: multistep methods), representing the kernel in terms of the inverse Laplace transform, and using the so-called Z-transform. In more detail but proceeding purely formally, we express $k$ in terms of the inverse Laplace transform and write with the Laplace transform $K = {\mathfrak L} k$ of the kernel $k$ 
\begin{align}
\label{eq:mel-y}
 y(t)&:= (K(\partial_t) u)(t)= \int_{\tau=0}^t k(t-\tau) u(\tau)\,d\tau = \frac{1}{2 \pi \bi} \int_{s \in \sigma + \bi {\mathbb R}} K(s) \int_{\tau=0}^t e^{s(t-\tau)} u(s)\,d\tau\,ds \\
\nonumber 
& = \frac{1}{2 \pi \bi} \int_{s \in \sigma + \bi {\mathbb R}} K(s) v(t; s)\, ds,
\end{align}
where the function $t \mapsto v(t; s)$ solves the initial value problem (IVP) $v^\prime(t; s) = s v(t; s) + u(t)$, $v(0,s) = 0$ (variation of constants formula). Next, this IVP is solved with a numerical method to get an approximation $v_n(s) \approx v(t_n; s)$, which leads to the approximation $y(t_n) \approx y_n:= (2 \pi \bi)^{-1} \int_{s \in \sigma + \bi {\mathbb R}} K(s) v_n(s)\, ds$. To compute $v_n(s)$, one may employ a multistep method with $k$-steps (determined by the coefficients $\alpha_i$, $\beta_i$), which takes the form
\begin{align*}
\frac{1}{\dt} \sum_{j=n-k}^n \alpha_{n-j} v_j & = \sum_{j=n-k}^n \beta_{n-j} (s v_j + u(t_j)),
\end{align*}
where the required additional initial data $v_j$, $j=-k+1,\ldots,-1$ are taken to be zero. The values $v_j$ can be computed recursively. A closed form of the sought values $v_j$ is obtained with the aid of the Z-transform, which maps a sequence $(a_i)_{i \in {\mathbb N}_0}$ to a formal power series ${\mathcal A}(z):=  \sum_{i=0}^\infty a_i z^i$. With the transforms $\alpha(z) =\sum_{i=0}^k \alpha_i z^i$,  $\beta(z) =\sum_{i=0}^k \beta z^i$,  ${\mathcal V}(z):= \sum_{i=0}^\infty v_i z^i$, ${\mathcal U}(z) = \sum_{i=0}^\infty u(t_i) z^i$, one arrives at
\begin{align*}
\frac{1}{\dt} \alpha(z) {\mathcal V}(z) &= \beta(z) \left(s {\mathcal V}(z) + {\mathcal U}(z)\right).
\end{align*}
Solving for ${\mathcal V}$ gives with the generating function $\delta(z):=\alpha(z)/\beta(z)$ of the multistep method
\begin{align*}
{\mathcal V}(z) & = (\delta(z)/\dt - s) ^{-1} {\mathcal U}(z) .
\end{align*}
Hence, the Z-transform ${\mathcal Y}(z) = \sum_{i=0}^\infty y_i z^i$ satisfies
\begin{align*}
{\mathcal Y}(z) & = \frac{1}{2 \pi \bi} \int_{s \in \sigma + \bi {\mathbb R}} K(s) (\delta(z)/\dt - s)^{-1} \, ds\, {\mathcal U}(z)  = K(\delta(z)/\dt) {\mathcal U}(z).
\end{align*}
Expanding
\begin{align*}
K(\delta(z)/\dt) & = \sum_{j=0}^\infty \omega_j z^j
\end{align*}
and recognizing that the power series of $K(\delta(z)/\dt) {\mathcal U}(z)$ involves the convolution of the sequences $(\omega_j)_{j\in{\mathbb N}_0}$ and $(u(t_j))_{j\in{\mathbb N}_0}$ gives the desired relation
\begin{align}
\label{eq:mel-yn}
y_n = \sum_{i=0}^n \omega_{n-i} u(t_i), \qquad n=0,1,\ldots
\end{align}
for the Convolution Quadrature. 
\eex
\end{ex}
The formal calculation in Example~\ref{example:mel-CQ} shows how to design different CQ methods by choosing the function $\delta$, which in turn is given by the choice of the underlying ODE solver. To formulate the method, only the Laplace transform $K$ of the kernel $k$ is required. Even for rather general $K$, the weights $\omega_j$ can be obtained numerically by applying the trapezoidal rule to the Cauchy integral formula for derivatives, \cite[Sec.~{3.1}]{banjai-sayas22}. The underlying ODE solver (e.g., the multistep method in Example~\ref{example:mel-CQ}) is often required to be A-stable. This stipulation derives from the observation that the domain of holomorphy of $K$ often contains the whole right half-plane (e.g., the function $K(s) = s^\alpha$, $\alpha \not\in {\mathbb N}_0$, which is only defined on the slit plane ${\mathbb C} \setminus (-\infty,0]$) and the observation  that for A-stable ODE solvers, the function $\delta$ maps the complex unit ball $B_1(0) \subset {\mathbb C}$ to a subset of the right half-plane $\{z \in {\mathbb C}\,|\,\operatorname{Re} z > 0\}$. However, the requirement of $A$-stability of the underlying ODE solver depends on the problem under consideration; for example, for (\ref{eq:mel-FPDE}) with $\alpha \in (0,1)$ and an elliptic operator $A$ weaker stability requirements suffice, \cite{lubich88-I}, \cite[Sec.~{3.1}]{jin-zhou23}, \cite[Sec.~{8}]{banjai-sayas22}.

For an application of CQ to fractional derivatives, one notes that the Laplace transform $K_\alpha (s) = \int_{t=0}^\infty t^\alpha e^{-st}\,dt$ of the kernel $k(t) = t^\alpha$ is explicitly given by $K_\alpha(s) = \Gamma(1+\alpha) s^{-(1+\alpha)}$. The Laplace transform of the differentiation operator is simply $s$. Hence, for $u$ with $u(0) = 0$, the Laplace transform of $\caputo{\alpha}{0,t} u$ is given by $s^\alpha {\mathfrak  L} u(s)$. That is, for $\alpha \in (0,1)$ the CQ realization of the Caputo derivative for functions vanishing at $t= 0$ can be done with the kernel $K_\alpha(s)$. Since the Riemann-Liouville and Caputo derivatives coincide for such functions, the same CQ techniques can be employed. There is a extensive convergence theory for CQ based on uniform meshes, see, e.g., \cite[Secs.~{2.5}, {5.3}]{banjai-sayas22}. Generally speaking, for sufficiently smooth function $u$, the difference $y(t_n)  - y_n$ of the values (\ref{eq:mel-y}) and (\ref{eq:mel-yn}) is $O(\dt^p)$ if
\begin{enumerate*}[label=\alph*)]
\item
\label{item:mel-dahlquist-1}
the order of the underlying ODE solver is $p$ and
\item
\label{item:mel-dahlquist-2}
$u$ vanishes to sufficient order at $t = 0$.
\end{enumerate*}
Concerning the requirement \ref{item:mel-dahlquist-1}: if one has to base CQ on an A-stable ODE solver, then multistep methods are limited to order $p = 2$ by Dahlquist's barrier. However, CQ can also be based on Runge-Kutta methods, \cite{lubich-ostermann93} and arbitrary order A-stable Runge-Kutta methods such as the Radau IIA formulas are available, \cite[Sec.~{5}]{banjai-sayas22}.  In the context of FODEs and FPDEs, requirement~\ref{item:mel-dahlquist-2} is typically a more limiting restriction as the solution $u$ of the FODE/FPDE cannot be expected to vanish to high order at $t= 0$. Order reduction is then visible. One remedy is to resort to a ``corrected'' CQ akin to that discussed above for the L1 method. For general CQ methods, see, e.g, \cite{lubich04} and \cite[Sec.~{3.2}]{jin-zhou23} for the case (\ref{eq:mel-FPDE}) of an elliptic $A$. An alternative remedy is to base the CQ on non-uniform meshes that are refined towards $t = 0$, see \cite{lopez-fernandez-guo21} for an recent account on such CQ methods. Finally, we point to \cite{banjai-makridakis22} for {\sl a posteriori} error estimates for CQ on general meshes.

\subsubsection{Fast Methods}
\label{sec:mel-fast-methods}

At first sight, the evaluation of fractional differential operators in the nodes $t_0,\ldots, t_n$ leads to complexity $O(n^2)$ due to the time-stepping nature; alternatively, the computation can be viewed as a matrix-vector multiplication with a lower triangular matrix. For uniform point distributions, the convolution structure leads to Toeplitz matrices and algorithms with $O(n \log n)$ complexity using FFT techniques are available for Toeplitz matrices, see \cite[Sec.~{3.3}]{banjai-sayas22} for such algorithms in the context of CQ. Especially in the context of FPDEs, which involves additionally discretizations of spatial operators, various ``fast'' methodologies have emerged that exploit the fact that the pertinent kernel $k(s) = s^\alpha$ is smooth with a singularity only at the origin. As such, compression techniques developed for integral operators (IO) such as fast multipole techniques \cite{greengard-rokhlin97}, wavelet compression, \cite{tausch-white03} and ${\mathcal H}$-matrices, \cite{hackbusch15} that lead to polylogarithmic-linear complexity $O(n \log^\beta n)$ (for some $\beta \ge 0$) for storage and matrix-vector multiplication are, in principle, available. The broad ${\mathcal H}$-matrix format consists of blockwise low-rank matrices, and, in the context of IOs, blockwise low-rank approximation are obtained by approximating the kernel function by, e.g., piecewise polynomials. Therefore, ${\mathcal H}$-matrix compression is well-suited for Galerkin or collocation discretizations of IOs on both uniform and nonuniform meshes; see, e.g., \cite[Sec.{10.6}]{hackbusch15} for a discussion of the case of Volterra IOs. ${\mathcal H}$-matrices come with an (approximate) arithmetic that includes addition, multiplication, and inverse of matrices, \cite{hackbusch15}.

In the context of fractional differential operators, a more popular technique to effectively realize a compression is to approximate the kernel by a (short) sum of exponentials (SOE) instead of piecewise polynomials. A detailed analysis is given in \cite[Sec.~{4.2}]{jin-zhou23} for the L1; see also the survey \cite{mclean18} for related schemes. To illustrate the main point, consider the L1 scheme, i.e., the realization of $\int_{s=0}^{t_n} k(t_n - s) (\Pi_1 u)(s)\,ds$, where $\Pi_1 u$ is the piecewise linear interpolant of $u$ on a uniform mesh with nodes $t_j = j \dt $. One assumes that the kernel $k(s)$ admits the approximation $k(s) \approx \sum_{\ell=0}^{N_k} \omega_\ell e^{-s_\ell s}$ for some small $N_k$ and suitable $\omega_\ell$, $s_\ell \in {\mathbb R}$; the smoothness of $k$ away from $s = 0$ ensures that this can be achieved for the argument $s$ bounded away from $0$. Since this approximation is poor for $s$ close to $0$, the following approximation is based on a decomposition into a ``history'' part and a  ``local'' part:
\begin{align*}
\int_{0}^{t_n} k(t_n - s) \Pi_1 u(s)\,ds  & =
\int_{0}^{t_{n-1}} k(t_n - s) \Pi_1 u(s)\,ds  +
\int_{t_{n-1}}^{t_{n}} k(t_n - s) \Pi_1 u(s)\,ds.
\end{align*}
The ``local'' part $\int_{t_{n-1}}^{t_{n}} k(t_n - s) \Pi_1 u(s)\,ds$ is evaluated exactly and requires only knowledge of the values $u(t_n)$ and $u(t_{n-1})$. The ``history'' part $\int_{0}^{t_{n-1}} k(t_n - s) \Pi_1 u(s)\,ds$ is realized  efficiently with a recursion for the $N_k$ terms $H_\ell(t_n):= \int_{0}^{t_n} e^{-s_\ell (t_n - s)} \Pi_1 u(s)\,ds$:
\begin{align*}
\int_{0}^{t_{n-1}} k(t_n - s) \Pi_1 u(s)\,ds  & \approx
\sum_{\ell=0}^{N_k} \omega_\ell \int_{0}^{t_{n-1}} e^{-s_\ell (t_n - s)} \Pi_1 u(s)\,ds  =
\sum_{\ell=0}^{N_k} \omega_\ell e^{-\dt s_\ell} H_\ell(t_{n-1}),
\end{align*}
and the update formula for the $N_k$ terms $H_\ell(t_n)$ is given by
\begin{align*}
H_\ell(t_n) & = H_\ell(t_{n-1}) + e^{-s_\ell \dt} \int_{t_{n-1}}^{t_n} e^{-s_\ell (t_{n}) - s} \Pi_1 u(s)\,ds,
\end{align*}
which again involves only knowledge of the values $u(t_{n-1})$  and $u(t_n)$ and can be done explicitly. In the context of FODEs, the values $u(t_j)$ are approximated by values $u_j$ and the initial value $u_0$ is prescribed. Correspondingly, one sets $H_\ell(t_0) = 0$. The computational complexity of algorithms based on SOE hinges on the number $N_k$ of terms in the approximation to be small. By \cite[Thm.~{4.2}]{jin-zhou23} one has $N_k = O(|\log (\varepsilon \dt)|^2)$ for a prescribed accuracy $\varepsilon$ so that for accuracy requirement $\varepsilon  = (\dt)^\beta$ ($\beta > 0$) one arrives at complexity $O(n \log^2 n)$.

A method somewhat similar to the technique of using sums of exponentials is the ``fast and oblivious'' CQ of \cite{schaedle-lopez-fernandez-lubich06}. It is a CQ that is applicable if the Laplace transform $K$ of the kernel $k$ is sectorial, i.e., it is holomorphic on a sector $\{z \in {\mathbb C}\,|\, |\operatorname{arg} z| < \pi/2+\delta \}$ for some $\delta > 0$ and satisfies certain decay conditions as $|z| \rightarrow \infty$. A typical representative of this setting are operators of the form $\caputo{\alpha}{0,t} u + A u$ where $A$ is an elliptic spatial operator. Fast and oblivious CQ also reduces the cost to $O(n \log^\beta n)$ ($\beta \ge 0$) by introducing a few auxiliary variables that track the history part in an agglomerated fashion (with controllable error); see \cite[Sec.~{8.2}]{banjai-sayas22} and the related \cite[Sec.~{3.5}]{jin-zhou23}, \cite{doelz-egger-shashkov21}, \cite{khristenko-wohlmuth23} for details.

\subsubsection{Spectral Galerkin Methods for FODEs}
\label{sec:mel-Galerkin-FODEs}

The techniques described in Sections~\ref{sec:mel-L1}, \ref{sec:mel-CQ} have the form of time-marching schemes. A different approach is obtained by spectral techniques, which come in variants such as spectral Galerkin, spectral Petrov-Galerkin, and spectral collocation methods. Starting points concerning the vast literature can be found in the contributions \cite{shen-handbook} and \cite{lischke-handbook} of the handbook \cite{handbook_vol_3},the monograph \cite{zayernouri-wang-shen-karniadakis23} as well as \cite[Sec.~{12}]{jin-zhou23}. See also \cite{chen-shen20, chen-shen-wang16, mao-shen18, zayernouri-karniadakis13} and references therein. While we discuss high order spectral Galerkin methods, Galerkin methods can also be formulated for fixed order methods similar to the SEM discussed below.

Consider (\ref{eq:mel-FODE}) with $\alpha \in (0,1)$ and $\lambda >  0$. The classical Spectral Method  seeks to approximate the solution from the space ${\mathcal P}_N:= \operatorname{span} \{t^i\,|\, i=0,\ldots,N\}$ of polynomials of degree $N$. In view of the expected (weak) singularity at the initial point $t = 0$ (cf.\ Remark~\ref{remk:mel-singularities}) better approximation spaces are required and this has led to the study of several classes of approximation spaces:
\begin{enumerate*}[label=(\roman*)]
\item
\label{item:mel-S1}
The first class, associated with the names of polyfractonomials \cite{zayernouri-karniadakis13} or Generalized Jacobi functions (GJF) \cite{shen-handbook, chen-shen-wang16} consists (in the present case of a singularity only at $t = 0$) of the space $t^{\beta} {\mathcal P}_N$, where the parameter $\beta$ can be chosen to match the singular behavior.
\item
\label{item:mel-S2}
The second class of functions are ``M\"untz polynomials'' $\{\pi(t^\beta)\,|\,  \pi \in {\mathcal P}_N\}$ \cite{hou-hasan-xu18}.
\item
\label{item:mel-S3}
The third class are the recent ``Log orthogonal polynomials'' $\{t^{\beta} \pi(\log t)\,|\, \pi \in {\mathcal P}_N\}$ for chosen $\beta \in {\mathbb R}$, \cite{chen-shen20}, \cite[Sec.~{12.2}]{jin-zhou23}.
\item
\label{item:mel-S4}
The last class is the very versatile class of spectral element methods (SEM), where approximation is effected by piecewise polynomials on suitable meshes, in particular the so-called geometric mesh that is strongly refined towards the singularities, see \cite[Sec.~{3.3.6}]{schwab98} in the context of PDEs and, in the context of fractional operators, \cite{mao-shen18}. Indeed, the geometric mesh is a proven tool to approximate functions with singularities and requires only knowledge of the location of the singularities but not of their precise form. The choices \ref{item:mel-S1}, \ref{item:mel-S2} lead to spectral approximation, i.e., a rate $O(N^{-\gamma})$ for any $\gamma>0$ when approximating functions of the form $u(t) = t^\beta \tilde u(t)$ for smooth $\tilde u$ and known $\beta$; the class \ref{item:mel-S2} can still have good approximation properties for a larger class of functions with singularities at $t = 0$ when the tuning parameter $\beta$ is chosen suitable.
\end{enumerate*}

The classes \ref{item:mel-S3}, \ref{item:mel-S4} are more versatile and do not even need {\sl a priori} knowledge of the precise structure of the singular behavior at $t = 0$ for spectral convergence. It should be noted that Remark~\ref{remk:mel-singularities} informs us of the expected solution structure and that the classes \ref{item:mel-S1} might only resolve the leading order singularity so that only fixed algebraic rates of convergence are achieved. 

Conceptually, a simple way to make use of the good approximation properties of the above spaces are Galerkin methods. The functional framework are the Sobolev spaces $H^\beta(0,T)$, $\beta \in (0,1)$, which are the completions of $C^\infty(0,T)$ under the $H^\beta(0,T)$-norm given by the Slobodeckij norm $\|u\|^2_{H^\beta(0,T)}:= \|u\|^2_{L^2(0,T)} + \int_{t = 0}^T \int_{\tau = 0}^T \frac{|u(t) - u(\tau)|^2}{|t - \tau|^{1+2 \beta}}\, dt\,d\tau$. Taking, e.g., the Riemann-Liouville derivative in (\ref{eq:mel-FODE}), multiplying the equation by a test function $v$ from $H^{\alpha/2}(0,T)$ and integrating over $(0,T)$ one arrives at the weak formulation to find $u \in H^{\alpha/2}(0,T)$ such that
\begin{equation}
\label{eq:mel-FODE-weak-formulation} 
\int_{0}^T \RL{\alpha}{0,t} u(t) v(t)\,dt   = \int_{0}^T \RL{\alpha/2}{0,t} u\, \RL{\alpha/2}{t,T} v\,dt =: b(u,v) = \ell(v):= \int_{0}^T f(t)v(t)\,dt
\end{equation}
for all $v \in H^{\alpha/2}(0,T)$, \cite{li-xu09}, \cite[Sec.~{12.1}]{jin-zhou23}. One should note that for $\alpha \in (0,1)$ functions from $H^{\alpha/2}(0,T)$ do not have a trace at $t = 0$ so that the implied initial condition $u(0) = 0$ is only satisfied in a weak sense, \cite[Rem.~{2.2}]{li-xu09}. A Galerkin discretization of (\ref{eq:mel-FODE-weak-formulation}) is immediately formulated by seeking, for a given a subspace $V_N \subset H^{\alpha/2}(0,T)$ of dimension $N:=\operatorname{dim} V_N$, an approximation $u_N \in V_N$ such that
\begin{equation}
\label{eq:mel_FODE-FEM}
b(u_N,v) = \ell(v) \qquad v \in V_N.
\end{equation}
Once a basis of $V_N$ is selected, the approximation $u_N$ is obtained as the solution of a linear system of equations. The key feature of the weak formulation (\ref{eq:mel-FODE-weak-formulation}) is that the bilinear form, although not symmetric (due to the presence of the different operators $\RL{\alpha}{0,t}$ and $\RL{\alpha}{t,T}$), is elliptic on $H^{\alpha/2}(0,T)$ by \cite{ervin-roop06, li-xu09} in the present case $\alpha \in (0,1)$. In fact,
\begin{equation*} 
c_1 \|u\|^2_{H^{\alpha/2}(0,T)} \leq b(u,u) \leq c_2 \|u\|^2_{H^{\alpha/2}(0,T)} , 
\qquad |b(u,v)|  \leq c_2 \|u\|_{H^{\alpha/2}(0,T)} \|v\|_{H^{\alpha/2}(0,T)} 
\qquad \forall u,v \in H^{\alpha/2}(0,T), 
\end{equation*}
for some $c_1$, $c_2$ depending only on $\alpha$ and $T$. This implies unique solvability of both the continuous problem (\ref{eq:mel-FODE-weak-formulation}) and the discrete problem (\ref{eq:mel_FODE-FEM}) for any $f \in (H^{\alpha/2}(0,T))^\prime$. Additionally, the Galerkin method is quasi-optimal, i.e,.
\begin{equation} 
\label{eq:mel-FODE-quasioptimal} 
\|u - u_N\|_{H^{\alpha/2}(0,T)} \leq \frac{c_2}{c_1} \inf_{v \in V_N} \|u - v\|_{H^{\alpha/2}(0,T)}. 
\end{equation}
We conclude that with good choices of approximation spaces $V_N$, high accuracy of Galerkin methods can be achieved. The Galerkin method leads to fully populated, non-symmetric system matrices. Its entries can be computed with suitable Gauss-Jacobi quadrature rules: For case \ref{item:mel-S4} of the SEM based on standard piecewise polynomials, these calculations are given in \cite{mao-shen18}. For case \ref{item:mel-S1}, i.e., polyfractonomials or GJF, it is observed in \cite{zayernouri-karniadakis13} and \cite[Lemma~{6}]{shen-handbook} or \cite{chen-shen-wang16} that fractional derivatives of GJFs can be expressed in terms of Jacobi polynomials. This motivates the choice of basis of the space GJF space $t^\beta {\mathcal P}_N$. Indeed, by \cite[Lemma~{6}]{shen-handbook}, \cite{chen-shen-wang16} one has (in certain parameter ranges) that $\RL{\alpha}{-1,t}  (1+t)^{\alpha+\alpha_2} P^{(\alpha_1-\alpha,\alpha_2+\alpha)}_n(t)  = \frac{\Gamma(n+\alpha_2+\alpha+1)}{\Gamma(n + \alpha_2+1)} P^{(\alpha_1,\alpha_2)}_n(t)$ with $P^{(\alpha_1,\alpha_2)}_n$ denoting a Jacobi polynomial defined, as usual, on the interval $(-1,1)$. That is, fractional derivatives of functions of the form $(1 + t)^{\alpha} P^{(\alpha,\beta)}_n$ can be expressed through Jacobi polynomials, albeit corresponding to different parameters. This allows one to evaluate the ``stiffness matrix'' arising from $\int_{0}^T \RL{\alpha/2}{0,t} u \RL{\alpha/2}{t,T} v\,dt$ explicitly by a special basis of the GJF space. For (\ref{eq:mel-FODE}) with $\lambda \ne 0$ or (\ref{eq:mel-FPDE}), also the ``mass matrix'' arising from $\lambda \int_{0}^T u  v$ has to be computed; it involves the computation of weighted $L^2$-inner products of polynomials and is thus accessible through Gauss-Jacobi quadratures. The computation of the matrix entries for case \ref{item:mel-S2} of the M\"untz polynomials is done in \cite{hou-hasan-xu18}. For case \ref{item:mel-S3}, i.e., the  ``log orthogonal'' case (with $\beta = 0$), basis functions of ${\mathcal P}_N$ related to the Laguerre are chosen and lead to efficient computations of the system matrix, \cite[Sec.~12.5]{jin-zhou23}.

Often, Spectral Methods are not applied as a Galerkin method but as a more general Petrov-Galerkin (PG) method. This generalization is based on two spaces, a trial space $U_N$ and a test space $V_N$ and takes the form: find $u_N \in U_N$ such
\begin{equation}
\label{eq:mel-PG}
b(u_N, v) = \ell(v) \qquad \forall v \in V_N. 
\end{equation}
The choice $V_N = U_N$ reverts to the Galerkin method. Good choices of the variational formulation (i.e., the form $b$) and $V_N$ can have a significant impact on the computational complexity. In simple cases, the resulting system matrix turns out to be sparse or even diagonal, \cite{zayernouri-karniadakis13} and \cite{shen-handbook}. The flexibility in the choice of the test space $V_N$ often leads to simplified ways of setting up the system matrix and can allow one to give it some structure. We refer to \cite{zayernouri-ainsworth-karniadakis15}, where symmetric system matrices result from a suitable choice of basis of $V_N$; one reason for preferring symmetric matrices over general non-symmetric ones is the availability of fast solver technology.

Especially when considering problems such as (\ref{eq:mel-FPDE}) with a discretization of the spatial operator $A$, a (Petrov-) Galerkin  discretization with its fully populated stiffness and mass matrix leads to very large problems. A common strategy is to try to ``decouple'' as described, e.g., in \cite[Sec.~{12.5}]{jin-zhou23}. To that end, one makes a change of basis of the temporal space $V_N$ so that the (temporal) stiffness matrix ${\mathbf A}$ and the (temporal) mass matrix ${\mathbf M}$ are simultaneously diagonalized, which leads to $N$ fully decoupled spatial problems. If a simultaneous diagonalization is not possible (or not advisable for stability reasons), e.g., because ${\mathbf A}$ and ${\mathbf M}$ are not symmetric, then bringing them simultaneous to upper triangular form \cite[Thm.~{7.7.2}]{golub-vanloan13} is an option and allows one to sequentially solve $N$ spatial problems.
\subsection{Numerical methods for Fractional Operators in Multi-d}
\label{ssec:numfracMulti}
Although the Riemann-Liouville and the Caputo derivative are 1D constructions, they can, in principle,  be employed to define various fractional operators in multi-d. Here, we focus instead on the fractional Laplacian $(-\Delta)^{\gamma/2}$ and restrict our discussion only to two different approaches to defining $(-\Delta)^{\gamma/2}$ for $\gamma \in (0,2)$ that we already introduced in Section~\ref{ssec:defspacefrac}, namely, the spectral fractional Laplacian (SFL) and the integral fractional Laplacian (IFL).
Recall that these operator coincide on ${\mathbb R}^d$, but they differ when restricted to bounded domains $\Omega$; for example, the boundary singularities
of the solutions of (\ref{eq:mel-SFL}) and (\ref{eq:mel-IFL}) differ, \cite{bonito-etal-survey18}.
Good surveys of these and related operators are \cite{lischke-etal20,bonito-etal-survey18,delia-etal20,zayernouri-wang-shen-karniadakis23}.

\subsubsection{Numerical Methods for the Spectral Fractional Laplacian}
Let $\Omega \subset {\mathbb R}^d$ be a bounded Lipschitz domain. Recall the definition of the SFL from Section~\ref{ssec:defspacefrac}. The standard model problem involving the SFL is then: given $f \in H^{-\gamma/2}(\Omega)$ find $u \in \widetilde{H}^{\gamma/2}(\Omega)$ such that
\begin{equation} 
\label{eq:mel-SFL}
(-\Delta)^{\gamma/2} u = f \mbox{ in $\Omega$}, 
\qquad u = 0 \quad \mbox{ on $\partial\Omega$.}
\end{equation}

\subsubsection*{Rational Approximations}

Very good discussions of various methods based on rational approximations of matrix functions are given in \cite{hofreither20} and \cite{burrage-hale-kay12}, which we follow here. A rather intuitive approach to treat (\ref{eq:mel-SFL}) numerically is to discretize the (positive) Laplacian $(-\Delta)$ and work with powers of the resulting system matrix. For example, as described in \cite[Sec.~{3.2.1}]{lischke-etal20}, one may compute discrete eigenfunctions and eigenvalues and mimick the definition of $(-\Delta)^{\gamma/2}$ by replacing the continuous eigenpairs $(\phi_i,\lambda_i)$ by the discrete ones. This approach, called the discrete eigenfunction method (DEM), is the basic version of methods that first discretize the Laplacian and then take matrix powers. To be specific, consider a FEM discretization based on a space $V_N \subset H^1_0(\Omega)$ with basis $\{\varphi_i\,|\, i=1,\ldots,N\}$ with the symmetric stiffness matrix ${\mathbf A} \in {\mathbb R}^{N \times N}$ with entries ${\mathbf A}_{ij} = (\nabla \varphi_j, \nabla \varphi_i)_{L^2(\Omega)}$ and the symmetric mass matrix ${\mathbf M} \in {\mathbb R}^{N \times N}$ with entries ${\mathbf M}_{ij} = (\varphi_j,\varphi_i)_{L^2(\Omega)}$. With the load vector ${\mathbf f} \in {\mathbb R}^N$ with entries $\mathbf f_i = (f,\varphi_i)_{L^2(\Omega)}$, the DEM amounts to approximating the solution $u$ of (\ref{eq:mel-SFL}) by $u^N = \sum_i {\mathbf u}^N_i \varphi_i$ with the vector ${\mathbf u}^N \in {\mathbb R}^N$ given by (see \cite[Sec.~2]{hofreither20} for details)
\begin{equation}
\label{eq:mel-SFL-discrete}
{\mathbf u}^N = \left({\mathbf M}^{-1} {\mathbf A}\right)^{-\gamma/2} {\mathbf f}. 
\end{equation}
(\ref{eq:mel-SFL-discrete}) indicates that, in contrast to the DEM, we may avoid computing all eigenpairs of the matrix pencil $({\mathbf A}, {\mathbf M})$ if we can find good approximations of $\left({\mathbf M}^{-1} {\mathbf A}\right)^{-\gamma/2}$. For example, we may use rational approximations $r(x) \approx x^{-\gamma/2}$ and approximate ${\mathbf u}^N  \approx {\mathbf u}^{N,rat}:= r({\mathbf M}^{-1} {\mathbf A}) {\mathbf f}$. Assuming that $r(x) = O(1/x)$ for $x \rightarrow \infty$, we can write write $r$ in terms of its poles $d_j \in {\mathbb C}$ as $r(x) = \sum_{j} c_j (x - d_j)^{-1}$ for suitable $c_j \in {\mathbb C}$;  thus the evaluation of $r({\mathbf M}^{-1} {\mathbf A}) {\mathbf f}$ reduces to solving (possibly in parallel) linear systems of equations of the form $({\mathbf A} - d_j {\mathbf M}) \widetilde{\mathbf u} = c_j^{-1} {\mathbf M} {\mathbf f}$. For large scale problems, such a system might require preconditioning; however, for $d_j$ with $\operatorname{Re} d_j > 0$, these problems correspond to discretizations of reaction-diffusion equations, for which preconditioners are available. By \cite[Thm.~{1}]{hofreither20}, the additional error $u^N - u^{N,rat}:= \sum_i {\mathbf u}^N_i \varphi_i  - \sum_i {\mathbf u}^{N,rat}_i \varphi_i$ incurred by the rational approximation is controlled by $\|u^N - u^{N,rat}\|_{L^2(\Omega)} \leq \|(\cdot)^{-\gamma/2} - r\|_{L^\infty([\lambda_{min},\lambda_{max}])} \|f\|_{L^2(\Omega)}$, where $\lambda_{min}$, $\lambda_{max}>0$ are the extremal eigenvalues of the pencil $({\mathbf A}, {\mathbf M})$. Algorithms for determining rational approximations $r$ proposed in the literature include proposed BURA, AAA (see \cite{hofreither20}) and BRASIL \cite{hofreither21}; see also \cite{hackbusch19}. Instead of trying to find best rational approximation, one may explicitly construct approximations. For example, starting from the formula  
\begin{equation}
\label{eq:mel-xs}
x^{-s} =  2 \frac{\sin (\pi s)}{\pi} \int_{-\infty}^\infty e^{2 s y} (1 + e^{2y} x)^{-1}\, dy, 
\end{equation}
 we may approximate the integral by a truncated trapezoidal rule (see \cite{hofreither20} for details) to get a rational approximation.  This method is closely related to the ``Balakrishnan formula'' (\ref{eq:mel-balakrishnan-formula}) below. Rational approximations that bring in reduced basis techniques can be found in \cite{danczul-schoeberl22}. An established technique to compute the matrix $g({\mathbf K})$ for a holomorphic $g$ and a matrix ${\mathbf K} \in {\mathbb R}^{N \times N}$ is the contour integral method, described, e.g., in \cite{burrage-hale-kay12} and is applicable here with $g(x) = x^{-\gamma/2}$. Discretizing the contour integral by a quadrature results in having to solve several linear systems of equations and can, in fact, be understood as a specific way of constructing a rational approximation of the function $g$.

\subsubsection*{Extension Methods}

Surveys for these methods include \cite{bonito-etal-survey18} and the literature reviews in \cite{banjai-etal19,banjai-etal23}. The above rational approximation methods start from a discretization of the underlying operator $-\Delta$ and then compute fractional powers on the discrete level. Alternative approaches start from representations of $(-\Delta)^{\gamma/2}$ or its inverse that can then be discretized. One way of achieving this is to use the so-called Caffarelli-Silvestre extension \cite{caffarelli-silvestre07}, which for the present bounded domain case is due to \cite{cabre-tan10,stinga-torrea10}. The extension $U = U(x,y)$ (with $x \in \Omega$, $y \in (0,\infty)$) of the function $u$ solves the following problem on the semi-infinite cylinder ${\mathcal C}:= \Omega \times (0,\infty)$:
\begin{align}
\label{eq:mel-CST}
\nabla_{(x,y)} \cdot \left( y^{1-\gamma} \nabla_{(x,y)} U\right) & = 0 \quad \mbox{ in ${\mathcal C}$},  \\
\quad U & = 0 \quad \mbox {on $\partial\Omega \times (0,\infty)$},  &
\quad \\ \lim_{y \rightarrow 0+} - y^{1-\gamma} \partial_y U(\cdot,y)  & = 2^{\gamma-1} \frac{\Gamma(\gamma/2)}{\Gamma(1-\gamma/2)} f \quad \mbox{ on $y = 0$.}
\end{align}
The sought solution $u$ of (\ref{eq:mel-SFL}) turns out to be $U|_{y = 0}$, \cite{caffarelli-silvestre07,cabre-tan10,stinga-torrea10} so that $U$ is indeed an extension of $u$ to ${\mathcal C}$. Hence, we may approximate $u$ by discretizing the (degenerate) elliptic problem (\ref{eq:mel-CST}). This can be achieved by established methods, e.g., a FEM, \cite{nochetto-otarola-salgado15,banjai-etal19,banjai-etal23}. Since $U$ decays exponentially as $y \rightarrow 0$, one can replace the infinite cylinder ${\mathcal C}$ by a truncated one. The main difficulty with (\ref{eq:mel-CST}) is that the coefficient $y^{1-\gamma}$ degenerates as $y \rightarrow 0$, which requires mesh refinement near $y = 0$, see \cite{nochetto-otarola-salgado15,banjai-etal19,banjai-etal23} for details. An analysis of the solution $U$ of (\ref{eq:mel-CST}) also reveals that $U$ typically has a singularity at $\partial\Omega \times \{0\}$ so that also its trace $u$ has a singularity at $\partial\Omega$, \cite{bonito-etal-survey18}. In the interest of efficiency, especially for higher order methods, this mandates (anisotropic) mesh refinement near $\partial\Omega$, see \cite{banjai-etal19,banjai-etal23} where exponential convergence of high order methods on suitably refined meshes is proved. Discretizations of the fractional Laplacian can be combined with time discretizations in time-dependent problems, see, e.g., \cite{melenk-rieder21,melenk-rieder23}, where again exponential convergence is achieved.

\subsubsection*{Balakrishnan Formula}

Somewhat surprisingly, using the Dunford-Taylor functional calculus, a formula for the solution $u$ of (\ref{eq:mel-SFL}) can be derived, the so-called Balakrishnan formula \cite{bonito-pasciak15}:
\begin{equation}
\label{eq:mel-balakrishnan-formula}
u = \frac{\sin (\pi \gamma/2)}{\pi} \int_{-\infty}^\infty e^{-(\gamma/2) y} \left(\operatorname{I} + e^{-y} (-\Delta) \right)^{-1}\, dy.
\end{equation}
(Note that this can formally be obtained from (\ref{eq:mel-xs}) by setting
$x = (-\Delta)$ and changing the integration variable.)

As is worked out in more detail in Section~\ref{sec:FPDEs-in-pde2path}, 
this integral can be discretized with a quadrature rule as proposed in \cite{bonito-pasciak15,banjai-etal23}, e.g., the truncated trapezoidal rule leading to a sum of elliptic boundary value problems.  The error incurred by this semidiscretization is root exponential in the number of quadrature points, \cite[Thm.~{3.1}]{bonito-lei-pasciak19}. The resulting elliptic problems can be approximated by established discretization techniques, e.g., FEM. Algebraic \cite{bonito-pasciak15} convergence for low order methods and exponential convergence of high order methods on suitable meshes \cite{banjai-etal23} is achieved. The Balakrishnan formula has natural extensions to the time(-fractional) problems, see, e.g.,~\cite{bonito-lei-pasciak17,melenk-rieder23}.

\subsubsection*{Discretization based on Parabolic Tools}

In addition to representations in terms of the extension or the Balakrishnan formula, the inverse of the SFL can be expressed in terms of the heat semigroup. Thus, numerical methods can be based on parabolic solvers, \cite{cusimano-del-teso-gerardo-giorda20}. Further methods based on parabolic solvers include methods due to Vabishchevich described in \cite[Sec.~{8}]{hofreither20}.

\subsubsection{The Integral Fractional Laplacian}
\label{sec:mel-IFL}

The integral fractional Laplacian (IFL) is an operator defined on ${\mathbb R}^d$. Several equivalent definitions are available, \cite{kwasnicki17}. One may define it in terms of the Fourier transformation ${\mathcal F}$  by $(-\Delta)^{\gamma/2} u:= {\mathcal F}^{-1} (|\xi|^{\gamma} {\mathcal F} u)$. In principle, this suggests already a numerical realization via FFT, if the IFL is considered on a (sufficiently large) periodic box. Another class of methods suitable for the IFL in full space are specially designed spectral methods, \cite{mao-shen17,tang-yuan-zhou18,tang-wang-yuan20,sheng-shen-tang-wang20,wang-etal21,papadopoulos2023frame}, and, at least in principle, methods based on approximating by so-called sinc-functions,  \cite{antil-dondl-striet21,antil-dondl-striet22}, which have also been used on bounded domains. We consider the model problem
\begin{equation}
\label{eq:mel-IFL} 
(-\Delta)^{\gamma/2} u = f \quad \mbox{ in $\Omega$}, 
\qquad u = 0 \quad \mbox{ on ${\mathbb R}^d \setminus \overline{\Omega}$.}
\end{equation}
A weak formulation of (\ref{eq:mel-IFL}), which is derived from the representation (\ref{fractionalLaplacian_multiD:singular_integral}), is to find $u \in \widetilde{H}^{\gamma/2}(\Omega)$ such that
\begin{equation}
a(u,v):= \frac{c_{d,\gamma}}{2}\int_{{\mathbb R}^d} \int_{{\mathbb R}^d} \frac{(u(x) - u(y))(v(x) - v(y))}{|x - y|^{d+\gamma}} \, dx\,dy = \ell(v):=  \int_{\Omega} f(x) v(x)\, dx
 \end{equation}
for all $v \in \widetilde{H}^{\gamma/2}(\Omega)$.
(The functions $u$, $v$ defined on $\Omega$ are extended by $0$ outside $\Omega$.)
The bilinear form $a$ is continuous and coercive on $\widetilde{H}^{\gamma/2}(\Omega)$, \cite{nochetto-otarola-salgado15}. Galerkin methods are then obtained by selecting a closed subspace $V_N \subset \widetilde{H}^{\gamma/2}(\Omega)$ and seeking $u_N \in V_N$ such that $a(u_N,v) = \ell(v)$ for all $v \in V_N$. C{\'e}a's lemma immediately gives the quasi-optimality result $\|u - u_N\|_{\widetilde{H}^{\gamma/2}(\Omega)} \leq C \inf_{v \in V_N} \|u - v\|_{\widetilde{H}^{\gamma/2}(\Omega)}$. The spaces $V_N$ may be chosen as spaces of piecewise polynomials on a mesh. Algebraic convergence for low order methods is shown in~\cite{nochetto-otarola-salgado15,banjai-etal19,borthagaray-li-nochetto20,borthagaray-nochetto22}. As in the case of the SFL, the solution $u$ of (\ref{eq:mel-IFL}) has strong singularities at $\partial\Omega$ \cite{faustmann-marcati-melenk-schwab22} so that (anisotropic) refinement is advantageous; exponential convergence of high order methods on suitably anisotropically refined meshes is obtained in~\cite{FMMS-hp}. 

The realization of the Galerkin method faces many challenges that are structurally similar to issues encountered and overcome in the boundary element method (BEM) \cite{sauter-schwab11}: Setting up the system matrix requires the evaluation of (hyper-)singular integrals for which we refer to \cite[Sec.~{5}]{sauter-schwab11} for the simpler case of the BEM, to \cite{chernov-von-petersdorff-schwab15}  and the realizations \cite{acosta-bersetche-borthagaray17,ainsworth-glusa17} in the context of the IFL. Furthermore, the Galerkin matrix is fully populated. However, the matrix compression techniques that were already described in Section~\ref{sec:mel-fast-methods} such as ${\mathcal H}$-matrices \cite{hackbusch15} can be brought to bear and lead to data-sparse representations of the matrix with logarithmic-linear complexity for storage and matrix-vector multiplication; see also~\cite{karkulik-melenk19}. A rather comprehensive discussion of various numerical aspects of Galerkin methods for the IFL can be found in \cite{ainsworth-glusa17}. Optimal preconditioning of the resulting stiffness matrices based on multilevel techniques is addressed in \cite{faustmann-melenk-parvizi21,borthagaray-nochetto-wu-xu23} and in \cite{ gimperlein-stocek-urzua-torres21} based on operator preconditioning. Adaptivity of FEM is analyzed in \cite{faustmann-melenk-praetorius21}.

The IFL admits a Caffarelli-Sivestre extension very similar to (\ref{eq:mel-CST}) where the semi-infinite cylinder is replaced with a half-space ${\mathbb R}^d \times {\mathbb R}^+$ (see also the discussion
in Section~\ref{sec:FA-extension-problem}).
Numerical methods based on this extension have recently been proposed in \cite{faustmann2023fembem}. The Dunford-Taylor calculus is also applicable and leads to representations of the IFL similar to the Balakrishnan formula (\ref{eq:mel-balakrishnan-formula}). Numerical methods based on this have been proposed in \cite{bonito-lei-pasciak19a}.
We finally point to numerical methods that exploit the connection of the IFL with $\alpha$-stable L{\'e}vy processes, \cite{kyprianou-osojnik-shardlow18,sheng-su-xu23}.

\section{Numerical Continuation: Fractional PDEs}
\label{sec:continuation}

In addition to theoretical investigations, it is also useful to look at the bifurcation structure of steady states to understand how solutions behave far from the homogeneous state or critical points when a parameter is varied. Moreover, by varying further parameters, the bifurcation structure itself also changes, which could reveal the effect of crucial quantities on the onset of non-trivial steady states. To compute the bifurcation structure of stationary steady states, a spatial discretization is usually combined with numerical continuation techniques~\cite{DankowiczSchilder, Doedel_AUTO2007, uecker_pde2path_2014}, i.e., one wants to combine spatial discretization methods from Section~\ref{sec:discrete} with analytical bifurcation formulas~\cite{Kielhoefer,Kuznetsov}. For classical reaction--diffusion systems on bounded domains (and~$d=1,2,3$), the global structure can be computed numerically with, for instance, the advanced continuation/bifurcation software \texttt{pde2path}~\cite{uecker_pde2path_2014}. Recently, it has also been used beyond its standard setting, namely to treat cross-diffusion systems~\cite{CKCS, MBCKCS}, optimal control problems~\cite{uecker2019infinite}, network dynamics~\cite{uecker2021pde2path}, and differential geometric PDEs in parametric form~\cite{meiners2023differential}.

It is also quite natural to extend the class of problems treated by \texttt{pde2path} to fractional diffusion PDEs to exploit the software continuation capabilities and routines. In fact, it allows for branch point and Hopf point continuation, continuation of relative equilibria (e.g., traveling waves and rotating waves), branch switching from periodic orbits (Hopf pitchfork/transcritical bifurcation, and period doubling). Therefore, in \cite{ehstand2021numerical}, new capabilities for non-standard diffusion problems involving the fractional Laplacian in the spectral definition on 1D domains have been developed and tested for three important benchmark problems in the non-local setting, aiming at better understanding the effect of fractional diffusion on the steady state bifurcation structure of generic fractional reaction--diffusion systems.

In the following, we show how to implement a discretization of the spectral fractional Laplacian and its embedding within the numerical continuation package \texttt{pde2path}, and we demonstrate results for the Swift--Hohenberg equation in the non-local setting \eqref{eq:SH-frac} obtained thanks to the new continuation capabilities~\cite{ehstand2021numerical}.

\subsection{Fractional PDEs in \texttt{pde2path}}
\label{sec:FPDEs-in-pde2path}
In this section, the FEM discretization of the fractional Laplacian and its implementation in \texttt{pde2path} are presented.
We do not aim to provide a complete guide to the software for beginner users here; for the basic setup and for the notation used in the following we refer to~\cite{dohnal2014pde2path, rademacher2018oopde, uecker_pde2path_2014, uecker2021numerical}. 

The continuation software \texttt{pde2path} in its standard setting is based on the FEM discretization of the stationary elliptic problem exploiting the package \texttt{OOPDE}~\cite{OOPDE} for the FEM discretization. In order to adapt it to handle (spectral) fractional reaction--diffusion equations~\eqref{eq:reaction-diffusion-ss-frac} endowed with homogeneous Dirichlet boundary conditions (we refer to~\cite{ehstand2021numerical} for the case of homogeneous Neumann boundary conditions), we need to discretize the spectral fractional Laplacian~\eqref{eq:spectral_definition} to converted the problem~\eqref{eq:reaction-diffusion-ss-frac} into the algebraic system
\begin{equation}\label{eq:algebraic_system_for_G}
G(\texttt{u},p)=K_{tot}\mathtt{u} - F_{tot}(\mathtt{u})=0,
\end{equation}
where $\texttt{u}\in \mathbb{R}^{n_p N}$ contains the nodal values of $u$ with $n_p$ mesh points, the matrix $K_{tot}\in \mathbb{R}^{n_pN\times n_pN}$ and the vector $F_{tot}\in \mathbb{R}^{n_pN}$ correspond to the finite element discretization of the diffusion and reaction terms (and of the boundary conditions) from equations~\eqref{eq:reaction-diffusion-ss-frac}. To this end, in~\cite{ehstand2021numerical} the Balakrishnan formula representation \cite{balakrishnan1960, kwasnicki17} in the version appearing in \cite[formula (4), p.260]{yosida1968functional} has been discretized via a sinc quadrature for the integral as well as a finite element discretization in space; see also Section~\ref{ssec:numfracMulti} in particular (\ref{eq:mel-balakrishnan-formula}) and the discussion following it. This method has been proposed in~\cite{dohr2019fem, bonito-etal-survey18} for the inverse formula (namely for the inverse fractional Laplacian representation, $\Delta^{-\gamma/2}u$), while in~\cite{ehstand2021numerical} this method has been adapted to the direct formula~\eqref{eq:Balakrishnan_formula-bis}, that is, for $\gamma\in(0, 2)$ and $u\in \widetilde{H}^{\gamma/2}(\Omega)$
\begin{equation}\label{eq:Balakrishnan_formula-bis}
\Delta^{\gamma/2} u(x)=\frac{\sin(\gamma\pi/2)}{\pi}\int_0^\infty \Delta(\nu I-\Delta)^{-1} u(x)\nu^{\gamma/2-1}\textnormal{d}\nu.
\end{equation}
Here are the steps leading to the discretization. Thanks to the change of variable $\nu = \txte^\eta$ ($\textnormal{d}\nu = \txte^\eta \textnormal{d}\eta$) in equation~\eqref{eq:Balakrishnan_formula-bis} and then approximating the integral using sinc quadrature \cite{lund1992sinc} one obtains
\begin{align*} 
\Delta^{\gamma/2}u(x)&=\frac{\sin(\gamma\pi/2)}{\pi}\int_{-\infty}^\infty (\txte^\eta I-\Delta)^{-1}\Delta u(x)\txte^{\eta \gamma/2} ~\textnormal{d}\eta \\
&\approx \frac{\sin(\gamma\pi/2)}{\pi} \kappa \sum_{l=-\infty}^{\infty}\txte^{\kappa l\gamma/2}(\txte^{\kappa l}I-\Delta)^{-1}\Delta u(x),
\end{align*}
where $\kappa$ is the stepsize parameter coming from the sinc quadrature.
Choosing suitable values of $n_+$ and $n_-$ for the truncation, the sum is approximated by
\begin{equation}\label{eq:quadrature_approx}
\Delta^{\gamma/2} u(x) \approx \frac{\sin(\gamma\pi/2)}{\pi} \kappa \sum_{l=-n_{-}}^{n_{+}}\txte^{\kappa l\gamma/2}v_l(x),
\end{equation}
where $v_l \in\mathbb{U}(\mathcal{\tau})$ is the solution of the Galerkin variational problem
\begin{equation*}
\txte^{\kappa l}\int_\Omega v_lw+\int_\Omega \nabla v_l\cdot \nabla w=\int_\Omega \nabla u \cdot \nabla w \qquad \forall w\in \mathbb{U}(\mathcal{\tau})
\end{equation*}
and $\mathbb{U}(\mathcal{\tau})$ is the finite element space satisfying homogeneous Dirichlet boundary conditions. 
The variational problem is then equivalent to the matrix formulation
\begin{equation}\label{eq:linsys}
(\txte^{\kappa l}M+K)\mathtt{v_l}=K\mathtt{u},
\end{equation}
where $M$ and $K$ are the classical FEM mass and stiffness matrices respectively, and $\mathtt{v_l}$ denotes the node vector for $v_l$ and can be obtained as the solution of the linear system \eqref{eq:linsys}. Therefore, the discretization requires solving~$n_-+n_++1$ independent linear systems to obtain an approximation of $\Delta^{\gamma/2} u$.

However, to interface the fractional Laplacian via a discretization of the Balakrishnan formula with the \texttt{pde2path} setup, a matrix approximation of the operator $\Delta^{\gamma/2}$ itself is needed. This matrix,  denoted by $K_\gamma$ to emphasize its dependence on the fractional order $\gamma$, can be built column by column by \emph{testing} the proposed approximation method on the standard basis of Lagrangian $\mathbb{P}_1$ elements. Finally, once the matrix $K_\gamma$ has been built, the standard stiffness matrix $K$ has to be replaced by $-MK_\gamma$ in the expression of $K_{tot}$ in the Matlab routines relevant to the algebraic system of the form~\eqref{eq:algebraic_system_for_G}. We refer to~\cite{ehstand2021numerical} for the detailed description of the implementation in \texttt{pde2path}, the Matlab code and the choice of the quantities $\kappa,\,n_+,\,n_-$. In addition, the Matlab code is also freely available in the GitHub folder \cite{VideoFolder}.

\subsection{The Fractional Swift--Hohenberg Equation}\label{sec:SHe-all}

Thanks to the new capabilities of the continuation software, it is possible to investigate the effect of the fractional order~$\gamma$ on the steady state bifurcation structure of benchmark problems. As in Section~\ref{sec:spectral2}, we consider as a case study the fractional Swift--Hohenberg equation with competing cubic--quintic nonlinearities~\eqref{eq:SH-frac} and homogeneous Dirichlet boundary conditions. For the classical Swift--Hohenberg equation (obtained with $\gamma=2$), it is known that bifurcation diagrams exhibit the so-called homoclinic snaking~\cite{BURKE2007681, avitabile2010snake, houghton2011swift}. A non-trivial branch of stationary solutions bifurcates from a branch point and, for a suitable range of the main bifurcation parameter, undergoes several fold bifurcations leading to multi-stability of the solutions along the branch that are often important localized patterns. Looking at the bifurcation diagram, this branch presents the characteristic snaking-type structure. Recently, in~\cite{ehstand2021numerical}, also the bifurcation structure of the Swift--Hohenberg equation with nonlocal reaction terms (but also other well-known PDEs) has been investigated with a focus on understanding how the snaking scenario changes under the nonlocal influence (namely decreasing the fractional power $\gamma$) as already carried out for nonlocal reaction terms in~\cite{MorganDawes}. As a first step, the Swift--Hohenberg equation can be converted into a system of second-order PDEs by substituting $(u_1,u_2) = (u, \Delta^{\gamma/2} u)$ in equation~\eqref{eq:SH-frac}. The resulting system is
\begin{equation*}
\begin{pmatrix} 1 & 0\\ 0 & 0 \end{pmatrix} \partial_t \begin{pmatrix} u_1\\ u_2\end{pmatrix}=
    \begin{pmatrix}
    -\Delta^{\gamma/2} u_2 -2u_2 -(1-p)u_1 + q u_1^3 - u_1^5\\
    \Delta^{\gamma/2} u_1 - u_2
    \end{pmatrix},
\end{equation*}
which fits the \texttt{pde2path} extended framework. As in~\cite{ehstand2021numerical}, we consider equation~\eqref{eq:SH-frac} defined on a bounded interval $\Omega=(-5\pi, 5\pi) \subset \mathbb{R}$ of length $L=10\pi$ together with homogeneous Dirichlet boundary conditions $u(-5\pi) = u(5 \pi) = 0$ and $\Delta^{\gamma/2}u(-5 \pi) =\Delta^{\gamma/2}u(5 \pi) = 0$. In the system two parameters are involved: we fix the parameter~$q=2$, while parameter $p$ is taken as the bifurcation parameter. For the space discretization, we choose the mesh size $h=0.04$, corresponding to $n_p=786$ mesh points. The quantities $\kappa$, $n_+$ and $n_-$ in the quadrature approximation~\eqref{eq:quadrature_approx} of the fractional Laplacian have been chosen in~\cite{ehstand2021numerical} following~ \cite[Corollary 2]{dohr2019fem} where a similar discretization method was used. In detail, the value of $\kappa$ is related to the mesh size $h$ and the values $n_+$ and $n_-$ are chosen proportional to $1/\kappa^2$ as follows
\begin{equation*}
\kappa = \frac{1}{|\log(h)|},\qquad n_{+} \coloneqq \left\lceil\frac{\pi^2}{4(1-\gamma/2)\kappa^2}\right\rceil,\qquad n_{-} \coloneqq \left\lceil\frac{\pi^2}{2\gamma\kappa^2}\right\rceil.
\end{equation*}
It is then possible to investigate some of the effects of super-diffusion on the steady state bifurcation structure and solution profiles. We highlight here the main effects that were observed for the bifurcation points on the homogeneous branch, the snaking branch in the bifurcation structure, and the solution profile.

\textbf{Bifurcation points on the homogeneous branch.}
As already mentioned in Section~\ref{sec:spectral2} from equation~\eqref{eq:subsequent-bif}, the bifurcation points on the homogeneous branch accumulate at the value of the bifurcation parameter at which a change of stability of the homogeneous solution occurs, namely at $p_c=0$, as the fractional order $\gamma$ decreases.

\begin{figure}
 \centering
 \begin{subfigure}[b]{0.45\textwidth}
         \centering
         \begin{overpic}[width=\textwidth]{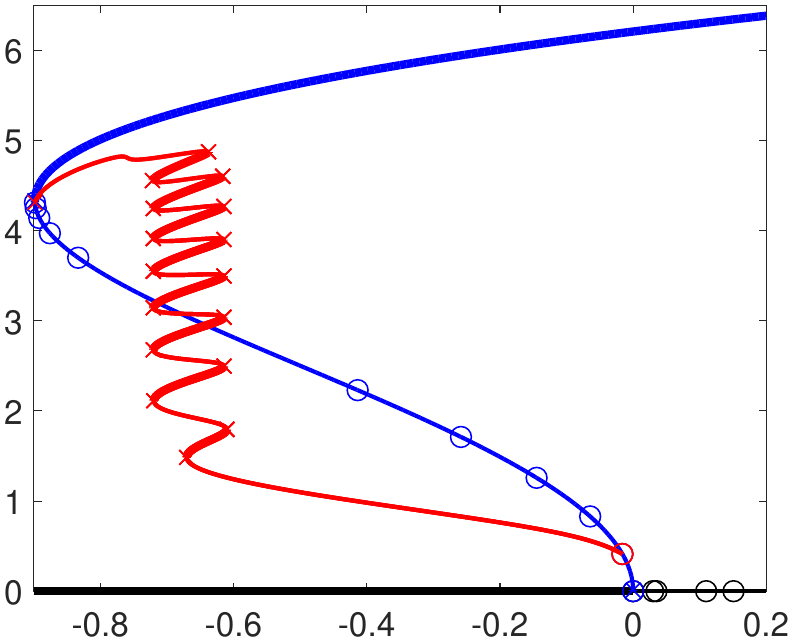}
         \put(-8,45){\rotatebox{90}{$\|u\|_{L^2}$}}
         \put(99,5){$p$}
         \end{overpic}
         \caption{$\gamma=1.8$}\label{fig:sh-1Dbif-s09}
 \end{subfigure}
 \hspace{1cm}
 \begin{subfigure}[b]{0.45\textwidth}
         \centering
         \begin{overpic}[width=\textwidth]{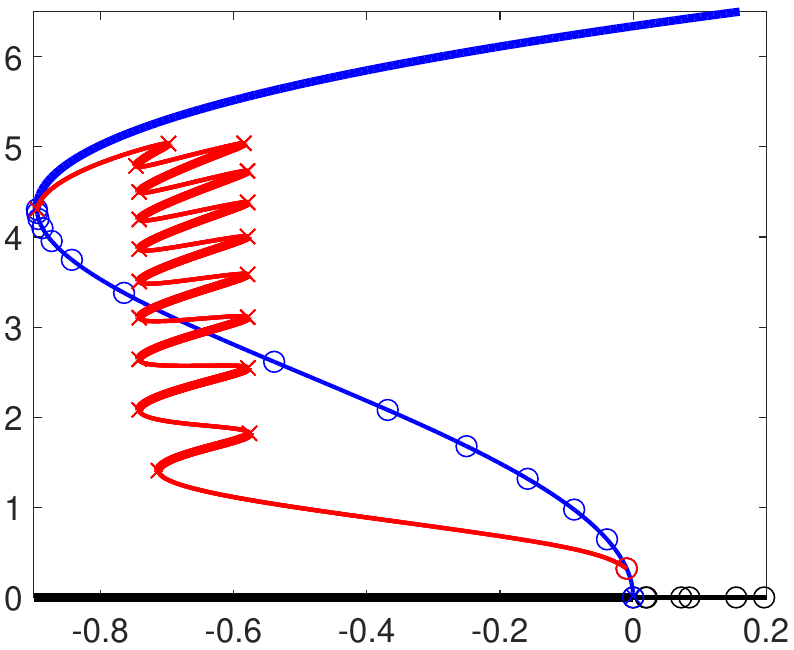}
         \put(-8,45){\rotatebox{90}{$\|u\|_{L^2}$}}
         \put(99,5){$p$}
         \end{overpic}
         \caption{$\gamma=1.4$}\label{fig:sh-1Dbif-s07}
 \end{subfigure} \\[0.5cm]
 \begin{subfigure}[b]{0.45\textwidth}
         \centering
         \begin{overpic}[width=\textwidth]{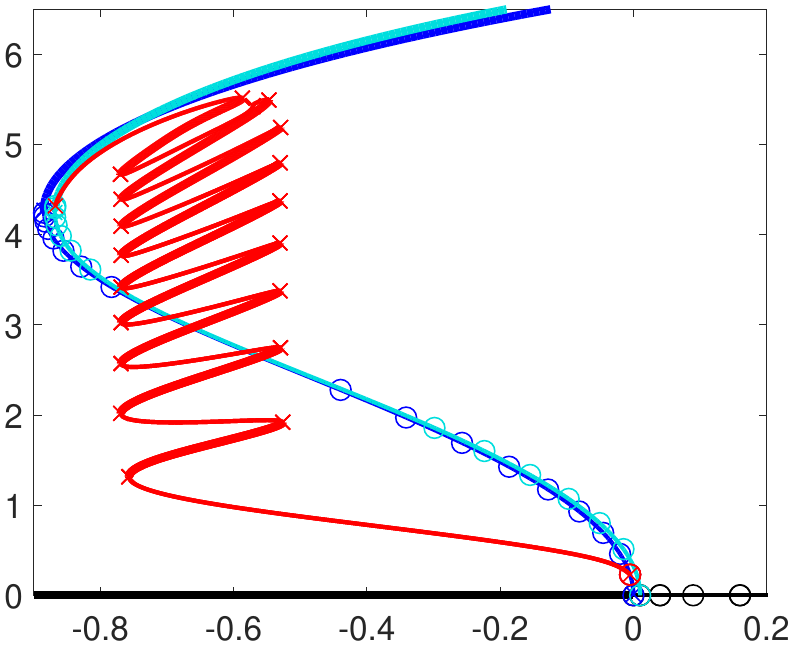}
         \put(-8,45){\rotatebox{90}{$\|u\|_{L^2}$}}
         \put(99,5){$p$}
         \end{overpic}
         \caption{$\gamma=1$}\label{fig:sh-1Dbif-s05}
 \end{subfigure}
  \hspace{1cm}
  \begin{subfigure}[b]{0.45\textwidth}
         \centering
         \begin{overpic}[width=\textwidth]{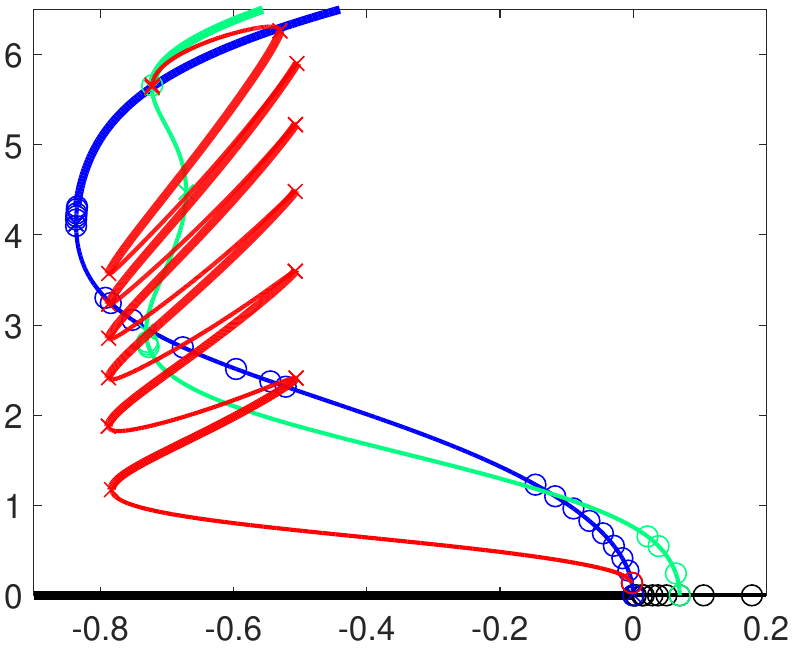}
         \put(-8,45){\rotatebox{90}{$\|u\|_{L^2}$}}
         \put(99,5){$p$}
         \end{overpic}
         \caption{$\gamma=0.6$}\label{fig:sh-1Dbif-s03}
 \end{subfigure}
 \caption{The snaking branch of the fractional Swift--Hohenberg equation~\eqref{eq:SH-frac} gets stretched (``bloated'') as the fractional order~$\gamma$ decreases. The bifurcation diagrams with respect to the bifurcation parameter~$p$ have been computed thanks to the new capabilities of the continuation software \texttt{pde2path}. The equation is considered on a bounded domain~$\Omega=(-5\pi, 5\pi)$ with homogeneous Dirichlet boundary conditions, and parameter $q=2$. Thick and thin lines denote stable and unstable solutions, while circles and crosses indicate branch and fold points respectively. This figure was published in~\cite{ehstand2021numerical}, Copyright Elsevier (2021).}
\label{fig:sh-bif-firstImpression}
\end{figure}

\begin{figure}[!ht]
\centering
 \begin{overpic}[width=0.4\textwidth]{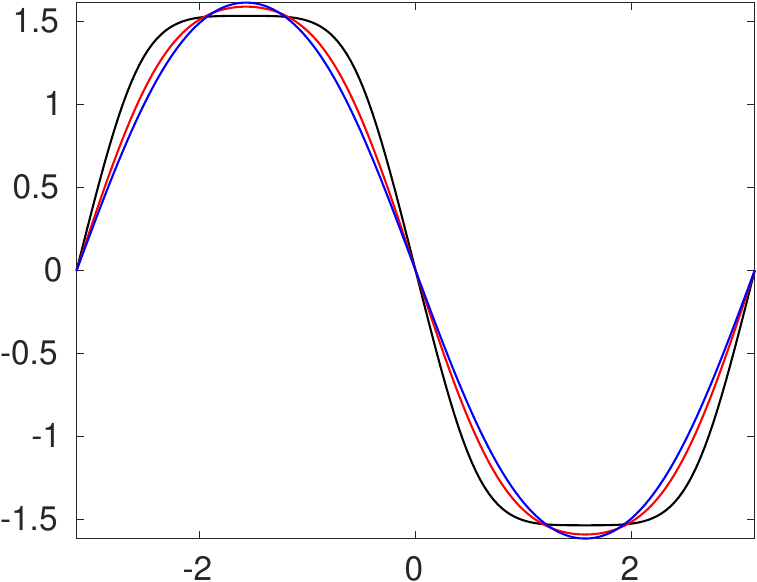}
 \put(-5,35){\rotatebox{90}{$u(x)$}}
 \put(100,5){$x$}
 \end{overpic}
\caption{Change in the solution profile to the fractional Swift--Hohenberg equation~\eqref{eq:SH-frac} as the fractional order~$\gamma$ decreases.  The equation is considered on a bounded domain~$\Omega=(-5\pi, 5\pi)$ with homogeneous Dirichlet boundary conditions, and parameter $q=2$. In the bifurcation diagrams, the solutions correspond to $p=-0.4$ and are located on the first (blue) branch bifurcating from the homogenous branch. Different colors denote different fractional orders: blue for $\gamma=1.8$, red for $\gamma=1$ and black for $\gamma=0.6$.  This figure was published in~\cite{ehstand2021numerical}, Copyright Elsevier (2021).}
\label{fig:sh-sol-profiles-0907005}
\end{figure}

\begin{figure}
\centering
\begin{subfigure}[t]{0.35\textwidth}
\centering
\hspace{-0.5cm}
\begin{overpic}[width=0.9\textwidth]{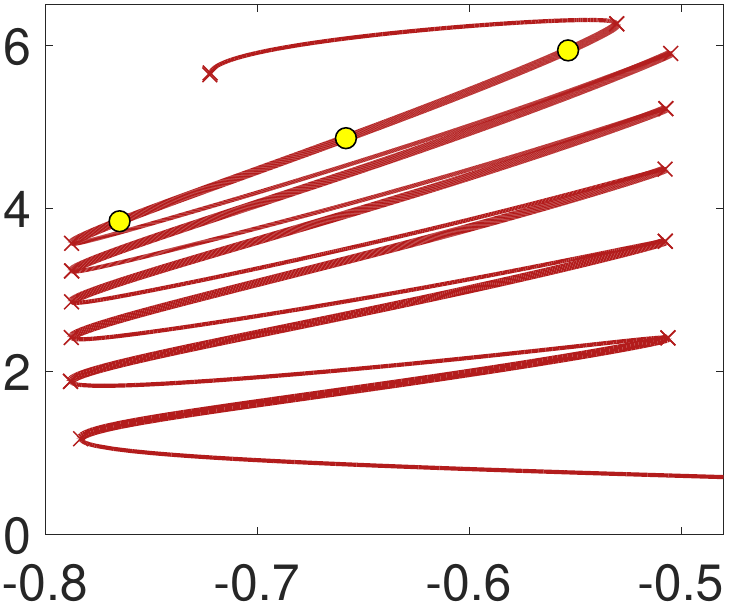}
\put(-10,45){\rotatebox{90}{$\|u\|_{L^2}$}}
\put(100,5){$p$}
\put(10,58){(a)}
\put(40,67){(b)}
\put(65,75){(c)}
\end{overpic}
\label{fig:sh-s03-tilted-labeled}
\end{subfigure}

\begin{subfigure}[t]{0.2\textwidth}
 \centering
 \begin{overpic}[width=\textwidth]{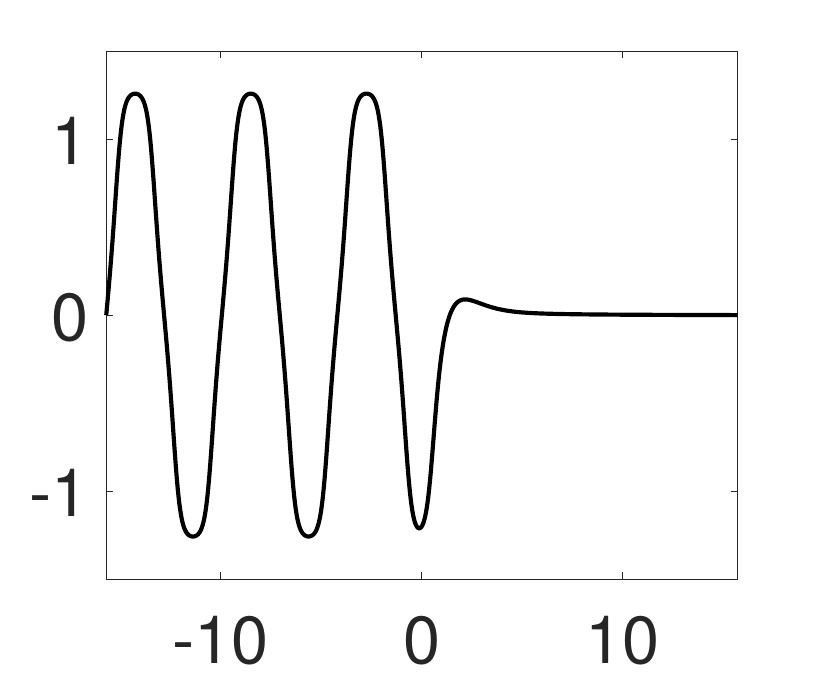}
 \put(-7,45){\rotatebox{90}{$u(x)$}}
 \put(95,5){$x$}
 \end{overpic}
 \caption{}\label{fig:sh-s03-pt5300}
\end{subfigure}
\begin{subfigure}[t]{0.2\textwidth}
 \centering
 \begin{overpic}[width=\textwidth]{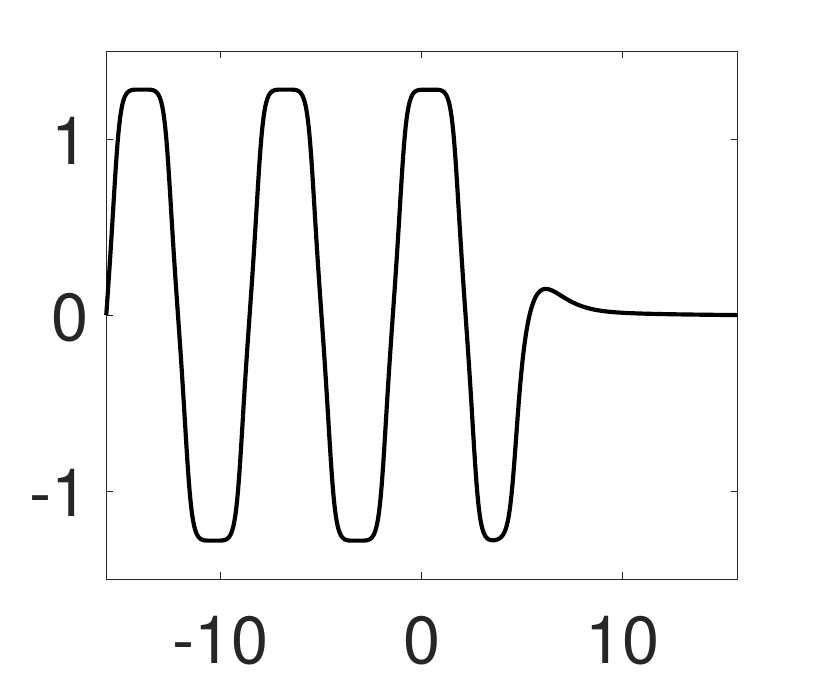}
 \put(95,5){$x$}
 \end{overpic}
 \caption{}\label{fig:sh-s03-pt5600}
\end{subfigure}
\begin{subfigure}[t]{0.2\textwidth}
 \centering
 \begin{overpic}[width=\textwidth]{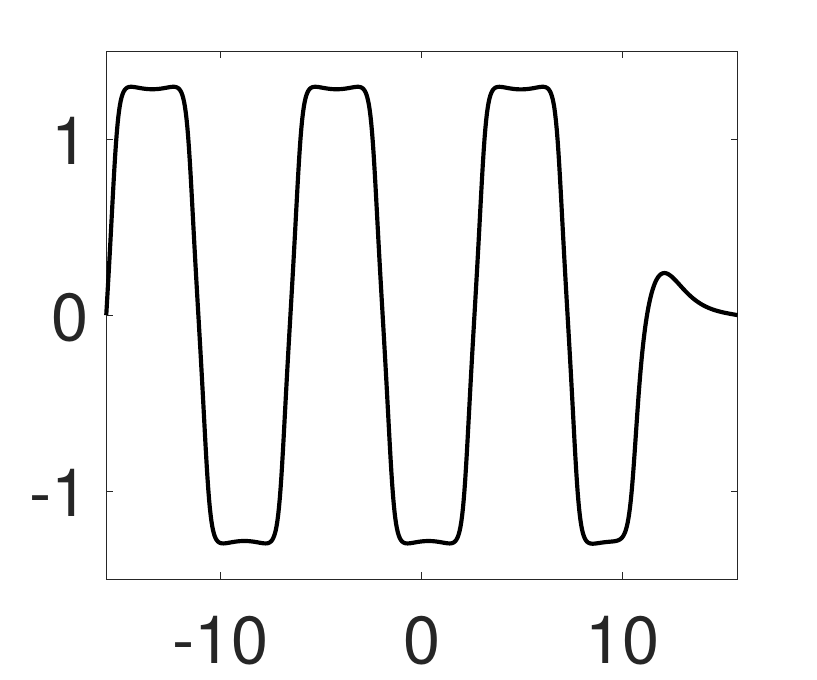}
 \put(95,5){$x$}
 \end{overpic}
 \caption{}\label{fig:sh-s03-pt5900}
\end{subfigure}
\caption{Change in the solution profile to the fractional Swift--Hohenberg equation~\eqref{eq:SH-frac} along the upper part of the snaking branch. The equation is considered on a bounded domain~$\Omega=(-5\pi, 5\pi)$ with homogeneous Dirichlet boundary conditions, and parameter $q=2$, while the fractional order is $\gamma=0.6$.  This figure was published in~\cite{ehstand2021numerical}, Copyright Elsevier (2021).}
\label{fig:sh-tilted-branch-s03}
\end{figure}

\textbf{Bifurcation structure and snaking branch.} The fractional order also affects the bifurcation structure of problem~\eqref{eq:SH-frac}. In order to get a general impression of the effect, Figure~\ref{fig:sh-bif-firstImpression} compares the bifurcation diagram for decreasing fractional orders (namely $\gamma = 1.8$, $\gamma = 1.4$,  $\gamma = 1$ and $\gamma = 0.6$). The homogeneous branch is shown in black, the first branch bifurcating at $p=0$ from the homogeneous branch in blue and the first (homoclinic) snaking branch in red. For $\gamma=1.8$ the bifurcation structure is very similar to the one of the standard Swift--Hohenberg equation but it deforms as $\gamma$ varies. The most evident changes regard the snaking branch. In fact, while the width of the snaking branch significantly increases (it gets stretched or ``bloated'') as the fractional order decreases, the back-and-forth oscillations forming the snaking tend to tilt on the left for smaller~$\gamma$. Further, also the number of turns on the snaking branches varies: we pass from $8$ pair of turns for $\gamma=1.8$ to $9$ pairs for $\gamma=1.4$ and $\gamma=1$, decreasing further to $6$ pairs for $\gamma=0.6$. Since the first snaking branch corresponds to front solutions as shown in Figure~\ref{fig:sh-tilted-branch-s03}, at each snaking turn one oscillation is added to the solution profile. This means that  the number of turns determines the final number of oscillations in the solution and consequently that the snaking branch has to reconnect to other periodic branches originating from subsequent bifurcations on the homogeneous instead of the first periodic branch. In particular, for $\gamma=1$ and $\gamma=0.6$ it reconnects to the branches originating from the third and ninth bifurcation points, respectively, on the homogeneous black branch (also with different wave numbers). Note that the deformation of the snaking branch is known to occur for Swift--Hohenberg equations with non-local reaction \cite{MorganDawes} and that the re-connection of the snaking branch not to the first branch also occurs for standard (Laplacian) Swift--Hohenberg equations varying the domain length \cite{bergeon2008eckhaus}. A less evident deformation is that the non-homogeneous steady states have slightly larger $L^2$-norm values as the fractional order $\gamma$ decreases and has an effect on the solution profiles.

\textbf{Solution profile.} Finally, it has been observed that, for a fixed value of the bifurcation parameter, the non-homogeneous stationary solutions tend to have flatter peaks or valleys as the fractional order of the Laplacian decreases. This effect has been named as ``sharpening the teeth'' of the spatial solution profiles. In Figure~\ref{fig:sh-sol-profiles-0907005}, solution profiles located on the first (blue) branch bifurcating from the homogenous branch at $p=4$ for different values of the fractional order $\gamma$ (see Figure~\ref{fig:sh-bif-firstImpression}) are shown. Whereas the solution profiles for $\gamma=1.8$ and $\gamma=1.4$ resemble the basic periodic pattern specific to the branch (i.e., to the eigenfunction associated to the corresonding mode), we clearly see that at $\gamma=1$ the profile significantly ``flattens'' towards the maximum of oscillation. The same happens to the solution profile along the upper part of the snaking branch for $\gamma=0.6$, as illustrated in Figure~\ref{fig:sh-tilted-branch-s03}. Moving up the snaking branch, namely increasing the parameter $p$, the oscillations in the solution profile significantly enlarge while the peaks and valleys are squared off. 

Moreover, in~\cite{ehstand2021numerical} it has been shown that all these effects are also common to the fractional Allen--Cahn equation and the fractional Schnakenberg system, contributing to a better understanding of the effects of the fractional Laplacian on the steady solutions of generic reaction--diffusion systems on bounded domains. Overall, there are currently still no additional detailed numerical continuation case studies available for various fractional dissipative PDEs but we anticipate this research direction to naturally grow in the coming years.

\section{Conclusion \& Outlook}
\label{sec:outlook}

In this chapter, we have summarized several important research directions for fractional dissipative PDEs. As mentioned in the introduction, it is impossible at this point to provide a completely comprehensive overview of the area as it has grown so rapidly in recent years. Yet, even in a relatively condensed format, several conclusions for fractional dissipative PDEs have become apparent. First, the general dynamics approach known for classical dissipative PDEs can be transferred in many cases. One starts to analyze the invariant sets, studies their local stability via spectral theory, gets a more global analytical view via gradient/Lyapunov-function structures, repeats the strategy for moving patterns, and finally transitions to numerical discretization and numerical continuation to fill out the dynamical picture. Second, the technical approach often differs significantly. Many methods from classical derivative calculus are not available for fractional derivatives. This often necessitates very significant technical extensions, e.g., for transform methods, spectral theory, energy estimates, bifurcation theory, sub/super-solutions, new integral discretizations, and updated numerical continuation packages. For some problems, entirely new methods have been employed such as converting by extension back to local derivatives. Third, once methods are available, one notices often significant quantitative changes of the dynamics. This includes, e.g., changes for the spectrum, space-time asymptotics, bifurcation points, numerical cost/accuracy or spatial structure of patterns. Yet, there are also numerous phenomena from classical dissipative PDEs that are qualitatively robust. This is re-assuring as one wants physical dissipative systems to behave robustly under small structural perturbations of the model. 

Currently, the field of fractional dissipative PDEs still poses immense challenges as only the dynamics of some basic phenomena have been understood, so we can only outline a few possible future directions here: (I) We have mostly focused on purely time-fractional or purely space-fractional PDEs with Riemann-Liouville, Caputo and Riesz--Feller derivatives. Although there is a good motivation for this viewpoint regarding diffusive processes, it is certainly plausible for applications to directly think in terms of space-time modelling leading to nonlocal operators in space and time. For example, for subdiffusive processes on complex networks, one expects the mean-field limit PDEs to have nonlocal operators in space and time, potentially even new coupled operators beyond current standard examples. (II) Since we have mostly focused on low-dimensional, or even scalar, problems, the analysis of the existence and stability theory for steady states was still somewhat tractable. Once spatial multi-component systems are considered, spectral problems are likely going to become much more complicated, which is an effect already well-known for classical dissipative PDEs. (III) The last issue also connects to a topic we have not covered here concerning invariant/inertial manifolds for fractional dissipative PDEs. We conjecture that existing techniques will transfer once suitable spectral gap conditions are satisfied but verifying these conditions could be much more challenging as already indicated in our discussion of spectral theory for time-fractional PDEs. (IV) On a global level, the theory of nonlocal gradient flows is extremely active currently in multiple directions, so we anticipate that fractional PDEs will be embedded within these efforts as key benchmark examples. (V) One the level of pattern formation, we only touched upon traveling waves, which have a relatively simple (translation) symmetry group and upon Turing patterns generated at local bifurcation points. Yet, for classical reaction-diffusion PDEs, there are many other phenomena of key interest such as chaos, chimeras, modulated patterns, and spirals, just to name a few. It is likely that fractional derivatives can fundamentally alter the mechanisms for these in analogy to pattern formation results observed when one transitions from lattice or all-to-all network models to intermediate-range nonlocal coupling. (VI) From a numerical perspective, some excellent discretization schemes do exist but we are still very far from having practical software packages available that would allow for a seamless use of these schemes to study fractional nonlinear dynamics. Hence, a major software development effort in this direction would be highly welcome. (VII) Lastly, the PDEs we focused on had relatively simple polynomial nonlinearities and mostly semilinear structures providing a connection to several classical parabolic PDEs. Yet, many complex systems applications are bound to lead us far beyond these mathematical base cases.\medskip 

\textbf{Acknowledgments:} 
FA was supported by the Austrian Science Fund (FWF) via the FWF-funded SFB \# F65 and FWF \# P33151. GA was supported by JSPS KAKENHI Grant Numbers JP21KK0044, JP21K18581, JP20H01812 and JP20H00117 as well as by the Research Institute for Mathematical Sciences, an International Joint Usage/Research Center located in Kyoto University.  CK would like to thank the VolkswagenStiftung for support via a Lichtenberg Professorship. JMM was supported by Austrian Science Fund (FWF) via the SFB ``Taming complexity in partial differential systems'' (grant SFB F65). CS is a member of the INdAM-GNFM Group and was partially supported by the OeAD project FR 03/2022.

\bibliographystyle{unsrt}
\bibliography{refs,JRrefs,FArefs,JMMRefs,GArefs,CSrefs}

\begin{thebibliography}{100}

\bibitem{Britton2}
N.F. Britton.
\newblock {\em Reaction-Diffusion Equations and Their Applications to Biology}.
\newblock Academic Press, 1986.

\bibitem{Grindrod}
P.~Grindrod.
\newblock {\em Patterns and Waves: The Theory and Applications of
  Reaction-Diffusion Equations}.
\newblock Clarendon Press, 1991.

\bibitem{Murray2}
J.D. Murray.
\newblock {\em Mathematical Biology II: Spatial Models and Biomedical
  Applications}.
\newblock Springer, 3rd edition, 2003.

\bibitem{Amann4}
H.~Amann.
\newblock {\em Linear and Quasilinear Parabolic Problems: Volume I: Abstract
  Linear Theory}.
\newblock Springer, 1995.

\bibitem{Evans}
L.C. Evans.
\newblock {\em Partial Differential Equations}.
\newblock AMS, 2002.

\bibitem{Lunardi}
A.~Lunardi.
\newblock {\em Analytic Semigroups and Optimal Regularity in Parabolic
  Problems}.
\newblock Springer, 1995.

\bibitem{Rothe}
F.~Rothe.
\newblock {\em Global Solutions of Reaction-Diffusion Systems}.
\newblock Springer, 1984.

\bibitem{Smoller}
J.~Smoller.
\newblock {\em Shock Waves and Reaction-Diffusion Equations}.
\newblock Springer, 1994.

\bibitem{Robinson1}
J.C. Robinson.
\newblock {\em Infinite-Dimensional Dynamical Systems}.
\newblock CUP, 2001.

\bibitem{Temam}
R.~Temam.
\newblock {\em Infinite-Dimensional Dynamical Systems in Mechanics and
  Physics}.
\newblock Springer, 1997.

\bibitem{Henry}
D.~Henry.
\newblock {\em Geometric Theory of Semilinear Parabolic Equations}.
\newblock Springer, Berlin Heidelberg, Germany, 1981.

\bibitem{KuehnBook1}
C.~Kuehn.
\newblock {\em PDE Dynamics: An Introduction}.
\newblock SIAM, 2019.

\bibitem{SchneiderUecker}
G.~Schneider and H.~Uecker.
\newblock {\em Nonlinear {PDEs}: A Dynamical Systems Approach}.
\newblock AMS, 2017.

\bibitem{ChowHale}
S.-N. Chow and J.K. Hale.
\newblock {\em Methods of Bifurcation Theory}.
\newblock Springer, 1982.

\bibitem{Kielhoefer}
H.~Kielhoefer.
\newblock {\em Bifurcation Theory: An Introduction with Applications to PDEs}.
\newblock Springer, 2004.

\bibitem{Kuznetsov}
Yu.A. Kuznetsov.
\newblock {\em Elements of Applied Bifurcation Theory}.
\newblock Springer, New York, NY, 3rd edition, 2004.

\bibitem{Zeidler3}
E.~Zeidler.
\newblock {\em Nonlinear Functional Analysis and its Applications {III}:
  Variational Methods and Optimization}.
\newblock Springer, 2013.

\bibitem{Kato}
T.~Kato.
\newblock {\em Perturbation Theory for Linear Operators}.
\newblock Springer, 1980.

\bibitem{KapitulaPromislow}
T.~Kapitula and K.~Promislow.
\newblock {\em Spectral and Dynamical Stability of Nonlinear Waves}.
\newblock Springer, 2013.

\bibitem{Sandstede1}
B.~Sandstede.
\newblock Stability of travelling waves.
\newblock In B.~Fiedler, editor, {\em Handbook of Dynamical Systems}, volume~2,
  pages 983--1055. Elsevier, 2001.

\bibitem{DankowiczSchilder}
H.~Dankowicz and F.~Schilder.
\newblock {\em Recipes for Continuation}.
\newblock SIAM, 2013.

\bibitem{Doedel_AUTO2007}
E.J. Doedel, A.~Champneys, F.~Dercole, T.~Fairgrieve, Y.~Kuznetsov, B.~Oldeman,
  R.~Paffenroth, B.~Sandstede, X.~Wang, and C.~Zhang.
\newblock Auto 2007p: Continuation and bifurcation software for ordinary
  differential equations (with homcont).
\newblock {\em http://cmvl.cs.concordia.ca/auto}, 2007.

\bibitem{uecker_pde2path_2014}
H.~Uecker, D.~Wetzel, and J.D.M. Rademacher.
\newblock \texttt{pde2path} - a {M}atlab package for continuation and
  bifurcation in 2{D} elliptic systems.
\newblock {\em Numerical Mathematics: Theory, Methods and Applications},
  7(1):58–106, 2014.

\bibitem{AllenCahn}
S.M. Allen and J.W. Cahn.
\newblock A microscopic theory for antiphase boundary motion and its
  application to antiphase domain coarsening.
\newblock {\em Acta Metallurgica}, 27(6):1085--1905, 1979.

\bibitem{ChafeeInfante}
N.~Chafee and E.F. Infante.
\newblock A bifurcation problem for a nonlinear partial differential equation
  of parabolic type.
\newblock {\em Applic. Anal.}, 4(1):17--37, 1974.

\bibitem{Nagumo}
J.~Nagumo, S.~Arimoto, and S.~Yoshizawa.
\newblock An active pulse transmission line simulating nerve axon.
\newblock {\em Proc. IRE}, 50:2061--2070, 1962.

\bibitem{CrossHohenberg}
M.C. Cross and P.C. Hohenberg.
\newblock Pattern formation outside of equilibrium.
\newblock {\em Rev. Mod. Phys.}, 65(3):851--1112, 1993.

\bibitem{Schloegl}
F.~Schl{\"o}gl.
\newblock Chemical reaction models for non-equilibrium phase transitions.
\newblock {\em Zeitschrift f{\"u}r Physik}, 253(2):147--161, 1972.

\bibitem{Hairer3}
M.~Hairer.
\newblock Regularity structures and the dynamical {$\Phi^4_3$} model.
\newblock {\em arXiv:1508.05261}, pages 1--46, 2015.

\bibitem{KirrmannSchneiderMielke}
P.~Kirrmann, G.~Schneider, and A.~Mielke.
\newblock The validity of modulation equations for extended systems with cubic
  nonlinearities.
\newblock {\em Proc. R. Soc. Edinburgh A}, 122(1):85--91, 1992.

\bibitem{FiedlerRocha}
B.~Fiedler and C.~Rocha.
\newblock Heteroclinic orbits of semilinear parabolic equations.
\newblock {\em J. Differen. Equat.}, 125(1):239--281, 1996.

\bibitem{FiedlerScheel}
Bernold Fiedler and Arnd Scheel.
\newblock Spatio-temporal dynamics of reaction-diffusion patterns.
\newblock In {\em Trends in nonlinear analysis}, pages 23--152. Springer,
  Berlin, 2003.

\bibitem{Sandstede2}
B.~Sandstede.
\newblock Stability of {$N$}-fronts bifurcating from a twisted heteroclinic
  loop and an application to the {FitzHugh-Nagumo} equation.
\newblock {\em SIAM J. Math. Anal.}, 29(1):183--207, 1998.

\bibitem{BeynThuemmler}
W.J. Beyn and V.~Th\"ummler.
\newblock Freezing solutions of equivariant evolution equations.
\newblock {\em SIAM J. Appl. Dyn. Syst.}, 3(2):85--116, 2004.

\bibitem{MetzlerKlafter}
R.~Metzler and J.~Klafter.
\newblock The random walk's guide to anomalous diffusion: a fractional dynamics
  approach.
\newblock {\em Phys. Rep.}, 339(1):1--77, 2000.

\bibitem{samko93}
Stefan~G. Samko, Anatoly~A. Kilbas, and Oleg~I. Marichev.
\newblock {\em Fractional integrals and derivatives}.
\newblock Gordon and Breach Science Publishers, Yverdon, 1993.
\newblock Theory and applications, Edited and with a foreword by S. M.
  Nikol{'}ski\u{\i}, Translated from the 1987 Russian original, Revised by the
  authors.

\bibitem{BeScMeWh01}
David~A. Benson, Rina Schumer, Mark~M. Meerschaert, and Stephen~W. Wheatcraft.
\newblock Fractional dispersion, {L}\'{e}vy motion, and the {MADE} tracer
  tests.
\newblock {\em Transp. Porous Media}, 42(1-2):211--240, 2001.

\bibitem{DelCastilloNegrete:1998}
D.~Del Castillo-Negrete.
\newblock {Asymmetric transport and non-Gaussian statistics of passive scalars
  in vortices in shear}.
\newblock {\em Physics of Fluids}, 10(3):576, 1998.

\bibitem{dCNCaLy04}
D.~del Castillo-Negrete, B.~A. Carreras, and V.~E. Lynch.
\newblock {Fractional diffusion in plasma turbulence}.
\newblock {\em Physics of Plasmas}, 11(8):3854--3864, 07 2004.

\bibitem{Hanert+etal:2011}
E.~Hanert, E.~Schumacher, and E.~Deleersnijder.
\newblock {Front dynamics in fractional-order epidemic models.}
\newblock {\em Journal of theoretical biology}, 279(1):9--16, June 2011.

\bibitem{MeJeChBa14}
Ralf Metzler, Jae-Hyung Jeon, Andrey~G. Cherstvy, and Eli Barkai.
\newblock Anomalous diffusion models and their properties: non-stationarity{,}
  non-ergodicity{,} and ageing at the centenary of single particle tracking.
\newblock {\em Phys. Chem. Chem. Phys.}, 16:24128--24164, 2014.

\bibitem{KlaRadSok08}
Rainer Klages, G.~(Günter) Radons, and Igor~M. Sokolov.
\newblock {\em Anomalous transport : foundations and applications}.
\newblock WILEY-VCH Verlag GmbH \& Co. KGaA, Weinheim, Germany, 2008.

\bibitem{ScMeBa09}
Rina Schumer, Mark~M. Meerschaert, and Boris Baeumer.
\newblock Fractional advection-dispersion equations for modeling transport at
  the earth surface.
\newblock {\em Journal of Geophysical Research: Earth Surface}, 114(F4), 2009.

\bibitem{delia-etal20}
Marta D'Elia, Qiang Du, Christian Glusa, Max Gunzburger, Xiaochuan Tian, and
  Zhi Zhou.
\newblock Numerical methods for nonlocal and fractional models.
\newblock {\em Acta Numer.}, 29:1--124, 2020.

\bibitem{Mendez+etal:2010}
V.~M{\'e}ndez, S.~Fedotov, and W.~Horsthemke.
\newblock {\em Reaction-transport systems}.
\newblock Springer Series in Synergetics. Springer, Heidelberg, 2010.
\newblock Mesoscopic foundations, fronts, and spatial instabilities.

\bibitem{Sato:1999}
K.~Sato.
\newblock {\em L\'evy processes and infinitely divisible distributions},
  volume~68 of {\em Cambridge Studies in Advanced Mathematics}.
\newblock Cambridge University Press, Cambridge, 1999.
\newblock Translated from the 1990 Japanese original, Revised by the author.

\bibitem{DiROSeVa22}
Serena Dipierro, Xavier Ros-Oton, Joaquim Serra, and Enrico Valdinoci.
\newblock Non-symmetric stable operators: regularity theory and integration by
  parts.
\newblock {\em Adv. Math.}, 401:Paper No. 108321, 100, 2022.

\bibitem{kwasnicki17}
Mateusz Kwa\'snicki.
\newblock Ten equivalent definitions of the fractional {L}aplace operator.
\newblock {\em Fract. Calc. Appl. Anal.}, 20(1):7--51, 2017.

\bibitem{mclean00}
William McLean.
\newblock {\em Strongly elliptic systems and boundary integral equations}.
\newblock Cambridge University Press, Cambridge, 2000.

\bibitem{diethelm04}
Kai Diethelm.
\newblock {\em The analysis of fractional differential equations}, volume 2004
  of {\em Lecture Notes in Mathematics}.
\newblock Springer-Verlag, Berlin, 2010.
\newblock An application-oriented exposition using differential operators of
  Caputo type.

\bibitem{Mendez2014}
Vicen\c{c} M\'{e}ndez, Daniel Campos, and Frederic Bartumeus.
\newblock {\em Stochastic Foundations in Movement Ecology: Anomalous Diffusion,
  Front Propagation and Random Searches}.
\newblock Springer Series in Synergetics. Springer-Verlag Berlin Heidelberg,
  1st edition, 2014.

\bibitem{Jin21}
Bangti Jin.
\newblock {\em Fractional differential equations---an approach via fractional
  derivatives}, volume 206 of {\em Applied Mathematical Sciences}.
\newblock Springer, Cham, 2021.

\bibitem{SK2005}
I.~M. Sokolov and J.~Klafter.
\newblock From diffusion to anomalous diffusion: A century after {Einstein’s
  Brownian} motion.
\newblock {\em Chaos: An Interdisciplinary Journal of Nonlinear Science},
  15(2):026103, 06 2005.

\bibitem{Zacher2019}
Rico Zacher.
\newblock Time fractional diffusion equations: solution concepts, regularity,
  and long-time behavior.
\newblock In Anatoly Kochubei and Yuri Luchko, editors, {\em Volume 2
  Fractional Differential Equations}, pages 159--180. De Gruyter, Berlin,
  Boston, 2019.

\bibitem{SW1989}
W.~R. Schneider and W.~Wyss.
\newblock Fractional diffusion and wave equations.
\newblock {\em Journal of Mathematical Physics}, 30(1):134--144, 01 1989.

\bibitem{MMPG2007}
Francesco Mainardi, Antonio Mura, Gianni Pagnini, and Rudolf Gorenflo.
\newblock Sub-diffusion equations of fractional order and their fundamental
  solutions.
\newblock In K.~Ta{\c{s}}, J.~A. Tenreiro~Machado, and D.~Baleanu, editors,
  {\em Mathematical Methods in Engineering}, pages 23--55. Springer
  Netherlands, Dordrecht, 2007.

\bibitem{Metzler2004}
Ralf Metzler and Joseph Klafter.
\newblock The restaurant at the end of the random walk: recent developments in
  the description of anomalous transport by fractional dynamics.
\newblock {\em Journal of Physics A: Mathematical and General},
  37(31):R161--R208, 2004.

\bibitem{Matignon96}
Denis Matignon.
\newblock Stability results for fractional differential equations with
  applications to control processing.
\newblock In {\em Computational engineering in systems applications}, volume~2,
  pages 963--968. 1996.

\bibitem{BGK21}
Oana Brandibur, Roberto Garrappa, and Eva Kaslik.
\newblock Stability of systems of fractional-order differential equations with
  {Caputo} derivatives.
\newblock {\em Mathematics}, 9(8), 2021.

\bibitem{Mainardi2001}
Francesco Mainardi, Yuri Luchko, and Gianni Pagnini.
\newblock The fundamental solution of the space-time fractional diffusion
  equation.
\newblock {\em Fractional Calculus and Applied Analysis}, 4(2):153--192, 2001.

\bibitem{Kilbas2006}
A.~A. Kilbas, H.~M. Srivastava, and J.~J. Trujillo.
\newblock {\em Theory and Applications of Fractional Differential Equations},
  volume 204 of {\em North-Holland Mathematics Studies}.
\newblock Elsevier Science, Amsterdam, 1st edition, 2006.

\bibitem{Kochubei1990}
A.~N. Kochubei.
\newblock Diffusion of fractional order.
\newblock {\em Differential Equations}, 26(4):485--492, 1990.

\bibitem{EK2004}
Samuil~D. Eidelman and Anatoly~N. Kochubei.
\newblock Cauchy problem for fractional diffusion equations.
\newblock {\em Journal of Differential Equations}, 199(2):211--255, 2004.

\bibitem{YR21}
Jichen Yang and Jens D.~M. Rademacher.
\newblock Reaction-subdiffusion systems and memory: spectra, {T}uring
  instability and decay estimates.
\newblock {\em IMA J. Appl. Math.}, 86(2):247--293, 2021.

\bibitem{MMAK2019}
William McLean, Kassem Mustapha, Raed Ali, and Omar Knio.
\newblock Well-posedness of time-fractional advection-diffusion-reaction
  equations.
\newblock {\em Fract. Calc. Appl. Anal.}, 22(4):918--944, 2019.

\bibitem{GLM2000}
Rudolf Gorenflo, Yuri Luchko, and Francesco Mainardi.
\newblock Wright functions as scale-invariant solutions of the diffusion-wave
  equation.
\newblock {\em Journal of Computational and Applied Mathematics},
  118(1):175--191, 2000.

\bibitem{MMP2010}
Francesco Mainardi, Antonio Mura, and Gianni Pagnini.
\newblock The {M-Wright} function in time-fractional diffusion processes: a
  tutorial survey.
\newblock {\em International Journal of Differential Equations}, 2010:1--29,
  2010.

\bibitem{Vergara2015}
Vicente Vergara and Rico Zacher.
\newblock Optimal decay estimates for time-fractional and other nonlocal
  subdiffusion equations via energy methods.
\newblock {\em SIAM Journal on Mathematical Analysis}, 47(1):210--239, 2015.

\bibitem{KSVZ2016}
Jukka Kemppainen, Juhana Siljander, Vicente Vergara, and Rico Zacher.
\newblock Decay estimates for time-fractional and other non-local in time
  subdiffusion equations in {$\mathbb{R}^d$}.
\newblock {\em Math. Ann.}, 366(3-4):941--979, 2016.

\bibitem{Henry2006}
B.~I. Henry, T.~A.~M. Langlands, and S.~L. Wearne.
\newblock Anomalous diffusion with linear reaction dynamics: from continuous
  time random walks to fractional reaction-diffusion equations.
\newblock {\em Physical Review E}, 74(3):031116, 2006.

\bibitem{SWT03}
Kazuhiko Seki, Mariusz Wojcik, and M~Tachiya.
\newblock Fractional reaction-diffusion equation.
\newblock {\em The Journal of chemical physics}, 119(4):2165--2170, 2003.

\bibitem{Seki2003}
Kazuhiko Seki, Mariusz Wojcik, and M~Tachiya.
\newblock Recombination kinetics in subdiffusive media.
\newblock {\em The Journal of chemical physics}, 119(14):7525--7533, 2003.

\bibitem{Yuste2004}
S.~B. Yuste, L.~Acedo, and Katja Lindenberg.
\newblock Reaction front in an {A + B $\to$ C} reaction-subdiffusion process.
\newblock {\em Physical Review E}, 69(3):036126, 2004.

\bibitem{LHW2007}
T.~A.~M. Langlands, B.~I. Henry, and S.~L. Wearne.
\newblock Turing pattern formation with fractional diffusion and fractional
  reactions.
\newblock {\em Journal of Physics: Condensed Matter}, 19(6):065115, 2007.

\bibitem{Nec2008}
Y.~Nec and A.~A. Nepomnyashchy.
\newblock Turing instability of anomalous reaction--anomalous diffusion
  systems.
\newblock {\em European Journal of Applied Mathematics}, 19(3):329--349, 2008.

\bibitem{VeZa17}
V.~Vergara and R.~Zacher.
\newblock Stability, instability, and blowup for time fractional and other
  nonlocal in time semilinear subdiffusion equations.
\newblock {\em J.~Evol.~Equ.}, 17:599--626, 2017.

\bibitem{KTT2020}
Tran~Dinh Ke, Nguyen~Nhu Thang, and Lam Tran~Phuong Thuy.
\newblock Regularity and stability analysis for a class of semilinear nonlocal
  differential equations in {Hilbert} spaces.
\newblock {\em Journal of Mathematical Analysis and Applications},
  483(2):123655, 2020.

\bibitem{DKT2023}
Nguyen~Van Dac, Tran~Dinh Ke, and Lam Tran~Phuong Thuy.
\newblock On stability and regularity for semilinear anomalous diffusion
  equations perturbed by weak-valued nonlinearities.
\newblock {\em Discrete and Continuous Dynamical Systems - S}, pages 1--19,
  2023.

\bibitem{NTN23}
Nhu~Thang Nguyen, Dinh~Ke Tran, and Van~Dac Nguyen.
\newblock Stability analysis for nonlocal evolution equations involving
  infinite delays.
\newblock {\em J. Fixed Point Theory Appl.}, 25(1):Paper No. 22, 33, 2023.

\bibitem{KP23}
Uttam Kumar and Subramaniam Pushpavanam.
\newblock The effect of subdiffusion on the stability of autocatalytic systems.
\newblock {\em Chemical Engineering Science}, 265:118230, 2023.

\bibitem{Henry2000}
Bruce~I. Henry and Susan~L. Wearne.
\newblock Fractional reaction--diffusion.
\newblock {\em Physica A: Statistical Mechanics and its Applications},
  276(3):448--455, 2000.

\bibitem{Fedotov2010}
Sergei Fedotov.
\newblock Non-markovian random walks and nonlinear reactions: Subdiffusion and
  propagating fronts.
\newblock {\em Physical Review E}, 81(1):011117, 2010.

\bibitem{Nec2007}
Y.~Nec and A.~A. Nepomnyashchy.
\newblock Linear stability of fractional reaction-diffusion systems.
\newblock {\em Mathematical Modelling of Natural Phenomena}, 2(2):77--105,
  2007.

\bibitem{Nepomnyashchy2016}
A.~A. Nepomnyashchy.
\newblock Mathematical modelling of subdiffusion-reaction systems.
\newblock {\em Mathematical Modelling of Natural Phenomena}, 11(1):26--36,
  2016.

\bibitem{Henry2002}
B.~I. Henry and S.~L. Wearne.
\newblock Existence of {Turing} instabilities in a two-species fractional
  reaction-diffusion system.
\newblock {\em SIAM Journal on Applied Mathematics}, 62(3):870--887, 2002.

\bibitem{Langlands2008}
T.~A.~M. Langlands, B.~I. Henry, and S.~L. Wearne.
\newblock Anomalous subdiffusion with multispecies linear reaction dynamics.
\newblock {\em Physical Review E}, 77(2):021111, 2008.

\bibitem{AL21}
Amanda~M. Alexander and Sean~D. Lawley.
\newblock Reaction-subdiffusion equations with species-dependent movement.
\newblock {\em SIAM J. Appl. Math.}, 81(6):2457--2479, 2021.

\bibitem{Sandstede2002}
Bj\"orn Sandstede.
\newblock Chapter 18 - stability of travelling waves.
\newblock In Bernold Fiedler, editor, {\em Handbook of Dynamical Systems},
  volume~2 of {\em Handbook of Dynamical Systems}, pages 983 -- 1055. Elsevier
  Science, 2002.

\bibitem{Henry2005}
B.~I. Henry, T.~A.~M. Langlands, and S.~L. Wearne.
\newblock Turing pattern formation in fractional activator-inhibitor systems.
\newblock {\em Physical Review E}, 72(2):026101, 2005.

\bibitem{Komura}
Y.~K\=omura.
\newblock Nonlinear semi-groups in {H}ilbert space.
\newblock {\em J.~Math.~Soc.~Japan}, 19:493--507, 1967.

\bibitem{CL}
M.G. Crandall and T.M. Liggett.
\newblock Generation of semi-groups of nonlinear transformations on general
  {B}anach spaces.
\newblock {\em Amer.~J.~Math.}, 93:265--298, 1971.

\bibitem{HB1}
H.~Br\'{e}zis.
\newblock {\em Op\'{e}rateurs maximaux monotones et semi-groupes de
  contractions dans les espaces de {H}ilbert}, volume No. 5 of {\em
  North-Holland Mathematics Studies}.
\newblock North-Holland Publishing Co., Amsterdam-London; American Elsevier
  Publishing Co., Inc., New York, 1973.
\newblock Notas de Matem\'{a}tica, No. 50. [Mathematical Notes].

\bibitem{B75}
V.~Barbu.
\newblock Nonlinear {V}olterra equations in a {H}ilbert space.
\newblock {\em SIAM J.~Math.~Anal.}, 6:728--741, 1975.

\bibitem{B77}
V.~Barbu.
\newblock On a nonlinear {V}olterra integral equation on a {H}ilbert space.
\newblock {\em SIAM J.~Math.~Anal.}, 8:346--355, 1977.

\bibitem{B-dn}
V.~Barbu.
\newblock Existence for nonlinear {V}olterra equations in {H}ilbert spaces.
\newblock {\em SIAM J.~Math.~Anal.}, pages 552--569, 1979.

\bibitem{CrNo78}
M.G. Crandall and J.A. Nohel.
\newblock An abstract functional differential equation and a related nonlinear
  {V}olterra equation.
\newblock {\em Israel J.~Math.}, 29:313--328, 1978.

\bibitem{Gri78}
G.~Gripenberg.
\newblock An existence result for a nonlinear {V}olterra integral equation in a
  {H}ilbert space.
\newblock {\em SIAM J.~Math.~Anal.}, 9:793--805, 1978.

\bibitem{LoSt78}
S.-O. Londen and O.J. Staffans.
\newblock A note on {V}olterra equations in a {H}ilbert space.
\newblock {\em Proc.~Amer.~Math.~Soc.}, 70:57--62, 1978.

\bibitem{KiSt79}
T.~Kiffe and M.~Stecher.
\newblock An abstract {V}olterra integral equation in a reflexive {B}anach
  space.
\newblock {\em J.~Differential Equations}, 34:303--325, 1979.

\bibitem{Cl80}
Ph. Cl\'ement.
\newblock On abstract {V}olterra equations with kernels having a positive
  resolvent.
\newblock {\em Israel J.~Math.}, 36:193--200, 1980.

\bibitem{Aiz93}
S.~Aizicovici.
\newblock An abstract doubly nonlinear {V}olterra integral equation.
\newblock {\em Funkcialaj Ekvacioj}, 36:479--497, 1993.

\bibitem{ClNo81}
Ph. Cl\'ement and J.A. Nohel.
\newblock Asymptotic behavior of solutions of nonlinear volterra equations with
  completely positive kernels.
\newblock {\em SIAM J.~Math.~Anal.}, 12:514--535, 1981.

\bibitem{ClNo79}
Ph. Cl\'ement and J.A. Nohel.
\newblock Abstract linear and nonlinear {V}olterra equations preserving
  positivity.
\newblock {\em SIAM J.~Math.~Anal.}, 10:365--388, 1979.

\bibitem{Gri79}
G.~Gripenberg.
\newblock An abstract nonlinear volterra equation.
\newblock {\em Israel J.~Math.}, 34:198--212, 1979.

\bibitem{Ki80}
T.~Kiffe.
\newblock A perturbation of an abstract {V}olterra equation.
\newblock {\em SIAM J.~Math.~Anal.}, 11:1036--1046, 1980.

\bibitem{Cl84}
Ph. Cl\'ement.
\newblock {\em On abstract {V}olterra equations in {B}anach spaces with
  completely positive kernels}, volume 1076 of {\em Lecture Notes in Math.},
  pages 32--40.
\newblock Springer, Berlin, 1984.

\bibitem{Gri85}
G.~Gripenberg.
\newblock Volterra integro-differential equations with accretive nonlinearity.
\newblock {\em J.~Differential Equations}, 60:57--79, 1985.

\bibitem{CP90}
Ph. Cl\'ement and J.~Pr\"uss.
\newblock Completely positive measures and {F}eller semigroups.
\newblock {\em Math.~Ann.}, 287:73--105, 1990.

\bibitem{Za09}
R.~Zacher.
\newblock Weak solutions of abstract evolutionary integro-differential
  equations in {H}ilbert spaces.
\newblock {\em Funkcial.~ Ekvac.}, 52:1--18, 2009.

\bibitem{VeZa08}
V.~Vergara and R.~Zacher.
\newblock Lyapunov functions and convergence to steady state for differential
  equations of fractional order.
\newblock {\em Math.~Z.}, 259:287--309, 2008.

\bibitem{Za08}
R.~Zacher.
\newblock Boundedness of weak solutions to evolutionary partial
  integro-differential equations with discontinuous coefficients.
\newblock {\em J.~Math.~Anal.~Appl.}, 348:137--149, 2008.

\bibitem{Za13}
R.~Zacher.
\newblock A {D}e {G}iorgi–{N}ash type theorem for time fractional diffusion
  equations.
\newblock {\em Math.~Ann.}, 356:99--146, 2013.

\bibitem{VeZa10}
V.~Vergara and R.~Zacher.
\newblock A priori bounds for degenerate and singular evolutionary partial
  integro-differential equations.
\newblock {\em Nonlinear Anal.}, 73:3572--3585, 2010.

\bibitem{Za10}
R.~Zacher.
\newblock A weak {H}arnack inequality for fractional evolution equations with
  discontinuous coefficients.
\newblock {\em Ann.~Sc.~Norm.~Super.~Pisa Cl.~Sci.~(5)}, 12:903--940, 2013.

\bibitem{AV19}
E.~Affili and E.~Valdinoci.
\newblock Decay estimates for evolution equations with classical and fractional
  time-derivatives.
\newblock {\em J.~Differential Equations}, 266:4027--4060, 2019.

\bibitem{Za19}
R.~Zacher.
\newblock {\em Time fractional diffusion equations: solution concepts,
  regularity, and long-time behavior}, volume~2, pages 159--179.
\newblock De {G}ruyter, Berlin, 2019.

\bibitem{Za05}
R.~Zacher.
\newblock Maximal regularity of type $l_p$ for abstract parabolic {V}olterra
  equations.
\newblock {\em J.~Evol.~Equ.}, 5:79--103, 2005.

\bibitem{Za06}
R.~Zacher.
\newblock Quasilinear parabolic integro-differential equations with nonlinear
  boundary conditions.
\newblock {\em Differential Integral Equations}, 19:1129--1156, 2006.

\bibitem{GW20}
C.G. Gal and M.~Warma.
\newblock {\em Fractional-in-time semilinear parabolic equations and
  applications}, volume~84 of {\em Mathematics \& Applications}.
\newblock Springer, Cham, 2020.

\bibitem{KRY20}
A.~Kubica, K.~Ryszewska, and M.~Yamamoto.
\newblock {\em Time-fractional differential equations—a theoretical
  introduction}.
\newblock Springer Briefs Math. Springer, Singapore, 2020.

\bibitem{A19}
G.~Akagi.
\newblock Fractional flows driven by subdifferentials in {H}ilbert spaces.
\newblock {\em Israel J.~Math.}, 234:809--862, 2019.

\bibitem{AB23}
C.O. Alves and T.~Boudjeriou.
\newblock Existence and uniqueness of solution for some time fractional
  parabolic equations involving the 1-{L}aplace operator.
\newblock {\em Partial Differential Equations and Applications}, 4, 2023.

\bibitem{WWZ21}
P.~Wittbold, P.~Wolejko, and R.~Zacher.
\newblock Bounded weak solutions of time-fractional porous medium type and more
  general nonlinear and degenerate evolutionary integro-differential equations.
\newblock {\em J.~Math.~Anal.~Appl.}, 499:20 pp, 2021.

\bibitem{capella2010regularity}
A.~Capella, J.~D{\'a}vila, L.~Dupaigne, and Y.~Sire.
\newblock {Regularity of radial extremal solutions for some non-local
  semilinear equations}.
\newblock {\em Communications in Partial Differential Equations},
  36(8):1353--1384, 2011.

\bibitem{abatangelo2019getting}
Nicola Abatangelo and Enrico Valdinoci.
\newblock Getting acquainted with the fractional {L}aplacian.
\newblock {\em Contemporary research in elliptic PDEs and related topics},
  pages 1--105, 2019.

\bibitem{servadei2014spectrum}
R.~Servadei and E.~Valdinoci.
\newblock On the spectrum of two different fractional operators.
\newblock {\em Proceedings of the Royal Society of Edinburgh Section A:
  Mathematics}, 144(4):831--855, 2014.

\bibitem{bucur2016nonlocal}
Claudia Bucur, Enrico Valdinoci, et~al.
\newblock {\em Nonlocal diffusion and applications}, volume~20.
\newblock Springer, 2016.

\bibitem{duo2019comparative}
Siwei Duo, Hong Wang, and Yanzhi Zhang.
\newblock A comparative study on nonlocal diffusion operators related to the
  fractional {L}aplacian.
\newblock {\em Discrete \& Continuous Dynamical Systems-Series B}, 24(1), 2019.

\bibitem{musina2014fractional}
Roberta Musina and Alexander~I Nazarov.
\newblock On fractional laplacians.
\newblock {\em Communications in partial differential equations},
  39(9):1780--1790, 2014.

\bibitem{owolabi2018modelling}
Kolade~M Owolabi and Abdon Atangana.
\newblock Modelling and formation of spatiotemporal patterns of fractional
  predation system in subdiffusion and superdiffusion scenarios.
\newblock {\em The European physical Journal Plus}, 133:1--13, 2018.

\bibitem{tian2015turing}
Canrong Tian.
\newblock Turing pattern formation in a semiarid vegetation model with
  fractional-in-space diffusion.
\newblock {\em Bulletin of mathematical biology}, 77:2072--2085, 2015.

\bibitem{Somathilake2018}
L.W. Somathilake and K.~Burrage.
\newblock A space-fractional-reaction-diffusion model for pattern formation in
  coral reefs.
\newblock {\em Cogent Mathematics \& Statistics}, 5(1):1426524, 2018.

\bibitem{golovin2008turing}
Alexander~A Golovin, Bernard~J Matkowsky, and Vladimir~A Volpert.
\newblock Turing pattern formation in the {B}russelator model with
  superdiffusion.
\newblock {\em SIAM Journal on Applied Mathematics}, 69(1):251--272, 2008.

\bibitem{zhang2014turing}
Lai Zhang and Canrong Tian.
\newblock Turing pattern dynamics in an activator-inhibitor system with
  superdiffusion.
\newblock {\em Physical Review E}, 90(6):062915, 2014.

\bibitem{schnakenberg1979simple}
J.~Schnakenberg.
\newblock Simple chemical reaction systems with limit cycle behaviour.
\newblock {\em Journal of Theoretical Biology}, 81(3):389--400, 1979.

\bibitem{alfaro2008singular}
M.~Alfaro, D.~Hilhorst, and H.~Matano.
\newblock The singular limit of the {A}llen--{C}ahn equation and the
  {F}itz{H}ugh--{N}agumo system.
\newblock {\em Journal of Differential Equations}, 245(2):505--565, 2008.

\bibitem{alama1997stationary}
S.~Alama, L.~Bronsard, and C.~Gui.
\newblock Stationary layered solutions in $\mathbb{R}^2$ for an {A}llen--{C}ahn
  system with multiple well potential.
\newblock {\em Calculus of Variations and Partial Differential Equations},
  5(4):359--390, 1997.

\bibitem{rabinowitz2003mixed}
P.H. Rabinowitz and E.W. Stredulinsky.
\newblock Mixed states for an {A}llen--{C}ahn type equation.
\newblock {\em Communications on Pure and Applied Mathematics},
  56(8):1078--1134, 2003.

\bibitem{ward1996metastable}
M.J. Ward.
\newblock Metastable bubble solutions for the {A}llen--{C}ahn equation with
  mass conservation.
\newblock {\em SIAM Journal on Applied Mathematics}, 56(5):1247--1279, 1996.

\bibitem{aranson2002world}
I.~S Aranson and L.~Kramer.
\newblock The world of the complex {G}inzburg--{L}andau equation.
\newblock {\em Reviews of Modern Physics}, 74(1):99--143, 2002.

\bibitem{mielke2002ginzburg}
A.~Mielke.
\newblock The {G}inzburg--{L}andau equation in its role as a modulation
  equation.
\newblock In {\em Handbook of Dynamical Systems}, volume~2, pages 759--834.
  Elsevier, 2002.

\bibitem{schneider1996validity}
G.~Schneider.
\newblock The validity of generalized {G}inzburg--{L}andau equations.
\newblock {\em Mathematical Methods in the Applied Sciences}, 19(9):717--736,
  1996.

\bibitem{kapitula1998instability}
T.~Kapitula and B.~Sandstede.
\newblock Instability mechanism for bright solitary-wave solutions to the
  cubic--quintic {G}inzburg--{L}andau equation.
\newblock {\em Journal of the Optical Society of America B}, 15(11):2757--2762,
  1998.

\bibitem{kuehn2015numerical}
C.~Kuehn.
\newblock Numerical continuation and {SPDE} stability for the 2{D}
  cubic--quintic {A}llen--{C}ahn equation.
\newblock {\em SIAM/ASA Journal on Uncertainty Quantification}, 3(1):762--789,
  2015.

\bibitem{AKMR21}
Franz Achleitner, Christian Kuehn, Jens~M. Melenk, and Alexander Rieder.
\newblock Metastable speeds in the fractional {A}llen-{C}ahn equation.
\newblock {\em Appl. Math. Comput.}, 408:Paper No. 126329, 18, 2021.

\bibitem{akagi2016fractional}
G.~Akagi, G.~Schimperna, and A.~Segatti.
\newblock Fractional {C}ahn--{H}illiard, {A}llen--{C}ahn and porous medium
  equations.
\newblock {\em Journal of Differential Equations}, 261(6):2935--2985, 2016.

\bibitem{bueno2014fourier}
A.~Bueno-Orovio, D.~Kay, and K.~Burrage.
\newblock Fourier spectral methods for fractional-in-space reaction-diffusion
  equations.
\newblock {\em BIT Numerical mathematics}, 54(4):937--954, 2014.

\bibitem{ehstand2021numerical}
N.~Ehstand, C.~Kuehn, and C.~Soresina.
\newblock Numerical continuation for fractional pdes: sharp teeth and bloated
  snakes.
\newblock {\em Communications in Nonlinear Science and Numerical Simulation},
  98:105762, 2021.

\bibitem{swift1977hydrodynamic}
J.~Swift and P.C. Hohenberg.
\newblock Hydrodynamic fluctuations at the convective instability.
\newblock {\em Physical Review A}, 15(1):319, 1977.

\bibitem{Schneider5}
G.~Schneider.
\newblock Diffusive stability of spatial periodic solutions of the
  {Swift--Hohenberg} equation.
\newblock {\em Communications in Mathematical Physics}, 178(3):679--702, 1996.

\bibitem{Thieleetal}
U.~Thiele, A.J. Archer, M.J. Robbins, H.~Gomez, and E.~Knobloch.
\newblock Localized states in the conserved {Swift--Hohenberg} equation with
  cubic nonlinearity.
\newblock {\em Physical Review E}, 87(4):042915, 2013.

\bibitem{ColletEckmann1}
P.~Collet and J.P. Eckmann.
\newblock The time dependent amplitude equation for the {Swift--Hohenberg}
  problem.
\newblock {\em Communications in Mathematical Physics}, 132(1):139--153, 1990.

\bibitem{vanHarten}
A.~van Harten.
\newblock On the validity of the {Ginzburg--Landau} equation.
\newblock {\em Journal of Nonlinear Science}, 1(4):397--422, 1991.

\bibitem{KuehnThrom}
C.~Kuehn and S.~Throm.
\newblock Validity of amplitude equations for non-local non-linearities.
\newblock {\em Journal of Mathematical Physics}, 59:071510, 2018.

\bibitem{MorganDawes}
D.~Morgan and J.H.P. Dawes.
\newblock The {Swift--Hohenberg} equation with a nonlocal nonlinearity.
\newblock {\em Physica D}, 270:60--80, 2014.

\bibitem{Uecker2014}
H.~Uecker and D.~Wetzel.
\newblock Numerical results for snaking of patterns over patterns in some 2{D}
  {Selkov}--{Schnakenberg} reaction-diffusion systems.
\newblock {\em {SIAM} Journal on Applied Dynamical Systems}, 13(1):94--128,
  2014.

\bibitem{VolNecNep13}
V.~A. Volpert, Y.~Nec, and A.~A. Nepomnyashchy.
\newblock Fronts in anomalous diffusion-reaction systems.
\newblock {\em Philos. Trans. R. Soc. Lond. Ser. A Math. Phys. Eng. Sci.},
  371(1982):20120179, 18, 2013.

\bibitem{NecVolNep10}
Yana Nec, Vladimir~A. Volpert, and Alexander~A. Nepomnyashchy.
\newblock Front propagation problems with sub-diffusion.
\newblock {\em Discrete Contin. Dyn. Syst.}, 27(2):827--846, 2010.

\bibitem{Ishii23_arxiv}
Hiroshi Ishii.
\newblock Propagating front solutions in a time-fractional fisher-kpp equation,
  2023.

\bibitem{AchCueHit14}
F.~Achleitner, C.~M. Cuesta, and S.~Hittmeir.
\newblock Travelling waves for a non-local {K}orteweg--de {V}ries--{B}urgers
  equation.
\newblock {\em J. Differential Equations}, 257(3):720--758, 2014.

\bibitem{CueAch17}
C.~M. Cuesta and F.~Achleitner.
\newblock Addendum to ``{T}ravelling waves for a non-local {K}orteweg--de
  {V}ries-{B}urgers equation'' [{J}. {D}ifferential {E}quations 257 (3) (2014)
  720--758] [{MR}3208089].
\newblock {\em J. Differential Equations}, 262(2):1155--1160, 2017.

\bibitem{AK15AC}
Franz Achleitner and Christian Kuehn.
\newblock Traveling waves for a bistable equation with nonlocal diffusion.
\newblock {\em Adv. Differential Equations}, 20(9-10):887--936, 2015.

\bibitem{Aronson+Weinberger:1974}
D.G. Aronson and H.F. Weinberger.
\newblock Nonlinear diffusion in population genetics, combustion, and nerve
  pulse propagation.
\newblock In {\em Partial Differential Equations and Related Topics}, volume
  446 of {\em Lecture Notes in Mathematics}, pages 5--49. Springer, 1974.

\bibitem{Fife+McLeod:1977}
Paul~C. Fife and J.~B. McLeod.
\newblock The approach of solutions of nonlinear diffusion equations to
  travelling front solutions.
\newblock {\em Arch. Ration. Mech. Anal.}, 65(4):335--361, 1977.

\bibitem{Cabre+Sire:2015}
Xavier Cabr{\'e} and Yannick Sire.
\newblock Nonlinear equations for fractional {L}aplacians {II}: {E}xistence,
  uniqueness, and qualitative properties of solutions.
\newblock {\em Trans. Amer. Math. Soc.}, 367(2):911--941, 2015.

\bibitem{Chmaj:2013}
A.~Chmaj.
\newblock {Existence of traveling waves in the fractional bistable equation}.
\newblock {\em Archiv der Mathematik}, 100(5):473--480, May 2013.

\bibitem{Gui+Zhao:2014}
Changfeng Gui and Mingfeng Zhao.
\newblock Traveling wave solutions of {A}llen-{C}ahn equation with a fractional
  {L}aplacian.
\newblock {\em Ann. Inst. H. Poincar\'{e} C Anal. Non Lin\'{e}aire},
  32(4):785--812, 2015.

\bibitem{Palatucci+etal:2013}
Giampiero Palatucci, Ovidiu Savin, and Enrico Valdinoci.
\newblock Local and global minimizers for a variational energy involving a
  fractional norm.
\newblock {\em Ann. Mat. Pura Appl. (4)}, 192(4):673--718, 2013.

\bibitem{Nec+etal:2008}
Y.~Nec, A.A. Nepomnyashchy, and A.A. Golovin.
\newblock Front-type solutions of fractional {A}llen-{C}ahn equation.
\newblock {\em Phys. D}, 237(24):3237--3251, 2008.

\bibitem{Volpert+etal:2010}
V.~A. Volpert, Y.~Nec, and A.~A. Nepomnyashchy.
\newblock Exact solutions in front propagation problems with superdiffusion.
\newblock {\em Phys. D}, 239(3-4):134--144, 2010.

\bibitem{Zanette:1997}
D.H. Zanette.
\newblock {Wave fronts in bistable reactions with anomalous L\'{e}vy-flight
  diffusion}.
\newblock {\em Physical Review E}, 55(1):1181--1184, 1997.

\bibitem{McKean:1970}
H.P. McKean.
\newblock {Nagumo's equation}.
\newblock {\em Advances in mathematics}, 4:209--223, 1970.

\bibitem{MaNiWa19}
Luyi Ma, Hong-Tao Niu, and Zhi-Cheng Wang.
\newblock Global asymptotic stability of traveling waves to the {A}llen-{C}ahn
  equation with a fractional {L}aplacian.
\newblock {\em Commun. Pure Appl. Anal.}, 18(5):2457--2472, 2019.

\bibitem{HoYu22}
Hongmei Cheng and Rong Yuan.
\newblock The stability of traveling waves for {A}llen-{C}ahn equations with
  fractional {L}aplacian.
\newblock {\em Appl. Anal.}, 101(1):263--273, 2022.

\bibitem{Bates+etal:1997}
Peter~W. Bates, Paul~C. Fife, Xiaofeng Ren, and Xuefeng Wang.
\newblock Traveling waves in a convolution model for phase transitions.
\newblock {\em Arch. Rational Mech. Anal.}, 138(2):105--136, 1997.

\bibitem{Chen:1997}
Xinfu Chen.
\newblock Existence, uniqueness, and asymptotic stability of traveling waves in
  nonlocal evolution equations.
\newblock {\em Adv. Differential Equations}, 2(1):125--160, 1997.

\bibitem{caffarelli-silvestre07}
Luis Caffarelli and Luis Silvestre.
\newblock An extension problem related to the fractional {L}aplacian.
\newblock {\em Comm. Partial Differential Equations}, 32(7-9):1245--1260, 2007.

\bibitem{Kw22}
Mateusz Kwa\'{s}nicki.
\newblock Harmonic extension technique for non-symmetric operators with
  completely monotone kernels.
\newblock {\em Calc. Var. Partial Differential Equations}, 61(6):Paper No. 202,
  40, 2022.

\bibitem{We21}
Maria~G. Westdickenberg.
\newblock On the metastability of the 1-{$d$} {A}llen-{C}ahn equation.
\newblock {\em J. Dynam. Differential Equations}, 33(4):1853--1879, 2021.

\bibitem{AK15AC2}
Franz Achleitner and Christian Kuehn.
\newblock Analysis and numerics of traveling waves for asymmetric fractional
  reaction-diffusion equations.
\newblock {\em Commun. Appl. Ind. Math.}, 6(2):e--532, 25, 2015.

\bibitem{Xin00}
Jack Xin.
\newblock Front propagation in heterogeneous media.
\newblock {\em SIAM Rev.}, 42(2):161--230, 2000.

\bibitem{Ha16}
Fran\c{c}ois Hamel.
\newblock Bistable transition fronts in {$\mathbb{R}^N$}.
\newblock {\em Adv. Math.}, 289:279--344, 2016.

\bibitem{Ta21}
Masaharu Taniguchi.
\newblock {\em Traveling front solutions in reaction-diffusion equations},
  volume~39 of {\em MSJ Memoirs}.
\newblock Mathematical Society of Japan, Tokyo, 2021.

\bibitem{Xin92}
J.~X. Xin.
\newblock Multidimensional stability of traveling waves in a bistable
  reaction-diffusion equation. {I}.
\newblock {\em Comm. Partial Differential Equations}, 17(11-12):1889--1899,
  1992.

\bibitem{LevXin92}
C.~D. Levermore and J.~X. Xin.
\newblock Multidimensional stability of traveling waves in a bistable
  reaction-diffusion equation. {II}.
\newblock {\em Comm. Partial Differential Equations}, 17(11-12):1901--1924,
  1992.

\bibitem{Jon83}
Christopher K. R.~T. Jones.
\newblock Asymptotic behaviour of a reaction-diffusion equation in higher space
  dimensions.
\newblock {\em Rocky Mountain J. Math.}, 13(2):355--364, 1983.

\bibitem{Kap97}
Todd Kapitula.
\newblock Multidimensional stability of planar travelling waves.
\newblock {\em Trans. Amer. Math. Soc.}, 349(1):257--269, 1997.

\bibitem{Volpert+etal:1994}
A.I. Volpert, V.A. Volpert, and V.A. Volpert.
\newblock {\em Traveling wave solutions of parabolic systems}, volume 140 of
  {\em Translations of Mathematical Monographs}.
\newblock American Mathematical Society, Providence, RI, 1994.
\newblock Translated from the Russian manuscript by James F. Heyda.

\bibitem{WaTa23}
Wah Wah and Masaharu Taniguchi.
\newblock Traveling front solutions for perturbed reaction-diffusion equations.
\newblock {\em Math. J. Okayama Univ.}, 65:125--143, 2023.

\bibitem{CheGuoHamNinRoq07}
Xinfu Chen, Jong-Shenq Guo, Fran\c{c}ois Hamel, Hirokazu Ninomiya, and
  Jean-Michel Roquejoffre.
\newblock Traveling waves with paraboloid like interfaces for balanced bistable
  dynamics.
\newblock {\em Ann. Inst. H. Poincar\'{e} C Anal. Non Lin\'{e}aire},
  24(3):369--393, 2007.

\bibitem{BeHa12}
Henri Berestycki and Fran\c{c}ois Hamel.
\newblock Generalized transition waves and their properties.
\newblock {\em Comm. Pure Appl. Math.}, 65(5):592--648, 2012.

\bibitem{FarVal11}
Alberto Farina and Enrico Valdinoci.
\newblock Rigidity results for elliptic {PDE}s with uniform limits: an abstract
  framework with applications.
\newblock {\em Indiana Univ. Math. J.}, 60(1):121--141, 2011.

\bibitem{DipVal19}
Serena Dipierro and Enrico Valdinoci.
\newblock Long-range phase coexistence models: recent progress on the
  fractional {A}llen-{C}ahn equation.
\newblock In {\em Topics in applied analysis and optimisation---partial
  differential equations, stochastic and numerical analysis}, CIM Ser. Math.
  Sci., pages 121--138. Springer, Cham, [2019] \copyright 2019.

\bibitem{FarVal09}
Alberto Farina and Enrico Valdinoci.
\newblock The state of the art for a conjecture of {D}e {G}iorgi and related
  problems.
\newblock In {\em Recent progress on reaction-diffusion systems and viscosity
  solutions}, pages 74--96. World Sci. Publ., Hackensack, NJ, 2009.

\bibitem{Sav10}
O.~Savin.
\newblock Phase transitions, minimal surfaces and a conjecture of {D}e
  {G}iorgi.
\newblock In {\em Current developments in mathematics, 2009}, pages 59--113.
  Int. Press, Somerville, MA, 2010.

\bibitem{ChaWei18}
Hardy Chan and Juncheng Wei.
\newblock On {D}e {G}iorgi's conjecture: recent progress and open problems.
\newblock {\em Sci. China Math.}, 61(11):1925--1946, 2018.

\bibitem{GraRyz15}
I.~S. Gradshteyn and I.~M. Ryzhik.
\newblock {\em Table of integrals, series, and products}.
\newblock Elsevier/Academic Press, Amsterdam, eighth edition, 2015.
\newblock Translated from the Russian, Translation edited and with a preface by
  Daniel Zwillinger and Victor Moll.

\bibitem{ChaWei17}
Hardy Chan and Juncheng Wei.
\newblock Traveling wave solutions for bistable fractional {A}llen-{C}ahn
  equations with a pyramidal front.
\newblock {\em J. Differential Equations}, 262(9):4567--4609, 2017.

\bibitem{MaWan21}
Luyi Ma and Zhi-Cheng Wang.
\newblock On the existence of cylindrically symmetric traveling fronts of
  fractional {A}llen-{C}ahn equation in {$\mathbb R^3$}.
\newblock {\em Differential Integral Equations}, 34(9-10):467--490, 2021.

\bibitem{MaNiuWan22}
Luyi Ma, Hong-Tao Niu, and Zhi-Cheng Wang.
\newblock V-shaped traveling fronts of fractional {A}llen-{C}ahn equations.
\newblock {\em J. Math. Phys.}, 63(2):Paper No. 021507, 29, 2022.

\bibitem{CheYua19}
Hongmei Cheng and Rong Yuan.
\newblock The stability of the equilibria of the {A}llen-{C}ahn equation with
  fractional diffusion.
\newblock {\em Appl. Anal.}, 98(3):600--610, 2019.

\bibitem{DipSerVal20}
Serena Dipierro, Joaquim Serra, and Enrico Valdinoci.
\newblock Improvement of flatness for nonlocal phase transitions.
\newblock {\em Amer. J. Math.}, 142(4):1083--1160, 2020.

\bibitem{DipVal23}
Serena Dipierro and Enrico Valdinoci.
\newblock Some perspectives on (non)local phase transitions and minimal
  surfaces.
\newblock {\em Bull. Math. Sci.}, 13(1):Paper No. 2330001, 77, 2023.

\bibitem{handbook_vol_3}
George~Em Karniadakis, editor.
\newblock {\em Handbook of fractional calculus with applications. {V}ol. 3}.
\newblock De Gruyter, Berlin, 2019.
\newblock Numerical methods.

\bibitem{li-chen18}
Changpin Li and An~Chen.
\newblock Numerical methods for fractional partial differential equations.
\newblock {\em Int. J. Comput. Math.}, 95(6-7):1048--1099, 2018.

\bibitem{jin-zhou23}
Bangti Jin and Zhi Zhou.
\newblock {\em Numerical treatment and analysis of time-fractional evolution
  equations}, volume 214 of {\em Applied Mathematical Sciences}.
\newblock Springer, Cham, [2023] \copyright 2023.

\bibitem{stynes-handbook}
Martin Stynes.
\newblock Singularities.
\newblock In {\em Handbook of fractional calculus with applications. {V}ol. 3},
  pages 287--305. De Gruyter, Berlin, 2019.

\bibitem{stynes-oriordan-gracia17}
Martin Stynes, Eugene O'Riordan, and Jos\'{e}~Luis Gracia.
\newblock Error analysis of a finite difference method on graded meshes for a
  time-fractional diffusion equation.
\newblock {\em SIAM J. Numer. Anal.}, 55(2):1057--1079, 2017.

\bibitem{mclean-mustapha15}
William McLean and Kassem Mustapha.
\newblock Time-stepping error bounds for fractional diffusion problems with
  non-smooth initial data.
\newblock {\em J. Comput. Phys.}, 293:201--217, 2015.

\bibitem{brunner04}
Hermann Brunner.
\newblock {\em Collocation methods for {V}olterra integral and related
  functional differential equations}, volume~15 of {\em Cambridge Monographs on
  Applied and Computational Mathematics}.
\newblock Cambridge University Press, Cambridge, 2004.

\bibitem{lubich88-I}
C.~Lubich.
\newblock Convolution quadrature and discretized operational calculus. {I}.
\newblock {\em Numer. Math.}, 52(2):129--145, 1988.

\bibitem{lubich88-II}
C.~Lubich.
\newblock Convolution quadrature and discretized operational calculus. {II}.
\newblock {\em Numer. Math.}, 52(4):413--425, 1988.

\bibitem{banjai-sayas22}
L.~Banjai and F.-J. Sayas.
\newblock {\em Integral equation methods for evolutionary PDE}, volume~53 of
  {\em Springer Series in Computational Mathematics}.
\newblock Springer Verlag, 2022.

\bibitem{lubich-ostermann93}
Ch. Lubich and A.~Ostermann.
\newblock Runge-{K}utta methods for parabolic equations and convolution
  quadrature.
\newblock {\em Math. Comp.}, 60(201):105--131, 1993.

\bibitem{lubich04}
Christian Lubich.
\newblock Convolution quadrature revisited.
\newblock {\em BIT}, 44(3):503--514, 2004.

\bibitem{lopez-fernandez-guo21}
Jing Guo and Maria Lopez-Fernandez.
\newblock Generalized convolution quadrature for the fractional integral and
  fractional diffusion equations, 2023.

\bibitem{banjai-makridakis22}
Lehel Banjai and Charalambos~G. Makridakis.
\newblock A posteriori error analysis for approximations of time-fractional
  subdiffusion problems.
\newblock {\em Math. Comp.}, 91(336):1711--1737, 2022.

\bibitem{greengard-rokhlin97}
Leslie Greengard and Vladimir Rokhlin.
\newblock A new version of the fast multipole method for the {L}aplace equation
  in three dimensions.
\newblock In {\em Acta numerica, 1997}, volume~6 of {\em Acta Numer.}, pages
  229--269. Cambridge Univ. Press, Cambridge, 1997.

\bibitem{tausch-white03}
Johannes Tausch and Jacob White.
\newblock Multiscale bases for the sparse representation of boundary integral
  operators on complex geometry.
\newblock {\em SIAM J. Sci. Comput.}, 24(5):1610--1629, 2003.

\bibitem{hackbusch15}
Wolfgang Hackbusch.
\newblock {\em Hierarchical matrices: algorithms and analysis}, volume~49 of
  {\em Springer Series in Computational Mathematics}.
\newblock Springer, Heidelberg, 2015.

\bibitem{mclean18}
William McLean.
\newblock Exponential sum approximations for {$t^{-\beta}$}.
\newblock In {\em Contemporary computational mathematics---a celebration of the
  80th birthday of {I}an {S}loan. {V}ol. 1, 2}, pages 911--930. Springer, Cham,
  2018.

\bibitem{schaedle-lopez-fernandez-lubich06}
Achim Sch\"{a}dle, Mar\'{\i}a L\'{o}pez-Fern\'{a}ndez, and Christian Lubich.
\newblock Fast and oblivious convolution quadrature.
\newblock {\em SIAM J. Sci. Comput.}, 28(2):421--438, 2006.

\bibitem{doelz-egger-shashkov21}
J.~D\"{o}lz, H.~Egger, and V.~Shashkov.
\newblock A fast and oblivious matrix compression algorithm for {V}olterra
  integral operators.
\newblock {\em Adv. Comput. Math.}, 47(6):Paper No. 81, 24, 2021.

\bibitem{khristenko-wohlmuth23}
Ustim Khristenko and Barbara Wohlmuth.
\newblock Solving time-fractional differential equations via rational
  approximation.
\newblock {\em IMA J. Numer. Anal.}, 43(3):1263--1290, 2023.

\bibitem{shen-handbook}
Jie Shen and Changtao Sheng.
\newblock Spectral methods for fractional differential equations using
  generalized {J}acobi functions.
\newblock In {\em Handbook of fractional calculus with applications. {V}ol. 3},
  pages 127--155. De Gruyter, Berlin, 2019.

\bibitem{lischke-handbook}
Anna Lischke, Mohsen Zayernouri, and Zhongqiang Zhang.
\newblock Spectral and spectral element methods for fractional
  advection-diffusion-reaction equations.
\newblock In {\em Handbook of fractional calculus with applications. {V}ol. 3},
  pages 157--183. De Gruyter, Berlin, 2019.

\bibitem{zayernouri-wang-shen-karniadakis23}
Mohsen Zayernouri, Li-Liang Wang, Jie Shen, and George~Em Karniadakis.
\newblock {\em Spectral and Spectral-Element Methods for Fractional Ordinary
  and Partial Diﬀerential Equations}.
\newblock Cambridge University Press, 2024.

\bibitem{chen-shen20}
Sheng Chen and Jie Shen.
\newblock Log orthogonal functions: approximation properties and applications.
\newblock {\em IMA J. Numer. Anal.}, 42(1):712--743, 2022.

\bibitem{chen-shen-wang16}
Sheng Chen, Jie Shen, and Li-Lian Wang.
\newblock Generalized {J}acobi functions and their applications to fractional
  differential equations.
\newblock {\em Math. Comp.}, 85(300):1603--1638, 2016.

\bibitem{mao-shen18}
Zhiping Mao and Jie Shen.
\newblock Spectral element method with geometric mesh for two-sided fractional
  differential equations.
\newblock {\em Adv. Comput. Math.}, 44(3):745--771, 2018.

\bibitem{zayernouri-karniadakis13}
Mohsen Zayernouri and George~Em Karniadakis.
\newblock Fractional {S}turm-{L}iouville eigen-problems: theory and numerical
  approximation.
\newblock {\em J. Comput. Phys.}, 252:495--517, 2013.

\bibitem{hou-hasan-xu18}
Dianming Hou, Mohammad~Tanzil Hasan, and Chuanju Xu.
\newblock M\"{u}ntz spectral methods for the time-fractional diffusion
  equation.
\newblock {\em Comput. Methods Appl. Math.}, 18(1):43--62, 2018.

\bibitem{schwab98}
Ch. Schwab.
\newblock {\em {$p$}- and {$hp$}-finite element methods}.
\newblock Numerical Mathematics and Scientific Computation. The Clarendon Press
  Oxford University Press, New York, 1998.
\newblock Theory and applications in solid and fluid mechanics.

\bibitem{li-xu09}
Xianjuan Li and Chuanju Xu.
\newblock A space-time spectral method for the time fractional diffusion
  equation.
\newblock {\em SIAM J. Numer. Anal.}, 47(3):2108--2131, 2009.

\bibitem{ervin-roop06}
Vincent~J. Ervin and John~Paul Roop.
\newblock Variational formulation for the stationary fractional advection
  dispersion equation.
\newblock {\em Numer. Methods Partial Differential Equations}, 22(3):558--576,
  2006.

\bibitem{zayernouri-ainsworth-karniadakis15}
Mohsen Zayernouri, Mark Ainsworth, and George~Em Karniadakis.
\newblock A unified {P}etrov-{G}alerkin spectral method for fractional {PDE}s.
\newblock {\em Comput. Methods Appl. Mech. Engrg.}, 283:1545--1569, 2015.

\bibitem{golub-vanloan13}
Gene~H. Golub and Charles~F. Van~Loan.
\newblock {\em Matrix computations}.
\newblock Johns Hopkins Studies in the Mathematical Sciences. Johns Hopkins
  University Press, Baltimore, MD, fourth edition, 2013.

\bibitem{bonito-etal-survey18}
Andrea Bonito, Juan~Pablo Borthagaray, Ricardo~H. Nochetto, Enrique
  Ot\'{a}rola, and Abner~J. Salgado.
\newblock Numerical methods for fractional diffusion.
\newblock {\em Comput. Vis. Sci.}, 19(5-6):19--46, 2018.

\bibitem{lischke-etal20}
Anna Lischke, Guofei Pang, Mamikon Gulian, and et~al.
\newblock What is the fractional {L}aplacian? {A} comparative review with new
  results.
\newblock {\em J. Comput. Phys.}, 404:109009, 62, 2020.

\bibitem{hofreither20}
Clemens Hofreither.
\newblock A unified view of some numerical methods for fractional diffusion.
\newblock {\em Comput. Math. Appl.}, 80(2):332--350, 2020.

\bibitem{burrage-hale-kay12}
Kevin Burrage, Nicholas Hale, and David Kay.
\newblock An efficient implicit {FEM} scheme for fractional-in-space
  reaction-diffusion equations.
\newblock {\em SIAM J. Sci. Comput.}, 34(4):A2145--A2172, 2012.

\bibitem{hofreither21}
Clemens Hofreither.
\newblock An algorithm for best rational approximation based on barycentric
  rational interpolation.
\newblock {\em Numer. Algorithms}, 88(1):365--388, 2021.

\bibitem{hackbusch19}
Wolfgang Hackbusch.
\newblock Computation of best {$L^\infty$} exponential sums for {$1/x$} by
  {R}emez' algorithm.
\newblock {\em Comput. Vis. Sci.}, 20(1-2):1--11, 2019.

\bibitem{danczul-schoeberl22}
Tobias Danczul and Joachim Sch\"{o}berl.
\newblock A reduced basis method for fractional diffusion operators {I}.
\newblock {\em Numer. Math.}, 151(2):369--404, 2022.

\bibitem{banjai-etal19}
Lehel Banjai, Jens~M. Melenk, Ricardo~H. Nochetto, Enrique Ot\'{a}rola,
  Abner~J. Salgado, and Christoph Schwab.
\newblock Tensor {FEM} for spectral fractional diffusion.
\newblock {\em Found. Comput. Math.}, 19(4):901--962, 2019.

\bibitem{banjai-etal23}
Lehel Banjai, Jens~M. Melenk, and Christoph Schwab.
\newblock Exponential convergence of {$hp$} {FEM} for spectral fractional
  diffusion in polygons.
\newblock {\em Numer. Math.}, 153(1):1--47, 2023.

\bibitem{cabre-tan10}
Xavier Cabr\'{e} and Jinggang Tan.
\newblock Positive solutions of nonlinear problems involving the square root of
  the {L}aplacian.
\newblock {\em Adv. Math.}, 224(5):2052--2093, 2010.

\bibitem{stinga-torrea10}
Pablo~Ra\'{u}l Stinga and Jos\'{e}~Luis Torrea.
\newblock Extension problem and {H}arnack's inequality for some fractional
  operators.
\newblock {\em Comm. Partial Differential Equations}, 35(11):2092--2122, 2010.

\bibitem{nochetto-otarola-salgado15}
Ricardo~H. Nochetto, Enrique Ot\'{a}rola, and Abner~J. Salgado.
\newblock A {PDE} approach to fractional diffusion in general domains: a priori
  error analysis.
\newblock {\em Found. Comput. Math.}, 15(3):733--791, 2015.

\bibitem{melenk-rieder21}
Jens~Markus Melenk and Alexander Rieder.
\newblock {$hp$}-{FEM} for the fractional heat equation.
\newblock {\em IMA J. Numer. Anal.}, 41(1):412--454, 2021.

\bibitem{melenk-rieder23}
Jens~Markus Melenk and Alexander Rieder.
\newblock An exponentially convergent discretization for space--time fractional
  parabolic equations using {$hp$}-{FEM}.
\newblock {\em IMA J. Numer. Anal.}, 43(4):2352--2376, 2023.

\bibitem{bonito-pasciak15}
Andrea Bonito and Joseph~E. Pasciak.
\newblock Numerical approximation of fractional powers of elliptic operators.
\newblock {\em Math. Comp.}, 84(295):2083--2110, 2015.

\bibitem{bonito-lei-pasciak19}
Andrea Bonito, Wenyu Lei, and Joseph~E. Pasciak.
\newblock On sinc quadrature approximations of fractional powers of regularly
  accretive operators.
\newblock {\em J. Numer. Math.}, 27(2):57--68, 2019.

\bibitem{bonito-lei-pasciak17}
Andrea Bonito, Wenyu Lei, and Joseph~E. Pasciak.
\newblock Numerical approximation of space-time fractional parabolic equations.
\newblock {\em Comput. Methods Appl. Math.}, 17(4):679--705, 2017.

\bibitem{cusimano-del-teso-gerardo-giorda20}
Nicole Cusimano, F\'{e}lix del Teso, and Luca Gerardo-Giorda.
\newblock Numerical approximations for fractional elliptic equations {$VIA$}
  the method of semigroups.
\newblock {\em ESAIM Math. Model. Numer. Anal.}, 54(3):751--774, 2020.

\bibitem{mao-shen17}
Zhiping Mao and Jie Shen.
\newblock Hermite spectral methods for fractional {PDE}s in unbounded domains.
\newblock {\em SIAM J. Sci. Comput.}, 39(5):A1928--A1950, 2017.

\bibitem{tang-yuan-zhou18}
Tao Tang, Huifang Yuan, and Tao Zhou.
\newblock Hermite spectral collocation methods for fractional {PDE}s in
  unbounded domains.
\newblock {\em Commun. Comput. Phys.}, 24(4):1143--1168, 2018.

\bibitem{tang-wang-yuan20}
Tao Tang, Li-Lian Wang, Huifang Yuan, and Tao Zhou.
\newblock Rational spectral methods for {PDE}s involving fractional {L}aplacian
  in unbounded domains.
\newblock {\em SIAM J. Sci. Comput.}, 42(2):A585--A611, 2020.

\bibitem{sheng-shen-tang-wang20}
Changtao Sheng, Jie Shen, Tao Tang, Li-Lian Wang, and Huifang Yuan.
\newblock Fast {F}ourier-like mapped {C}hebyshev spectral-{G}alerkin methods
  for {PDE}s with integral fractional {L}aplacian in unbounded domains.
\newblock {\em SIAM J. Numer. Anal.}, 58(5):2435--2464, 2020.

\bibitem{wang-etal21}
Changtao Sheng, Suna Ma, Huiyuan Li, Li-Lian Wang, and Lueling Jia.
\newblock Nontensorial generalised {H}ermite spectral methods for {PDE}s with
  fractional {L}aplacian and {S}chr\"{o}dinger operators.
\newblock {\em ESAIM Math. Model. Numer. Anal.}, 55(5):2141--2168, 2021.

\bibitem{papadopoulos2023frame}
Ioannis P.~A. Papadopoulos, Timon~S. Gutleb, José~A. Carrillo, and Sheehan
  Olver.
\newblock A frame approach for equations involving the fractional {L}aplacian,
  2023.
\newblock arXiv:2311.12451.

\bibitem{antil-dondl-striet21}
Harbir Antil, Patrick Dondl, and Ludwig Striet.
\newblock Approximation of integral fractional {L}aplacian and fractional
  {PDE}s via sinc-basis.
\newblock {\em SIAM J. Sci. Comput.}, 43(4):A2897--A2922, 2021.

\bibitem{antil-dondl-striet22}
Harbir Antil, Patrick Dondl, and Ludwig Striet.
\newblock Analysis of a sinc-{G}alerkin method for the fractional {L}aplacian,
  2023.
\newblock arXiv:2212.11581.

\bibitem{borthagaray-li-nochetto20}
Juan~Pablo Borthagaray, Wenbo Li, and Ricardo~H. Nochetto.
\newblock Linear and nonlinear fractional elliptic problems.
\newblock In {\em 75 years of mathematics of computation}, volume 754 of {\em
  Contemp. Math.}, pages 69--92. Amer. Math. Soc., [Providence], RI, [2020]
  \copyright 2020.

\bibitem{borthagaray-nochetto22}
Juan~Pablo Borthagaray, Wenbo Li, and Ricardo~H. Nochetto.
\newblock Fractional elliptic problems on {L}ipschitz domains: Regularity and
  approximation, 2022.

\bibitem{faustmann-marcati-melenk-schwab22}
Markus Faustmann, Carlo Marcati, Jens~Markus Melenk, and Christoph Schwab.
\newblock Weighted analytic regularity for the integral fractional {L}aplacian
  in polygons.
\newblock {\em SIAM J. Math. Anal.}, 54(6):6323--6357, 2022.

\bibitem{FMMS-hp}
M.~Faustmann, C.~Marcati, J.M. Melenk, and Ch. Schwab.
\newblock Exponential convergence of $hp$-{FEM} for the integral fractional
  {L}aplacian in polygons.
\newblock {\em arXiv:2209.11468}, 2022.

\bibitem{sauter-schwab11}
Stefan~A. Sauter and Christoph Schwab.
\newblock {\em Boundary element methods}, volume~39 of {\em Springer Series in
  Computational Mathematics}.
\newblock Springer-Verlag, Berlin, 2011.
\newblock Translated and expanded from the 2004 German original.

\bibitem{chernov-von-petersdorff-schwab15}
Alexey Chernov, Tobias von Petersdorff, and Christoph Schwab.
\newblock Quadrature algorithms for high dimensional singular integrands on
  simplices.
\newblock {\em Numer. Algorithms}, 70(4):847--874, 2015.

\bibitem{acosta-bersetche-borthagaray17}
G.~Acosta, F.M. Bersetche, and J.P. Borthagaray.
\newblock A short {FE} implementation for a 2d homogeneous {D}irichlet problem
  of a fractional {L}aplacian.
\newblock {\em Comput. Math. Appl.}, 74(4):784--816, 2017.

\bibitem{ainsworth-glusa17}
Mark Ainsworth and Christian Glusa.
\newblock Aspects of an adaptive finite element method for the fractional
  {L}aplacian: a priori and a posteriori error estimates, efficient
  implementation and multigrid solver.
\newblock {\em Comput. Methods Appl. Mech. Engrg.}, 327:4--35, 2017.

\bibitem{karkulik-melenk19}
Michael Karkulik and Jens~Markus Melenk.
\newblock {$\mathcal H$}-matrix approximability of inverses of discretizations
  of the fractional {L}aplacian.
\newblock {\em Adv. Comput. Math.}, 45(5-6):2893--2919, 2019.

\bibitem{faustmann-melenk-parvizi21}
Markus Faustmann, Jens~Markus Melenk, and Maryam Parvizi.
\newblock On the stability of {S}cott-{Z}hang type operators and application to
  multilevel preconditioning in fractional diffusion.
\newblock {\em ESAIM Math. Model. Numer. Anal.}, 55(2):595--625, 2021.

\bibitem{borthagaray-nochetto-wu-xu23}
Juan Borthagaray, Ricardo Nochetto, Shuonan Wu, and Jinchao Xu.
\newblock Robust {BPX} preconditioner for fractional {L}aplacians on bounded
  {L}ipschitz domains.
\newblock {\em Math. Comp.}, 92(344):2439--2473, 2023.

\bibitem{gimperlein-stocek-urzua-torres21}
Heiko Gimperlein, Jakub Stocek, and Carolina Urz\'{u}a-Torres.
\newblock Optimal operator preconditioning for pseudodifferential boundary
  problems.
\newblock {\em Numer. Math.}, 148(1):1--41, 2021.

\bibitem{faustmann-melenk-praetorius21}
Markus Faustmann, Jens~Markus Melenk, and Dirk Praetorius.
\newblock Quasi-optimal convergence rate for an adaptive method for the
  integral fractional {L}aplacian.
\newblock {\em Math. Comp.}, 90(330):1557--1587, 2021.

\bibitem{faustmann2023fembem}
Markus Faustmann and Alexander Rieder.
\newblock {FEM}-{BEM} coupling in fractional diffusion, 2023.

\bibitem{bonito-lei-pasciak19a}
Andrea Bonito, Wenyu Lei, and Joseph~E. Pasciak.
\newblock Numerical approximation of the integral fractional {L}aplacian.
\newblock {\em Numer. Math.}, 142(2):235--278, 2019.

\bibitem{kyprianou-osojnik-shardlow18}
Andreas~E. Kyprianou, Ana Osojnik, and Tony Shardlow.
\newblock Unbiased `walk-on-spheres' {M}onte {C}arlo methods for the fractional
  {L}aplacian.
\newblock {\em IMA J. Numer. Anal.}, 38(3):1550--1578, 2018.

\bibitem{sheng-su-xu23}
Changtao Sheng, Bihao Su, and Chenglong Xu.
\newblock Efficient {M}onte {C}arlo {M}ethod for {I}ntegral {F}ractional
  {L}aplacian in {M}ultiple {D}imensions.
\newblock {\em SIAM J. Numer. Anal.}, 61(5):2035--2061, 2023.

\bibitem{CKCS}
C.~Kuehn and C.~Soresina.
\newblock Numerical continuation for a fast reaction system and its
  cross-diffusion limit.
\newblock {\em SN Partial Differential Equations and Applications}, 1:7, 2020.

\bibitem{MBCKCS}
M.~Breden, C.~Kuehn, and C.~Soresina.
\newblock On the influence of cross-diffusion in pattern formation.
\newblock {\em Journal of Computational Dynamics}, 8(2), 2021.

\bibitem{uecker2019infinite}
H.~Uecker and H.~de~Witt.
\newblock Infinite time horizon spatially distributed optimal control problems
  with pde2path algorithms and tutorial examples.
\newblock {\em arXiv preprint arXiv:1912.11135}, 2019.

\bibitem{uecker2021pde2path}
H.~Uecker.
\newblock pde2path without finite elements, 2021.

\bibitem{meiners2023differential}
A.~Meiners and H.~Uecker.
\newblock Differential geometric bifurcation problems in pde2path--algorithms
  and tutorial examples.
\newblock {\em arXiv preprint arXiv:2309.03646}, 2023.

\bibitem{dohnal2014pde2path}
T.~Dohnal, J.D.M. Rademacher, H.~Uecker, and D.~Wetzel.
\newblock \texttt{pde2path}-version 2.0: faster {FEM}, multi-parameter
  continuation, nonlinear boundary conditions, and periodic domains --- a short
  manual, 2014.

\bibitem{rademacher2018oopde}
J.D.M. Rademacher and H.~Uecker.
\newblock The \texttt{OOPDE} setting of \texttt{pde2path} --- a tutorial via
  some {A}llen--{C}ahn models.
\newblock available at
  \url{http://www.staff.uni-oldenburg.de/hannes.uecker/pde2path/tuts/actut.pdf},
  2018.

\bibitem{uecker2021numerical}
Hannes Uecker.
\newblock {\em Numerical continuation and bifurcation in Nonlinear PDEs}.
\newblock SIAM, 2021.

\bibitem{OOPDE}
U.~Pr{\"u}fert.
\newblock \texttt{OOPDE} - an object oriented approach to finite elements in
  {MATLAB}.
\newblock Quickstart Guide, available at
  \url{http://www.mathe.tu-freiberg.de/nmo/mitarbeiter/uwe-pruefert/software},
  2014.

\bibitem{balakrishnan1960}
A.V. Balakrishnan.
\newblock Fractional powers of closed operators and the semigroups generated by
  them.
\newblock {\em Pacific Journal of Mathematics}, 10(2):419--437, 1960.

\bibitem{yosida1968functional}
K.~Yosida.
\newblock {\em Functional Analysis}.
\newblock Springer-Verlag Berlin Heidelberg, 6 edition, 1980.

\bibitem{dohr2019fem}
S.~Dohr, C.~Kahle, S.~Rogovs, and P.~Swierczynski.
\newblock A {FEM} for an optimal control problem subject to the fractional
  {L}aplace equation.
\newblock {\em Calcolo}, 56(4):37, 2019.

\bibitem{lund1992sinc}
J.~Lund and K.L. Bowers.
\newblock {\em Sinc Methods for Quadrature and Differential Equations},
  volume~32.
\newblock SIAM, 1992.

\bibitem{VideoFolder}
N.~Ehstand, C.~Kuehn, and C.~Soresina.
\newblock Supplementary material.
\newblock Figures and videos at
  \url{https://www-m8.ma.tum.de/bin/view/Allgemeines/CinziaSoresinaPublic},
  matlab scripts at \url{https://github.com/soresina/fractional_pde2path},
  2020.
\newblock Accessed December 20, 2020.

\bibitem{BURKE2007681}
J.~Burke and E.~Knobloch.
\newblock Snakes and ladders: {L}ocalized states in the {S}wift--{H}ohenberg
  equation.
\newblock {\em Physics Letters A}, 360(6):681--688, 2007.

\bibitem{avitabile2010snake}
D.~Avitabile, D.J.B. Lloyd, J.~Burke, E.~Knobloch, and B.~Sandstede.
\newblock To snake or not to snake in the planar {S}wift--{H}ohenberg equation.
\newblock {\em SIAM Journal on Applied Dynamical Systems}, 9(3):704--733, 2010.

\bibitem{houghton2011swift}
S.M. Houghton and E.~Knobloch.
\newblock Swift--{H}ohenberg equation with broken cubic--quintic nonlinearity.
\newblock {\em Physical Review E}, 84(1):016204, 2011.

\bibitem{bergeon2008eckhaus}
A.~Bergeon, J.~Burke, E.~Knobloch, and I.~Mercader.
\newblock Eckhaus instability and homoclinic snaking.
\newblock {\em Physical Review E}, 78(4):046201, 2008.

\end{thebibliography}

\end{document}